\documentclass[a4paper,10pt,reqno]{amsart}

\usepackage{amssymb,amsthm,amsmath,mathrsfs}
\usepackage{mathtools}
\usepackage{bbm}
\usepackage{stackengine}
\usepackage{enumitem}
\usepackage{soul}
\usepackage[foot]{amsaddr}

\usepackage{caption,color,graphicx}
\usepackage[dvipsnames]{xcolor}
\usepackage{soul}

\usepackage{hyperref}
\usepackage{doi}
\hypersetup{colorlinks  = true,
            urlcolor    = blue,
            linkcolor   = red,
            citecolor   = blue}

\usepackage{verbatim}

\usepackage{environ}
\mathtoolsset{showonlyrefs,showmanualtags}

\usepackage{graphicx}

\usepackage[numbers]{natbib}

\usepackage[centering, top=3cm, bottom=2.5cm]{geometry}

\usepackage{tikz}
\usetikzlibrary{arrows.meta}
\newcommand{\N}{\mathbb{N}}
\newcommand{\Z}{\mathbb{Z}}
\newcommand{\Q}{\mathbb{Q}}
\newcommand{\R}{\mathbb{R}}

\newcommand{\E}{\mathbb{E}}
\renewcommand{\P}{\mathbb{P}}
\newcommand{\1}{\mathbbm{1}}

\newcommand{\supp}{\operatorname{supp}\,}
\newcommand{\id}{\operatorname{id}}

\newcommand{\calP}{\mathcal{P}}

\newcommand{\scrP}{\mathscr P}

\newcommand{\dist}{\operatorname{dist}}

\newcommand{\vertiii}[1]{{\left\vert\kern-0.25ex\left\vert\kern-0.25ex\left\vert #1 \right\vert\kern-0.25ex\right\vert\kern-0.25ex\right\vert}}

\let\oldtocsection=\tocsection
\let\oldtocsubsection=\tocsubsection
\renewcommand{\tocsection}[2]{\hspace{0em}\oldtocsection{#1}{#2}}
\renewcommand{\tocsubsection}[2]{\hspace{1em}\oldtocsubsection{#1}{#2}}

\theoremstyle{plain}
\newtheorem*{theorem*}{Theorem}
\newtheorem{theorem}{Theorem}[section]
\newtheorem{lemma}[theorem]{Lemma}
\newtheorem{proposition}[theorem]{Proposition}
\newtheorem{corollary}[theorem]{Corollary}

\newtheorem{maintheorem}{Theorem}

\theoremstyle{definition}
\newtheorem{definition}[theorem]{Definition}
\newtheorem{remark}[theorem]{Remark}

\newcommand{\stepqedsymbol}{{$\blacksquare$}}
\newenvironment{stepproof}[1][\proofname]{%
  \renewcommand{\qedsymbol}{\stepqedsymbol}%
  \begin{proof}[#1]%
}{%
  \end{proof}%
}

\newcommand{\stepsuppressqed}{\renewcommand{\qedsymbol}{}}

\newtheorem{manualstep}{\bf Step}
\newenvironment{step}[1][]{%
    \renewcommand\themanualstep{#1}%
    \begin{manualstep}\itshape
}{\end{manualstep}}

\newtheorem{manualhyp}{Hypothesis}

\IfPackageLoadedTF{MnSymbol}{}{}

\newcommand{\norm}[1]{\lVert #1\rVert}
\newcommand{\abs}[1]{\lvert #1\rvert}
\newcommand{\htop}{h_{\mathrm{top}}}
\newcommand{\Sig}{\Sigma}

\newcommand{\calA}{\mathcal{A}}
\newcommand{\calF}{\mathcal{F}}
\newcommand{\calM}{\mathcal{M}}
\newcommand{\calN}{\mathcal{N}}
\newcommand{\calW}{\mathcal{W}}

\newcommand{\Lip}{\operatorname{Lip}}
\newcommand{\Law}{\operatorname{Law}}

\theoremstyle{plain}

\author[Matheus M.~Castro]{Matheus M. Castro$^{1}$}
\email{\href{mailto:mmcastro@ime.unicamp.br}{mmcastro@ime.unicamp.br}}

\address{$^{1}$Departamento de Matemática – Universidade Estadual de Campinas (UNICAMP), 13083-859, Campinas, SP, Brazil.}

\renewcommand{\P}{\mathbb{P}}

  \newcommand{\calB}{\mathcal{B}}
  \newcommand{\calC}{\mathcal{C}}

  \newcommand{\calH}{\mathcal{H}}
  \newcommand{\scrQ}{\mathscr{Q}}
  \newcommand{\calX}{\mathcal{X}}

  \newcommand{\aux}{\mathrm{aux}}      % auxiliary coordinates
\title[Exponentially mixing equilibrium states without the CLT]{Finite-order exponential mixing does not imply the central limit theorem for unique equilibrium states}

\date{\today}

\begin{document}

\begin{abstract}
Let $M$ be a compact smooth manifold of dimension at least three, $\varphi:M\to\R$ be continuous, and  $R\ge2$ be an integer. We construct a $\mathcal C^\infty$ diffeomorphism $T:M\to M$ and an observable $f\in\mathcal C^\infty(M,\R)$ such that $(T,\varphi)$ has a unique equilibrium state $\mu$ and $\mu(f)=0$. The system $(T,\mu)$ is exponentially mixing of order $R$ but not of order $R+1$ for H\"older observables, and $n^{-1/2}\sum_{i=0}^{n-1}f\circ T^i$ does not satisfy the central limit theorem. The main results have been formalised in Lean 4.
\end{abstract}


\vspace*{-1cm}
\maketitle
\setcounter{tocdepth}{1}
{\thispagestyle{empty}\tableofcontents}

\section{Introduction}
\label{sec:intro}

Let $T$ be a diffeomorphism of a compact Riemannian manifold $M$ preserving a Borel probability measure $\mu$, and let $f:M\to\R$ be an observable with $\mu(f)=0$. The central limit theorem (CLT) asks whether the normalised Birkhoff sums
\begin{equation}
\label{eq:clt}
\frac{S_nf}{\sqrt n}\quad\text{where}\quad S_nf:=\sum_{i=0}^{n-1}f\circ T^i,
\end{equation}
converge in distribution to a centred Gaussian law $\calN(0,\sigma_f^2)$, the degenerate case $\sigma_f^2=0$ included. 

In hyperbolic dynamics, there are many examples of systems for which H\"older observables satisfy the CLT and have \emph{exponential decay of correlations} (see \cite{Bowen1975,ParryPollicott1990,Young1998,BuzziCrovisierSarig2025} and the references therein), that is, for every $0<\beta\le1$ there are $C,c>0$ such that
\begin{equation}
\label{eq:decay}
\abs{\mu(g\,h\circ T^n)-\mu(g)\mu(h)}\le Ce^{-cn}\norm g_{\mathcal C^\beta}\norm h_{\mathcal C^\beta}
\end{equation}
for all $g,h\in\mathcal C^\beta(M)$ and all integers $n\ge1$. It is therefore natural to ask whether \eqref{eq:decay} alone implies the CLT \cite{DolgopyatDongKanigowskiNandori2022}. Bj\"orklund and Gorodnik \cite[Theorem 1.1]{BjorklundGorodnik2020} proved that every centred H\"older observable satisfies the CLT if $(T,\mu)$ is \emph{exponentially mixing of all orders}, that is, for every $s\ge2$ and $0<\beta\le1$ there are $C_s,c_s>0$ such that
\begin{equation}
\label{eq:allorders}
\Bigl|\mu\bigl(f_1\cdot f_2\circ T^{n_1}\cdots f_s\circ T^{n_1+\dots+n_{s-1}}\bigr)-\prod_{j=1}^s\mu(f_j)\Bigr|
\le C_se^{-c_s\min_{1\le k<s}n_k}\prod_{j=1}^s\norm{f_j}_{\mathcal C^\beta}
\end{equation}
for all $f_1,\dots,f_s\in \mathcal C^\beta(M)$ and all integers $n_1,\dots,n_{s-1}\ge1$. A system satisfying \eqref{eq:allorders} for a fixed $s$ is called \emph{exponentially mixing of order $s$}, so \eqref{eq:decay} is exponential mixing of order $2$. 

Recently, Bassols Cornudella and Castro \cite{BassolsCastro2026} showed that the theorem of Bj\"orklund and Gorodnik is sharp. For every $R\ge2$ they constructed a $\mathcal C^\infty$ surface diffeomorphism with a physical measure which is exponentially mixing of order $R$ and for which the CLT fails.

A distinguished class of invariant measures is given by the equilibrium states, which come from statistical mechanics \cite{Sinai1972,Ruelle1978,Bowen1975}. For a continuous potential $\varphi:M\to\R$, an \emph{equilibrium state} of $(T,\varphi)$ is a $T$-invariant Borel probability $\mu$ attaining the supremum in the variational principle
$$
P_T(\varphi)=\sup\Bigl\{h_\rho(T)+\int\varphi\,d\rho\ ;\ \rho\text{ is a $T$-invariant Borel probability}\Bigr\},
$$
where $P_T(\varphi)$ is the topological pressure \cite[Chapter 9]{Walters1982} (see Section \ref{sec:entropy-def} for details).
Measures of maximal entropy are the equilibrium states of $\varphi\equiv0$, and on hyperbolic attractors SRB measures are the equilibrium states of the geometric potential. In \cite[Question 1.3]{BassolsCastro2026} it is asked whether the same can happen for equilibrium states of H\"older potentials. In this paper we show that it happens for \emph{every} prescribed continuous potential, in particular for every H\"older potential, and for a \emph{unique} equilibrium state.

\begin{maintheorem}
\label{thm:A}
Let $M$ be a compact smooth manifold of dimension at least three, let $\varphi:M\to\R$ be continuous, and let $R\ge2$ be an integer. There exist a $\mathcal C^\infty$ diffeomorphism $T:M\to M$ and $f\in \mathcal C^\infty(M,\R)$ such that $(T,\varphi)$ has a unique equilibrium state $\mu$, $\mu(f)=0$, and the following hold.
\begin{enumerate}[label=\textup{(\alph*)},leftmargin=*]
\item\label{A:mix} $(T,\mu)$ is exponentially mixing of order $R$ for H\"older observables.
\item\label{A:nomix} $(T,\mu)$ is not exponentially mixing of order $R+1$ for H\"older observables.
\item\label{A:clt} The laws of $S_nf/\sqrt n$ under $\mu$ have no Gaussian weak subsequential limit. In particular $f$ does not satisfy the CLT.
\end{enumerate}
\end{maintheorem}
\begin{remark}
For every $K>0$, a single pair $(T,\mu)$ satisfies Theorem~\ref{thm:A} for all continuous $\varphi$ with $\sup_M\varphi-\inf_M\varphi<K$, since the potential enters the proof only through the choice of the integer $q$ in Section~\ref{sec:realisation}. By Lemma~\ref{lem:holderC1}, H\"older observables can be replaced by $\mathcal C^r$ observables, for any integer $r\ge1$.
\end{remark}
Theorem~\ref{thm:A} is proved in Section~\ref{sec:realisation}, as a consequence of the more precise Theorem~\ref{thm:Aprecise}.

The maps and equilibrium states constructed in Theorem~\ref{thm:A} reveal a gap in the literature. They have exponential decay of correlations, but none of the methods used to study equilibrium states in smooth dynamics applies to them. To the best of our knowledge, every proof of exponential decay of correlations for equilibrium states of H\"older potentials relies on a spectral gap of a transfer operator, either directly \cite{Bowen1975,ParryPollicott1990} or through a symbolic or tower model \cite{Young1998,Sarig2013,BuzziCrovisierSarig2025,LuoYang2026}, except for SRB measures, for which geometric methods are also available \cite{Dolgopyat2004}. A spectral gap also yields the CLT by the method of Nagaev and Guivarc'h \cite{Gouezel2015}, which fails here.  In our construction, $T$ is a smooth realisation of a subshift, and we exploit an independence property of the symbolic measure which gives exponential mixing of order $R$ but not of all orders.

Taking $\varphi\equiv0$, whose equilibrium states are the measures of maximal entropy, Theorem~\ref{thm:A} gives the following.

\begin{corollary}
\label{cor:mme}
Let $M$ be a compact smooth manifold of dimension at least three and let $R\ge2$ be an integer. There exist a $\mathcal C^\infty$ diffeomorphism $T:M\to M$ and $f\in\mathcal C^\infty(M,\R)$ such that $T$ has a unique measure of maximal entropy $\mu$, $\mu(f)=0$, and the following hold.
\begin{enumerate}[label=\textup{(\alph*)},leftmargin=*]
\item $(T,\mu)$ is exponentially mixing of order $R$ for H\"older observables.
\item $(T,\mu)$ is not exponentially mixing of order $R+1$ for H\"older observables.
\item The laws of $S_nf/\sqrt n$ under $\mu$ have no Gaussian weak subsequential limit.
\end{enumerate}
\end{corollary}

Before describing the proof we place the result in its context.

\subsection{Central limit theorems and mixing in smooth ergodic theory}
\label{sec:lit-clt}
 
The CLT for deterministic systems goes back to Sinai \cite{Sinai1960} for geodesic flows on surfaces of constant negative curvature and to Ratner \cite{Ratner1973} for geodesic flows in variable negative curvature. Three main methods emerged, namely the spectral method of Nagaev \cite{GuivarchHardy1988,Gouezel2015}, the martingale method of Gordin \cite{Gordin1969,Liverani1996}, and Young towers \cite{Young1998}, which cover many nonuniformly hyperbolic systems. These methods also yield further limit theorems, such as almost sure invariance principles, local limit theorems and large deviations (see \cite{DenkerPhilipp1984,Gouezel2015} and the references therein). For partially hyperbolic systems the CLT was obtained from the equidistribution of unstable leaves \cite{Dolgopyat2004}, and for toral and nilmanifold automorphisms in \cite{LeBorgne1999} and, through exponential mixing of all orders, in \cite{GorodnikSpatzier2014,GorodnikSpatzier2015}. In each of these settings the proof produces more than \eqref{eq:decay}, namely a spectral gap, a martingale approximation or multiple mixing estimates. Bj\"orklund and Gorodnik \cite{BjorklundGorodnik2020} isolated exponential mixing of all orders as a sufficient condition through the method of moments. Whether exponential mixing implies the CLT was raised for smooth invariant measures in \cite[Appendix C]{DolgopyatDongKanigowskiNandori2022}.
 
\subsection{Counterexamples to the CLT with prescribed independence}
\label{sec:lit-counter}
 
How much independence is compatible with the failure of the CLT is a classical theme in probability. Rosenblatt \cite{Rosenblatt1956} and Ibragimov \cite{Ibragimov1962} proved CLTs under strong mixing with rates. Without rates the CLT may fail. Janson \cite{Janson1988} constructed stationary pairwise independent counterexamples, Bradley constructed absolutely regular ones \cite{Bradley1989} and showed that his example is triplewise independent \cite{Bradley2007}, Pruss \cite{Pruss1998} constructed bounded $N$-tuplewise independent identically distributed counterexamples, and Bradley and Pruss \cite{BradleyPruss2009} made them strictly stationary. Weakley \cite{Weakley2013} obtained strictly stationary $N$-tuplewise independent counterexamples which are mixing, and Bradley \cite{Bradley2010} constructed a strictly stationary, causal, five-tuplewise independent Bernoulli counterexample whose normalised sums have no Gaussian weak subsequential limit. 
 
Bassols Cornudella and Castro \cite[Theorem C]{BassolsCastro2026} added quantitative mixing to this list, with exponential mixing of order $R$, failure at order $2R+2$, and no Gaussian weak subsequential limit. Their Theorem A realises the process on a fixed horseshoe of a linear Anosov diffeomorphism of the two-torus, and their Theorem B realises it as a physical measure of a $\mathcal C^\infty$ surface diffeomorphism. We use their Theorem C as a black box in this paper.
 
\subsection{Equilibrium states and measures of maximal entropy}
\label{sec:lit-eq}
 
The thermodynamic formalism of Sinai, Ruelle and Bowen \cite{Sinai1972,Ruelle1978,Bowen1975} shows that for a H\"older potential on a basic set of an Axiom A diffeomorphism the equilibrium state is unique, has the Gibbs property, has exponential decay of correlations for H\"older observables when the basic set is topologically mixing, and satisfies the CLT (see also Parry and Pollicott \cite{ParryPollicott1990}). For subshifts of finite type the measure of maximal entropy is the Parry measure \cite{Parry1964}, and uniqueness of equilibrium states extends to expansive systems with specification by Bowen \cite{Bowen1974}.

Beyond uniform hyperbolicity, equilibrium states have been studied through a number of methods, including symbolic codings by countable Markov shifts \cite{Sarig1999,Sarig2013,BenOvadia2018}, Young towers \cite{Young1998,PesinSentiZhang2016}, non-uniform specification \cite{ClimenhagaThompson2016,ClimenhagaFisherThompson2018}, and geometric measure theory \cite{ClimenhagaPesinZelerowicz2019,ClimenhagaPesinZelerowicz2020}.

Existence of measures of maximal entropy holds for all $\mathcal C^\infty$ diffeomorphisms by Newhouse \cite{Newhouse1989}, using Yomdin's theory \cite{Yomdin1987}, and fails in finite smoothness \cite{Misiurewicz1973,Buzzi2014}. For $\mathcal C^\infty$ surface diffeomorphisms of positive entropy, Buzzi, Crovisier and Sarig proved that there are finitely many ergodic measures of maximal entropy \cite{BuzziCrovisierSarig2022} and that they are strongly positively recurrent, hence exponentially mixing up to a period and satisfying the CLT \cite{BuzziCrovisierSarig2025}. Luo and Yang \cite{LuoYang2026} obtained exponential mixing for equilibrium states of H\"older potentials $\varphi$ with $\sup\varphi-\inf\varphi$ smaller than the topological entropy.

Theorem~\ref{thm:A} shows that in every dimension at least three exponential mixing of a prescribed finite order does not suffice, even for a unique equilibrium state. The mechanism is a zero Lyapunov exponent, which cannot occur for an ergodic measure of positive entropy of a surface diffeomorphism.

\subsection{Strategy}
\label{sec:strategy}

Here we outline the proof of Theorem~\ref{thm:A}. Let $\Sig:=\{-1,1\}^\Z$ with the left shift $\sigma$, let $Z_i:z\in\Sig\mapsto z_i\in\{-1,1\}$ be the coordinate projections, where $i\in\Z$, and let $\nu$ be the $\sigma$-invariant probability on $\Sig$ given by \cite[Theorem C]{BassolsCastro2026} (see Theorem~\ref{thm:BCC}). Under $\nu$ the process $Z=(Z_i)_{i\in\Z}$ has the following properties:
\begin{itemize}[leftmargin=*]
\item The variables $Z_i$ take the values $-1$ and $1$ with probability $\frac12$ each, and any $2R+1$ of them are independent.
\item $(\Sig,\sigma,\nu)$ is exponentially mixing of order $R$, but not of order $R+1$, for H\"older observables (see Lemma~\ref{lem:nomixR1}).
\item The normalised sums $S_n/\sqrt n$, where $S_n:=Z_0+\dots+Z_{n-1}$, have no Gaussian weak subsequential limit.
\item The measure $\nu$ has full support on $\Sigma$.
\end{itemize}

Since $\nu$ has full support, its support is the full shift $\Sig$, whose unique measure of maximal entropy is the product measure $(\frac12\delta_{-1}+\frac12\delta_1)^{\otimes\Z}$ (see \cite[Section 8.3]{Walters1982}). Under the measure of maximal entropy, the coordinates $Z_i$ are independent, so by the classical CLT the normalised sums $S_n/\sqrt n$ converge weakly to the standard Gaussian law, while under $\nu$ they have no Gaussian weak subsequential limit. Hence $\nu$ is not a measure of maximal entropy of its support.

We therefore realise $\nu$ as the marginal of a measure of maximal entropy on a larger subshift. We add auxiliary coordinates $\omega_i\in\{-1,1\}^q$ and construct a subshift $Y$ of $(\{-1,1\}\times\{-1,1\}^q)^\Z$ with an invariant measure $\mu_0$ whose marginal on the first coordinate is $\nu$. Under $\mu_0$ the auxiliary coordinates carry almost all of their possible entropy, while every other ergodic measure of $(Y,\sigma)$ has much smaller entropy. Explicitly, 
$$h_{\mu_0}(\sigma)\ge \left(\frac78q-1\right)\log2,\text{ and }h_\rho(\sigma)\le \left(2+\frac q8\right)\log2$$ for every ergodic $\rho\ne\mu_0$ on $Y$ (Propositions~\ref{prop:lower} and~\ref{prop:upper}). The construction also keeps exact independence of $R$ distant blocks of coordinates, which gives exponential mixing of order $R$. This gives the following symbolic theorem.

\begin{theorem*}[Symbolic theorem, informal]
For every integer $q\ge12$ there are a subshift $Y\subset(\{-1,1\}\times\{-1,1\}^q)^\Z$, with first-coordinate factor map $\pi_1:Y\to\{-1,1\}^\Z$, and an invariant measure $\mu_0$ on $Y$ with $(\pi_1)_*\mu_0=\nu$, such that $(\sigma,\mu_0)$ is ergodic and exponentially mixing of order $R$ for every H\"older class, and $\mu_0$ is the unique equilibrium state of $\psi|_Y$ for every continuous potential $\psi$ with $\sup_Y\psi-\inf_Y\psi<\frac q2\log2$. The precise statement is Theorem~\ref{thm:Bformal}.
\end{theorem*}

Theorem~\ref{thm:A} follows from the symbolic theorem by a smooth realisation. We realise $(Y,\mu_0)$ smoothly on an affine horseshoe of a disc diffeomorphism, and we embed the disc in a solid torus by the neutral isolation construction of Lin and Tian \cite{LinTian2026}, so that every other ergodic measure of the resulting diffeomorphism is a Dirac mass at a fixed point. The coding is Lipschitz, so exponential mixing of order $R$ passes to H\"older observables, and a smooth function which records the first symbolic coordinate has the law $\nu$. Choosing $q$ with $\frac q2\log2>\sup_M\varphi-\inf_M\varphi$ makes the equilibrium state unique. The construction of $Y$ and the proof of the symbolic theorem occupy Sections~\ref{sec:BCC}--\ref{sec:entropy}, with an outline at the beginning of Section~\ref{sec:extension}, and the smooth realisation is in Section~\ref{sec:realisation}.

\subsection{Organisation}
Section~\ref{sec:prelim} fixes the notation and recalls entropy, pressure, equilibrium states and exponential mixing of order $s$. Section~\ref{sec:BCC} recalls \cite[Theorem C]{BassolsCastro2026} and the nested marker sets of \cite[Proposition 4.2]{BassolsCastro2026}, and introduces the auxiliary shift $\calX$. Section~\ref{sec:products} constructs laws on products of finite sets with prescribed marginals. Section~\ref{sec:extension} outlines the construction of $Y$, constructs the laws of the blocks, the subshift $Y$ and the measure $\mu_0$, and proves that every invariant probability on $\calX$ has exactly one invariant lift and that $Y$ has no fixed point. Section~\ref{sec:mixing} proves that $(\sigma,\mu_0)$ is exponentially mixing of order $R$. Section~\ref{sec:entropy} proves the entropy gap and the symbolic theorem. Section~\ref{sec:realisation} contains the smooth realisation and the proof of Theorem~\ref{thm:A}. The main results of this article have been formalised in Lean 4.35.0, as explained in Section~\ref{sec:lean}.

\section{Notation, pressure and exponential mixing}
\label{sec:prelim}
In this section we fix the notation, recall entropy, pressure and equilibrium states, and recall exponential mixing of order $s$.

\subsection{Notation}
\label{sec:notation}
Throughout this paper $R\ge2$ is the integer of Theorem~\ref{thm:A}.

\begin{enumerate}[label=\textup{(N\arabic*)},leftmargin=*]
\item\label{not:intervals} \emph{Intervals.} An \emph{interval} is a finite nonempty set of consecutive integers $I=[a,b]:=\{a,a+1,\dots,b\}$, of length $|I|=b-a+1$. For intervals $I,J$ we write $$\dist(I,J):=\min\{|i-j|:i\in I,\ j\in J\},$$ and for an interval $I$ and $c\ge0$ we write $I^{+c}:=\{i\in\Z\ ;\ |i-i'|\le c\text{ for some }i'\in I\}$ and $I^{(c)}:=\{i\in I\ ;\ [i-c,i+c]\subset I\}$, so that
\begin{equation}
\label{eq:interior}
\#\bigl(I\setminus I^{(c)}\bigr)\le2c .
\end{equation}

\item\label{not:restrictions} \emph{Restrictions.} For a sequence $x=(x_n)_{n\in\Z}$ and a set $F\subset\Z$, $x|_F:=(x_n)_{n\in F}$, and for an interval $I$ we also write $x[I]:=x|_I$. For a finite set $F\subset\Z$ we identify a function $H:\calA^F\to\R$ with the function $x\mapsto H(x|_F)$ on $\calA^\Z$, and similarly for other alphabets.

\item \emph{Shifts.} For a finite alphabet $\calA$ the full shift $\calA^\Z$ carries the product topology and the left shift $\sigma$, $(\sigma x)_n:=x_{n+1}$. A \emph{subshift} is a closed $\sigma$-invariant subset. For $0<\theta<1$ the symbolic metric is $d_\theta(x,y):=\theta^{\min\{|n|:x_n\ne y_n\}}$ for $x\ne y$. We also write $\sigma$ for the left shift on the other spaces of sequences indexed by $\Z$ in this paper, and for products of such shifts, except on the spaces $\Sigma_\calH$ and $\calX$ of Section~\ref{sec:BCC}, where we write $\sigma_\calH$ and $\sigma_\calX$.

\item \emph{H\"older functions on subshifts.}  Let  $F:X\to\R$ be a function on a subshift $X$. We define
$$\Lip_\theta(F):=\sup_{x\ne y} \frac{|F(x)-F(y)|}{d_\theta(x,y)}\ \text{and }\norm F_\theta:=\norm F_\infty+\Lip_\theta(F).$$ The space $\calF_\theta(X)$ of functions with $\norm F_\theta<\infty$ consists of the $d_\theta$-Lipschitz functions. A function is $\gamma$-H\"older for $d_\theta$ exactly when it is Lipschitz for $d_{\theta^\gamma}$, and we call the spaces $\calF_\theta$ the \emph{H\"older classes}. If $F\in\calF_\theta(X)$ and $F_m$ is any function depending only on $x[-m,m]$ with $F_m(x)=F(x')$ for some $x'$ with $x'[-m,m]=x[-m,m]$, then
\begin{equation}
\label{eq:cylinderapprox}
\norm{F-F_m}_\infty\le\Lip_\theta(F)\,\theta^{m+1}.
\end{equation}

\item \emph{Hölder functions on manifolds.} For a compact metric space $X$ we write $\mathcal C^0(X)$ for the space of continuous functions $X\to\R$ with the norm $\norm\cdot_\infty$. On a compact Riemannian manifold $M$ we write $\mathcal C^\beta(M)$, $0<\beta\le1$, for the space of $\beta$-H\"older functions with the norm $$\norm\phi_{\mathcal C^\beta}:=\norm\phi_\infty+\sup_{p\ne p'}\frac{|\phi(p)-\phi(p')|}{d_M(p,p')^\beta}.$$ For an integer $r\ge1$ we write $\mathcal C^r(M)$ for the space of $r$ times continuously differentiable functions with the norm $\norm\phi_{\mathcal C^r}:=\sum_{j=0}^r\norm{\nabla^j\phi}_\infty$, where $\nabla^j\phi$ is the $j$-th covariant derivative of $\phi$. In particular, $\mathcal C^1(M)$ always denotes the space of continuously differentiable functions.

\item\label{not:dep} \emph{Processes and finite dependence.} For a stationary process $W=(W_n)_{n\in\Z}$ with values in a finite set $\calA$, we identify $\Law(W)$ with a shift-invariant probability on $\calA^\Z$. For $d\ge0$ a process $(W_n)_{n\in\Z}$ is \emph{$d$-dependent} if its restrictions to any two sets of indices at distance greater than $d$ are independent (see \cite[Definition 4.1]{BassolsCastro2026}). Applying this to the sets $F_1\cup\dots\cup F_{s-1}$ and $F_s$ for $s=2,\dots,r$ shows that the restrictions to sets $F_1,\dots,F_r$ whose pairwise distances are greater than $d$ are jointly independent.

\item \emph{Rademacher processes.} A process $W=(W_n)_{n\in\Z}$ is a \emph{Rademacher process} if $\P(W_n=1)=\P(W_n=-1)=\frac12$ for every $n$. We do not assume that the $W_n$ are independent.
\end{enumerate}

\subsection{Entropy, pressure and equilibrium states}
\label{sec:entropy-def}
Let $T$ be a homeomorphism of a compact metric space $X$. For a $T$-invariant Borel probability $\mu$ we denote by $h_\mu(T)$ its Kolmogorov--Sinai entropy, and for a continuous $\psi:X\to\R$ we denote by $P_T(\psi)$ the topological pressure (see \cite[Chapters 4 and 9]{Walters1982}). By the variational principle \cite[Chapter 9]{Walters1982},
\begin{equation}
\label{eq:varprinciple}
P_T(\psi)=\sup\Bigl\{h_\mu(T)+\int\psi\,d\mu\ \Big|\ \mu\text{ a $T$-invariant Borel probability on }X\Bigr\}.
\end{equation}
An \emph{equilibrium state} of $(T,\psi)$ is a maximiser in \eqref{eq:varprinciple}. For $\psi\equiv0$ the pressure is the topological entropy $\htop(T)$, and equilibrium states are measures of maximal entropy. The map $\mu\mapsto h_\mu(T)$ is affine and respects the ergodic decomposition \cite[Chapter 8]{Walters1982}. Hence almost every ergodic component of an equilibrium state is an equilibrium state. A $T$-invariant probability is ergodic exactly when it is an extreme point of the convex set of $T$-invariant probabilities \cite[Theorem 6.10]{Walters1982}.

\emph{Shannon entropy.} Let $\mu$ be a probability, let $\xi,\eta$ be finite measurable partitions and let $\mathcal G\subset\mathcal G'$ be $\sigma$-algebras. The entropy of $\xi$ and its conditional entropy given $\mathcal G$ are
$$
H_\mu(\xi):=-\sum_{C\in\xi}\mu(C)\log\mu(C)\quad\text{and}\quad H_\mu(\xi\mid\mathcal G):=-\int\sum_{C\in\xi}\mu(C\mid\mathcal G)\log\mu(C\mid\mathcal G)\,d\mu,
$$
where $\mu(C\mid\mathcal G):=\E_\mu(\mathbbm 1_C\mid\mathcal G)$ and $0\log0=0$. Throughout the paper we use the following well-known bounds
\begin{equation}
\label{eq:shannon}
H_\mu(\xi\mid\mathcal G')\le H_\mu(\xi\mid\mathcal G)\le H_\mu(\xi)\le\log\#\xi\ \text{and}\ H_\mu(\xi\vee\eta\mid\mathcal G)=H_\mu(\xi\mid\mathcal G)+H_\mu(\eta\mid\xi\vee\mathcal G),
\end{equation}
where $\xi\vee\mathcal G$ is the $\sigma$-algebra generated by $\xi$ and $\mathcal G$ \cite[Chapter 4]{Walters1982}. If $\mathcal G_1\subset\mathcal G_2\subset\cdots$ are $\sigma$-algebras generating a $\sigma$-algebra $\mathcal G_\infty$, then $H_\mu(\xi\mid\mathcal G_s)$ decreases to $H_\mu(\xi\mid\mathcal G_\infty)$ as $s\to\infty$, by the martingale convergence theorem (see \cite[Chapter 1]{Downarowicz2011}). Let $D\subset\Z$, let $A_i$ be a finite set for every $i\in D$, and let $X\subset\prod_{i\in D}A_i$. For a finite set $F\subset D$, let $\scrP_F$ be the partition of $X$ into the cylinders $[w]_F$, $w\in\prod_{i\in F}A_i$, where $[w]_F$ is the set of $x\in X$ with $x|_F=w$. We use these partitions for subshifts $X\subset\calA^\Z$, where $D=\Z$ and $A_i=\calA$, and for subsets $X$ of finite products, where $D$ is finite and $\scrP_D$ is the partition of $X$ into points. For a subshift $X$ and $n\ge1$ we write
$$
\scrP^n:=\scrP_{[0,n-1]}=\bigvee_{i=0}^{n-1}\sigma^{-i}\scrP_{\{0\}}
$$
for the partition into the cylinders over $[0,n-1]$.

\emph{Entropy of shifts.} Let $X\subset\calA^\Z$ be a subshift. The partition $\scrP_{\{0\}}$ is a generator for $\sigma$, since the partitions $\sigma^{-n}\scrP_{\{0\}}$, $n\in\Z$, generate the Borel $\sigma$-algebra of $X$. By the Kolmogorov--Sinai theorem \cite{Kolmogorov1958,Sinai1959} (see also \cite[Theorem 4.17]{Walters1982}), every shift-invariant probability $\mu$ on $X$ satisfies
$$
h_\mu(\sigma)=\lim_{n\to\infty}\frac1nH_\mu(\scrP^n)=\inf_{n\ge1}\frac1nH_\mu(\scrP^n).
$$
Cylinder sets are clopen, so each map $\mu\mapsto H_\mu(\scrP^n)$ is continuous for the weak$^*$ topology, and the entropy map, an infimum of continuous functions, is upper semicontinuous (see \cite[Theorem 8.2]{Walters1982}). Hence every continuous $\psi:X\to\R$ has an equilibrium state. We write $P_X(\psi)$ for the pressure of $\psi$ with respect to $\sigma|_X$.

\emph{Relative entropy.} Let $X'$ be a compact metric space with a homeomorphism $\sigma_{X'}$, let $\calA_\aux$ be a finite set, let $X\subset X'\times\calA_\aux^\Z$ be a compact set with $\widehat\sigma(X)=X$, where $\widehat\sigma(x',\omega):=(\sigma_{X'}x',\sigma\omega)$, and let $\pi(x',\omega):=x'$. For $n\ge1$ let $\scrP^n_\aux$ be the partition of $X$ according to $\omega[0,n-1]$. The entropy of a $\widehat\sigma$-invariant probability $\mu$ on $X$ relative to $X'$ is
\begin{equation}
\label{eq:relentropy}
h_\mu(\widehat\sigma\mid X'):=\lim_{n\to\infty}\frac1nH_\mu\bigl(\scrP^n_\aux\mid\pi^{-1}\calB_{X'}\bigr)=\inf_{n\ge1}\frac1nH_\mu\bigl(\scrP^n_\aux\mid\pi^{-1}\calB_{X'}\bigr),
\end{equation}
where $\calB_{X'}$ is the Borel $\sigma$-algebra of $X'$. The limit exists and equals the infimum because $\pi^{-1}\calB_{X'}$ is $\widehat\sigma$-invariant, so that \eqref{eq:shannon} makes the sequence $H_\mu(\scrP^n_\aux\mid\pi^{-1}\calB_{X'})$ subadditive. Since the partitions $\widehat\sigma^{-n}\scrP^1_\aux$, $n\in\Z$, and $\pi^{-1}\calB_{X'}$ generate the Borel $\sigma$-algebra of $X$, the Abramov--Rokhlin formula \cite{AbramovRokhlin1962} (see also \cite[Chapter 4]{Downarowicz2011}) states that
\begin{equation}
\label{eq:abramov}
h_\mu(\widehat\sigma)=h_{\pi_*\mu}(\sigma_{X'})+h_\mu(\widehat\sigma\mid X').
\end{equation}
Let $\eta:=\pi_*\mu$ and let $(\mu_{x'})_{x'\in X'}$ be the sample measures of $\mu$, that is, $\mu(dx',d\omega)=\mu_{x'}(d\omega)\,\eta(dx')$. The conditional entropy given the base is the average of the entropies of the sample measures (see \cite[Chapter II]{Kifer1986} and \cite[(3.19)]{DooleyZhang2015}), that is,
\begin{equation}
\label{eq:relsample}
H_\mu\bigl(\scrP^n_\aux\bigm|\pi^{-1}\calB_{X'}\bigr)=\int H_{\mu_{x'}}(\scrP^n)\,d\eta(x'),
\end{equation}
where on the right $\scrP^n$ is the partition of $\calA_\aux^\Z$ into the cylinders over $[0,n-1]$. We use these facts for $X'=\Sig$, with $\sigma_{X'}=\sigma$, and a subshift $X$ of $(\{-1,1\}\times\calA_\aux)^\Z=\Sig\times\calA_\aux^\Z$, in which case $\widehat\sigma$ is the shift $\sigma$ of $\Sig\times\calA_\aux^\Z$, and for the space $X'=\calX$ of Section~\ref{sec:base}, with $\sigma_{X'}=\sigma_\calX$.

\subsection{Exponential mixing of order \texorpdfstring{$s$}{s}}
\label{sec:mixing-def}
Let $(X,\sigma,m)$ be a probability-preserving system and let $\calB$ be a space of bounded functions with a norm $\norm\cdot_\calB$. Following \cite[(3)]{BassolsCastro2026}, $(X,\sigma,m)$ is \emph{exponentially mixing of order $s$} on $\calB$ if there are $C,c>0$ with
\begin{equation}
\label{eq:mixings}
\Bigl|\int\prod_{j=1}^sF_j\circ\sigma^{n_j}\,dm-\prod_{j=1}^s\int F_j\,dm\Bigr|\le Ce^{-c\Delta(n_1,\dots,n_s)}\prod_{j=1}^s\norm{F_j}_\calB
\end{equation}
for all $F_1,\dots,F_s\in\calB$ and all integers $0\le n_1<\dots<n_s$, where $\Delta(n_1,\dots,n_s):=\min_{1\le j<s}(n_{j+1}-n_j)$. On subshifts we take $\calB=\calF_\theta$, and on manifolds $\calB=\mathcal C^\beta$. We say that a subshift system is exponentially mixing of order $s$ \emph{for every H\"older class} if \eqref{eq:mixings} holds on $\calF_\theta$ for every $0<\theta<1$, with constants depending on $\theta$.

On manifolds, the class of observables does not matter for a fixed order. For $s=2$ the following lemma is \cite[Lemma B.1]{DolgopyatKanigowskiRodriguezHertz2024}, and its interpolation proof, which approximates H\"older observables by $\mathcal C^r$ observables obtained by convolution, applies verbatim to every $s\ge2$.

\begin{lemma}
\label{lem:holderC1}
Let $T$ be a measurable map of a compact Riemannian manifold $M$ preserving a Borel probability $\mu$, and let $s\ge2$ and $r\ge1$ be integers. Then $(T,\mu)$ is exponentially mixing of order $s$ for $\mathcal C^r$ observables if and only if it is exponentially mixing of order $s$ for H\"older observables (see \eqref{eq:allorders}).
\end{lemma}

\section{The process \texorpdfstring{$Z_n$}{Zn}, the marker sets and the auxiliary shift \texorpdfstring{$\mathcal X$}{X}}
\label{sec:BCC}
In this section we construct the auxiliary shift $\calX\subset \Sig\times\{0,1\}^{\Z\times\N_0}$ which is invariant under the shift $$\sigma_\calX(z,h):=(\sigma z,\sigma_\calH h),\text{ where }(\sigma z)_i:=z_{i+1}\text{ and }(\sigma_\calH h)_{i,k}:=h_{i+1,k},$$ together with a $\sigma_\calX$-invariant probability $\zeta$ on $\calX$. In Section~\ref{sec:extension} we extend $\calX$ by auxiliary coordinates in $\{-1,1\}^q$, obtaining $\calX\times(\{-1,1\}^q)^\Z$, and we construct the subshift $$Y\subset(\{-1,1\}\times\{-1,1\}^q)^\Z\text{ as the image }Y:=\widehat\pi(\widehat Y)$$ of a closed $\widehat\sigma$-invariant set $\widehat Y\subset\calX\times(\{-1,1\}^q)^\Z$ under the projection
$$
\widehat\pi:(z,h,\omega)\in \calX\times(\{-1,1\}^q)^\Z\mapsto (z,\omega)\in (\{-1,1\}\times\{-1,1\}^q)^\Z .
$$

\subsection{The process \texorpdfstring{$Z$}{Z} and its law \texorpdfstring{$\nu$}{nu}}
\label{sec:nu}
Let $\Sig:=\{-1,1\}^\Z$ with the left shift $\sigma$. For $i\in\Z$ we write $Z_i$ for the $i$-th coordinate, that is,
$$
Z_i:z\in \Sig\mapsto z_i\in \{-1,1\} ,
$$
and for a finite interval $I=[u,v]$ we write $Z[I]:=(Z_u,\dots,Z_v)$. We write $S_n:=\sum_{i=0}^{n-1}Z_i$, as in \cite[Section 3.1]{BassolsCastro2026}.
We use the following theorem as a black box, for the integer $R$ of Theorem~\ref{thm:A}, and $\nu$ always denotes the measure it gives. It is proved in \cite[Sections 4--8]{BassolsCastro2026}.
 
\begin{theorem}[{\cite[Theorem C]{BassolsCastro2026}}]
\label{thm:BCC}
For every integer $R\ge2$, there exists a full-support $\sigma$-invariant Borel probability measure $\nu$ on $\Sig$ with the following properties.
\begin{enumerate}[label=\textup{(\roman*)},leftmargin=*]
\item\label{BCC:i} Under $\nu$, the coordinate process $(Z_i)_{i\in\Z}$ is stationary, is $\{-1,1\}$-valued, satisfies $\nu(Z_0=-1)=\nu(Z_0=1)=\frac12$, and is $(2R+1)$-tuplewise independent. Thus any collection of at most $2R+1$ coordinates with distinct indices is jointly independent.
\item\label{BCC:ii} There exists $K_R<\infty$ such that, for all disjoint nonempty finite intervals $I_1,\dots,I_R$ ordered from left to right, if
$$
\min_{i<R}\dist(I_i,I_{i+1})>K_R\max_i|I_i| ,
$$
then $Z[I_1],\dots,Z[I_R]$ are jointly independent.
\item\label{BCC:iii} For every $0<\theta<1$, there exist constants $C_{R,\theta},c_\theta>0$ such that
$$
\Bigl|\nu\Bigl(\prod_{i=1}^RF_i\circ\sigma^{n_i}\Bigr)-\prod_{i=1}^R\nu(F_i)\Bigr|\le C_{R,\theta}\,e^{-c_\theta\Delta(n_1,\dots,n_R)}\prod_{i=1}^R\norm{F_i}_\theta
$$
for all $F_1,\dots,F_R\in\calF_\theta(\Sig)$ and all integers $0\le n_1<\dots<n_R$.
\item\label{BCC:iv} There exist $\delta>0$ and $n_0\ge1$ such that
$$
\E_\nu\Bigl[\Bigl(\frac{S_n}{\sqrt n}\Bigr)^{2R+2}\Bigr]\le(2R+1)!!-\delta\quad\text{for every }n\ge n_0 .
$$
In particular, the laws of $S_n/\sqrt n$ under $\nu$ have no Gaussian weak subsequential limit.
\item\label{BCC:v} There exist constants $c,C>0$ and $N_0\ge1$ such that, for every integer $N\ge N_0$, there are integers $0\le n_1<\dots<n_{2R+2}<N$ satisfying
$$
\Delta(n_1,\dots,n_{2R+2})\ge\frac NC\quad\text{and}\quad-\E_\nu\Bigl[\prod_{i=1}^{2R+2}Z_{n_i}\Bigr]\ge cN^{-R-1}.
$$
Hence, $(\Sig,\sigma,\nu)$ is not exponentially mixing of order $2R+2$ in $\calF_\theta(\Sig)$ for any $0<\theta<1$.
\end{enumerate}
\end{theorem}
 We also recall that $(\sigma,\nu)$ is mixing.
\begin{lemma}[{\cite[Corollary 3.2]{BassolsCastro2026}}]
\label{lem:ergodic}
The system $(\sigma,\nu)$ is mixing, hence ergodic.
\end{lemma}

\subsection{Nested marker sets}
\label{sec:markers}
In this subsection we recall the nested marker sets of \cite[Proposition 4.2]{BassolsCastro2026}. They give, for each level $k\ge0$, a random partition of $\Z$ into intervals, which we call the blocks of level $k$. These partitions are nested, that is, every block of level $k\ge1$ is a union of blocks of level $k-1$. We use them in Section~\ref{sec:base} to define the auxiliary shift $\calX$ and an invariant measure $\zeta$ on it. In Section~\ref{sec:extension} we extend $\zeta$ by auxiliary coordinates, whose laws are defined block by block, and we obtain the measure $\mu_0$. The independence of the marker sets on distant intervals is used to prove that $(\sigma,\mu_0)$ is exponentially mixing of order $R$ (Section~\ref{sec:mixing}). The following proposition gives the marker sets for an integer $R'\ge2$, which we fix in \eqref{eq:a} below.

\begin{proposition}[{\cite[Proposition 4.2 and Lemma 4.4]{BassolsCastro2026}}]
\label{prop:markers}
Let $R'\ge2$ be an integer and $L:=2R'+2$. There are a finite constant $D=D(L)\ge1$ and a random sequence $\calH=(\calM_k)_{k\ge0}$ of nested marker sets, with $\Z=\calM_0\supset\calM_1\supset\calM_2\supset\cdots$, such that the following hold, where for $i\in\Z$ and $k\ge0$ we write
$$
\calH^{(k)}_i:=\bigl(\1_{\{i\in\calM_0\}},\dots,\1_{\{i\in\calM_k\}}\bigr) .
$$
\begin{enumerate}[label=\textup{(\roman*)},leftmargin=*]
\item\label{H:length} Almost surely, for every $k\ge1$ the set $\calM_k$ is unbounded from above and from below, and all consecutive points $u<v$ of $\calM_k$ satisfy $(4L)^k<v-u\le3(4L)^k$.
\item\label{H:stat} For every $c\in\Z$ the sequences $(\calM_k+c)_{k\ge0}$ and $(\calM_k)_{k\ge0}$ have the same law, i.e.\ for all $k\ge0$, $p\ge1$ and $i_1,\dots,i_p\in\Z$, the vectors $$\bigl(\calH^{(k)}_{i_1+c},\dots,\calH^{(k)}_{i_p+c}\bigr)\ \text{and}\ \bigl(\calH^{(k)}_{i_1},\dots,\calH^{(k)}_{i_p}\bigr)\ \text{have the same law.}$$
\item\label{H:dep} For every $k\ge0$ the process $\bigl(\calH^{(k)}_i\bigr)_{i\in\Z}$ is $D(4L)^k$-dependent, i.e.\ for all sets $I,J\subset\Z$ with $\min\{|i-j|\ ;\ i\in I,\ j\in J\}>D(4L)^k$, the families $$\bigl(\calH^{(k)}_i\bigr)_{i\in I}\ \text{and}\ \bigl(\calH^{(k)}_j\bigr)_{j\in J}\ \text{are independent.}$$
\end{enumerate}
\end{proposition}

Items \ref{H:stat} and \ref{H:dep} of Proposition~\ref{prop:markers} are (H3) and (H4) of \cite[Proposition 4.2]{BassolsCastro2026}, and item \ref{H:length} is \cite[Lemma 4.4]{BassolsCastro2026} together with (H1) of \cite[Proposition 4.2]{BassolsCastro2026}, where the level-$k$ blocks are the intervals $[u,v-1]$ for consecutive points $u<v$ of $\calM_k$.

We fix the least integer $R'\ge2$ such that $L:=2R'+2$ satisfies
\begin{equation}
\label{eq:a}
L\ge\max\{144R,2304\}\quad\text{and}\quad4(2R+1)(24L)^{2R}e^{-L/1152}\le1 .
\end{equation}
Such an $R'$ exists because $(24L)^{2R}e^{-L/1152}\to0$ as $L\to\infty$. The integers $R'$ and $L$ depend only on $R$. The first condition in \eqref{eq:a} is used in Lemmas~\ref{lem:baseentropy} and~\ref{lem:blocklawsint}, in Proposition~\ref{prop:upper} and in Theorem~\ref{thm:Bformal}, and the second in Lemma~\ref{lem:blocklawsint}.

Let $\calH=(\calM_k)_{k\ge0}$ be as in Proposition~\ref{prop:markers} for the integer $R'$ of \eqref{eq:a}. We identify a sequence $(M_k)_{k\ge0}$ of subsets of $\Z$ with the array $\bigl(\1_{\{i\in M_k\}}\bigr)_{i\in\Z,\,k\ge0}\in\{0,1\}^{\Z\times\N_0}$, and we give $\{0,1\}^{\Z\times\N_0}$ the product topology. Thus $\calH$ is a random variable with values in $\{0,1\}^{\Z\times\N_0}$.

We write
\begin{itemize}[leftmargin=*]
\item $\P_\calH$ for the law of $\calH$ and $\Sigma_\calH$ for the support of $\P_\calH$,
\item $h=(h_{i,k})_{i\in\Z,\,k\ge0}$ for the points of $\{0,1\}^{\Z\times\N_0}$ and $\sigma_\calH$ for the shift $(\sigma_\calH h)_{i,k}:=h_{i+1,k}$,
\item $\calM_k(h):=\bigl\{i\in\Z\ ;\ h_{i,k}=1\bigr\}$ for $h\in\{0,1\}^{\Z\times\N_0}$ and $k\ge0$, so that $\calM_k(\calH)=\calM_k$,
\item $\calP_k(h):=\bigl\{[u,v-1]\ ;\ u<v\text{ are consecutive points of }\calM_k(h)\bigr\}$ for the set of \emph{level-$k$ blocks} of $h$, as in \cite[Section 4.1]{BassolsCastro2026}.
\end{itemize}

Proposition~\ref{prop:markers} describes the marker sets almost surely. The next lemma shows that the same properties hold at every point of the support $\Sigma_\calH$. 
\begin{lemma}
\label{lem:markers}
Let $h\in\Sigma_\calH$. Then the following hold.
\begin{enumerate}[label=\textup{(\roman*)},leftmargin=*]
\item\label{mk:length} $\calM_0(h)=\Z$, and for every $k\ge1$ we have $\calM_k(h)\subset\calM_{k-1}(h)$, the set $\calM_k(h)$ is unbounded from above and from below, and all consecutive points $u<v$ of $\calM_k(h)$ satisfy $(4L)^k<v-u\le3(4L)^k$. Hence, $\calP_0(h)$ is the partition of $\Z$ into singletons, and for $k\ge1$ the set $\calP_k(h)$ is a partition of $\Z$ into intervals $B$ with $(4L)^k<|B|\le3(4L)^k$.
\item\label{mk:children} Let $k\ge1$ and $B=[u,v-1]\in\calP_k(h)$, and let $u=u_0<u_1<\dots<u_n=v$ be the points of $\calM_{k-1}(h)\cap[u,v]$. Then
$$
B=B_1\cup\dots\cup B_n,\quad\text{where}\quad B_j:=[u_{j-1},u_j-1]\in\calP_{k-1}(h)\quad\text{for }1\le j\le n .
$$
We call $B_1,\dots,B_n$ the \emph{children} of $B$, and, as in \cite[Section 4.1]{BassolsCastro2026}, we write $N_k(B):=n$ for their number. They satisfy
$$
\frac{4L}3<n\le12L,\quad\max_{1\le j\le n}|B_j|\le\frac{3|B|}n\quad\text{and}\quad n-1<\frac{|B|}{(4L)^{k-1}} .
$$
\item\label{mk:dep} The measure $\P_\calH$ is $\sigma_\calH$-invariant and $\sigma_\calH(\Sigma_\calH)=\Sigma_\calH$. Moreover, for every $k\ge0$ the process $\bigl((h_{i,0},\dots,h_{i,k})\bigr)_{i\in\Z}$ under $\P_\calH$ is $D(4L)^k$-dependent.
\end{enumerate}
\end{lemma}

\begin{proof}
We start by proving \ref{mk:length}. Let $\Sigma'_\calH$ be the set of $h=(h_{i,k})_{i\in\Z,\,k\ge0}\in\{0,1\}^{\Z\times\N_0}$ such that, for all $i\in\Z$ and $k\ge1$,
\begin{enumerate}[label=\textup{(\alph*)},leftmargin=*]
\item $h_{i,0}=1$ and $h_{i,k}\le h_{i,k-1}$,
\item $\max\bigl\{h_{i',k}\ ;\ i\le i'\le i+3(4L)^k-1\bigr\}=1$,
\item $h_{i,k}\,h_{i',k}=0$ for every $i'\in\Z$ with $0<|i-i'|\le(4L)^k$.
\end{enumerate}
Thus $h\in\Sigma'_\calH$ exactly when $\calM_0(h)=\Z$, $\calM_k(h)\subset\calM_{k-1}(h)$, every interval of length $3(4L)^k$ contains a point of $\calM_k(h)$, and any two distinct points of $\calM_k(h)$ are at distance greater than $(4L)^k$, for every $k\ge1$. For fixed $i$ and $k$, each of these conditions depends on finitely many coordinates of $h$, so it defines a closed set, and $\Sigma'_\calH$ is closed as an intersection of closed sets.

By Proposition~\ref{prop:markers}\ref{H:length}, $\calH$ lies in $\Sigma'_\calH$ almost surely. Hence $\P_\calH(\Sigma'_\calH)=1$, and the support $\Sigma_\calH$ is contained in $\Sigma'_\calH$. For $h\in\Sigma'_\calH$ and $k\ge1$, a set which intersects every interval of length $3(4L)^k$ is unbounded from above and from below and its consecutive points satisfy $v-u\le3(4L)^k$, which gives \ref{mk:length}.

We now prove \ref{mk:children}. Since $\calM_k(h)\subset\calM_{k-1}(h)$, the points $u$ and $v$ belong to $\calM_{k-1}(h)$, so $u_0<\dots<u_n$ are consecutive points of $\calM_{k-1}(h)$. Hence $B_j=[u_{j-1},u_j-1]\in\calP_{k-1}(h)$ for $1\le j\le n$, and $B=B_1\cup\dots\cup B_n$ is a disjoint union, so $|B|=\sum_j|B_j|$. 

Recall that $B\in\calP_k(h)$ satisfies $(4L)^k<|B|\le3(4L)^k$ by \ref{mk:length}. If $k=1$, then every $B_j\in\calP_0(h)$ is a singleton, so $$n=|B|\in(4L,12L],\quad\max_j|B_j|=1=|B|/n\quad\text{and}\quad n-1<|B|=|B|/(4L)^0.$$ If $k\ge2$, then every $B_j\in\calP_{k-1}(h)$ satisfies $(4L)^{k-1}<|B_j|\le3(4L)^{k-1}$ by \ref{mk:length}. Hence
$$
\frac{4L}3=\frac{(4L)^k}{3(4L)^{k-1}}<\frac{|B|}{3(4L)^{k-1}}\le n<\frac{|B|}{(4L)^{k-1}}\le\frac{3(4L)^k}{(4L)^{k-1}}=12L ,
$$
and $\max_j|B_j|\le3(4L)^{k-1}<\frac{3|B|}n$.

Finally we prove \ref{mk:dep}. The invariance of $\P_\calH$ follows from Proposition~\ref{prop:markers}\ref{H:stat}, and since $\sigma_\calH$ is a homeomorphism, the support $\Sigma_\calH$ of the $\sigma_\calH$-invariant measure $\P_\calH$ satisfies $\sigma_\calH(\Sigma_\calH)=\Sigma_\calH$. Since $\P_\calH$ is the law of $\calH$, the process $\bigl((h_{i,0},\dots,h_{i,k})\bigr)_{i\in\Z}$ under $\P_\calH$ has the law of $\bigl(\calH^{(k)}_i\bigr)_{i\in\Z}$, so it is $D(4L)^k$-dependent by Proposition~\ref{prop:markers}\ref{H:dep}.
\end{proof}

By Theorem~\ref{thm:BCC}\ref{BCC:v}, $(\sigma,\nu)$ is not exponentially mixing of order $2R+2$. The next lemma lowers this order to $R+1$, which gives item \ref{A:nomix} of Theorem~\ref{thm:A}.

\begin{lemma}
\label{lem:nomixR1}
For every integer $k\ge2$ there are integers $0\le t_1<\dots<t_{R+1}$ with
\begin{equation}
\label{eq:nomixR1}
\Delta(t_1,\dots,t_{R+1})>2(8R+8)^{k-1}\quad\text{and}\quad-\E_\nu\Bigl[\prod_{\ell=1}^{R+1}Z_{t_\ell}Z_{t_\ell+1}\Bigr]\ge144^{-R-1}(8R+8)^{-2(R+1)k}.
\end{equation}
Hence, $(\sigma,\nu)$ is not exponentially mixing of order $R+1$ in $\calF_\theta(\Sig)$ for any $0<\theta<1$.
\end{lemma}

\begin{proof}
Let $F(z):=z_0z_1$. Since $F$ depends only on $z_0$ and $z_1$, we have $F\in\calF_\theta(\Sig)$ for every $0<\theta<1$, and $\nu(F)=0$ by Theorem~\ref{thm:BCC}\ref{BCC:i}. Since $\prod_{\ell}F\circ\sigma^{t_\ell}=\prod_{\ell}Z_{t_\ell}Z_{t_\ell+1}$, \eqref{eq:nomixR1} gives
$$
\Bigl|\nu\Bigl(\prod_{\ell=1}^{R+1}F\circ\sigma^{t_\ell}\Bigr)-\prod_{\ell=1}^{R+1}\nu(F)\Bigr|\ge144^{-R-1}(8R+8)^{-2(R+1)k}
$$
at times with $\Delta(t_1,\dots,t_{R+1})>2(8R+8)^{k-1}$. This lower bound decays only polynomially in $(8R+8)^k$, so \eqref{eq:mixings} with $s=R+1$ fails for $F$. It remains to prove \eqref{eq:nomixR1}.

Let $\calH'$ be the random nested marker sets of Proposition~\ref{prop:markers} with $R'=R$, that is, with $L$ replaced by $2R+2$, and let $\P'$ be its law. These are the marker sets used in \cite[Section 6]{BassolsCastro2026} to construct $\nu$. In what follows, $h$ denotes a realisation of $\calH'$, and we write $\calP_k(h)$ for its level-$k$ blocks, as in Section~\ref{sec:markers}. By \cite[Proposition 6.3]{BassolsCastro2026}, $\nu=\int\nu_h\,\P'(dh)$, where $(\nu_h)_h$ is a measurable family of Borel probabilities on $\Sig$. For $\P'$-almost every $h$ the following hold.
\begin{itemize}[leftmargin=*]
\item By \cite[(H1) and (H2)]{BassolsCastro2026}, every $B\in\calP_k(h)$ with $k\ge1$ satisfies $(8R+8)^k<|B|\le3(8R+8)^k$, and $B=B_1\cup\dots\cup B_r$ with $r>2R+2$, where $B_1,\dots,B_r\in\calP_{k-1}(h)$ are ordered from left to right.
\item By \cite[Lemma 6.4]{BassolsCastro2026}, $\E_{\nu_h}\bigl[\prod_{i=1}^{2R+2}Z_{j_i}\bigr]\le0$ for all integers $j_1<\dots<j_{2R+2}$, and if $j_i\in B_i$ for $1\le i\le2R+2$, for some $B\in\calP_k(h)$ as above, then
$$
\E_{\nu_h}\Bigl[\prod_{i=1}^{2R+2}Z_{j_i}\Bigr]=-\frac{2^{-3(R+1)}}{(|B_1|\cdots|B_{2R+2}|)^{1/2}} .
$$
\end{itemize}

Fix $k\ge2$, let $N:=6(8R+8)^k$, and let
$$
\mathcal T:=\bigl\{(t_1,\dots,t_{R+1})\in\Z^{R+1}\ ;\ 0\le t_1<\dots<t_{R+1}<N\ \text{and}\ \Delta(t_1,\dots,t_{R+1})>2(8R+8)^{k-1}\bigr\} .
$$

Fix $h$ as above. The block of $\calP_k(h)$ containing $0$ and the next block have length at most $3(8R+8)^k$, so there is $B\in\calP_k(h)$ with $B\subset[0,N)$. Write $$B_i=[u_{i-1},u_i-1]\ \text{and }t_\ell:=u_{2\ell-1}-1\ \text{for}\ 1\le\ell\le R+1.$$ Then $$t_\ell\in B_{2\ell-1}, \ t_\ell+1\in B_{2\ell}\text{ and }t_{\ell+1}-t_\ell=|B_{2\ell}|+|B_{2\ell+1}|>2(8R+8)^{k-1}.$$ Hence $(t_1,\dots,t_{R+1})\in\mathcal T$ and, since $|B_i|\le3(8R+8)^{k-1}$,
$$
\nu_h\Bigl(\prod_{\ell=1}^{R+1}F\circ\sigma^{t_\ell}\Bigr)\le-24^{-R-1}(8R+8)^{-(k-1)(R+1)} .
$$
For every $(t'_1,\dots,t'_{R+1})\in\mathcal T$ the $2R+2$ points $t'_\ell,t'_\ell+1$ are distinct, so $\nu_h\bigl(\prod_{\ell}F\circ\sigma^{t'_\ell}\bigr)\le0$. Since $\mathcal T$ does not depend on $h$, integrating over $h$ gives
$$
-\sum_{t'\in\mathcal T}\nu\Bigl(\prod_{\ell=1}^{R+1}F\circ\sigma^{t'_\ell}\Bigr)\ge24^{-R-1}(8R+8)^{-(k-1)(R+1)} .
$$
Since $|\mathcal T|\le N^{R+1}=6^{R+1}(8R+8)^{k(R+1)}$, some $(t_1,\dots,t_{R+1})\in\mathcal T$ satisfies $$-\nu\bigl(\prod_{\ell}F\circ\sigma^{t_\ell}\bigr)\ge144^{-R-1}(8R+8)^{-2(R+1)k},$$ which is \eqref{eq:nomixR1}.
\end{proof}

\subsection{The auxiliary shift \texorpdfstring{$\calX$}{X}}
\label{sec:base}
We define $\calX:=\Sig\times\Sigma_\calH$, the probability $\zeta:=\nu\times\P_\calH$ on $\calX$, and the shift
$$
\sigma_\calX:(z,h)\in \calX\mapsto (\sigma z,\sigma_\calH h)\in \calX .
$$
For $s\ge0$, $a\in\{-1,1\}^{[-s,s]}$ and $b\in\{0,1\}^{[-s,s]\times[0,s]}$ let
$$
C_s(a,b):=\bigl\{(z,h)\in\calX\ ;\ z[-s,s]=a\ \text{and}\ h_{i,k}=b_{i,k}\ \text{for all }i\in[-s,s]\text{ and }0\le k\le s\bigr\},
$$
and let
$$
\scrQ_s:=\bigl\{C_s(a,b)\ ;\ a\in\{-1,1\}^{[-s,s]},\ b\in\{0,1\}^{[-s,s]\times[0,s]}\ \text{and}\ C_s(a,b)\ne\emptyset\bigr\}.
$$
Thus $\scrQ_s$ is a finite partition of $\calX$ into clopen sets. The partitions $\scrQ_s$ increase with $s$, the family $\{\emptyset\}\cup\bigcup_{s\ge0}\scrQ_s$ is closed under intersections, and $\bigcup_{s\ge0}\scrQ_s$ generates the Borel $\sigma$-algebra of $\calX$.

\begin{lemma}
\label{lem:base}
For $s\ge0$ define $\varrho_s:=(K_R+1)(2s+1)+D(4L)^s$, where $K_R$ is the constant of Theorem~\ref{thm:BCC}\ref{BCC:ii} and $D$ is the constant of Proposition~\ref{prop:markers}. Then the following hold.
\begin{enumerate}[label=\textup{(\roman*)},leftmargin=*]
\item\label{base:indep} Let $s\ge0$, let $C,C'\in\scrQ_s$, and let $i,i'\in\Z$ with $|i-i'|>\varrho_s$. Then the functions $\1_C\circ\sigma_\calX^i$ and $\1_{C'}\circ\sigma_\calX^{i'}$ are independent under $\zeta$.
\item\label{base:mixing} The measure $\zeta$ is $\sigma_\calX$-invariant, and $(\sigma_\calX,\zeta)$ is mixing, hence ergodic.
\end{enumerate}
\end{lemma}

\begin{proof}
We start by proving \ref{base:indep}. We may assume $i<i'$, and we write $I:=[i-s,i+s]$ and $I':=[i'-s,i'+s]$. By the definition of $\scrQ_s$, we have $\1_C\circ\sigma_\calX^i(z,h)=f(z)\,g(h)$, where $f$ depends only on $z[I]$ and $g$ only on $(h_{l,k})_{l\in I,\,0\le k\le s}$, and similarly $\1_{C'}\circ\sigma_\calX^{i'}(z,h)=f'(z)\,g'(h)$ with $I'$ in place of $I$. Since $i'-i>\varrho_s$,
$$
\dist(I,I')=i'-i-2s>K_R(2s+1)+D(4L)^s .
$$
Hence $z[I]$ and $z[I']$ are independent under $\nu$ by Theorem~\ref{thm:BCC}\ref{BCC:ii}, applied to $I$, $I'$ and the intervals $I'+jd$, $1\le j\le R-2$, where $d:=(K_R+1)(2s+1)$. Moreover $(h_{l,k})_{l\in I,\,0\le k\le s}$ and $(h_{l,k})_{l\in I',\,0\le k\le s}$ are independent under $\P_\calH$ by Lemma~\ref{lem:markers}\ref{mk:dep}. Since $\zeta=\nu\times\P_\calH$, this gives
\begin{align*}
\zeta\bigl(\sigma_\calX^{-i}C\cap\sigma_\calX^{-i'}C'\bigr)&=\nu(ff')\,\P_\calH(gg')=\nu(f)\,\nu(f')\,\P_\calH(g)\,\P_\calH(g')=\zeta\bigl(\sigma_\calX^{-i}C\bigr)\,\zeta\bigl(\sigma_\calX^{-i'}C'\bigr) .
\end{align*}

Finally we prove \ref{base:mixing}. Since $\nu$ is $\sigma$-invariant and $\P_\calH$ is $\sigma_\calH$-invariant by Lemma~\ref{lem:markers}\ref{mk:dep}, the measure $\zeta$ is $\sigma_\calX$-invariant. Let $C\in\scrQ_s$ and $C'\in\scrQ_{s'}$, and let $m:=\max\{s,s'\}$. Then $C$ and $C'$ are unions of elements of $\scrQ_m$, so \ref{base:indep} with $i=0$ and $i'=n$ gives $\zeta(C\cap\sigma_\calX^{-n}C')=\zeta(C)\,\zeta(C')$ for every $n>\varrho_m$. The family $\{\emptyset\}\cup\bigcup_{s\ge0}\scrQ_s$ is a semi-algebra which generates the Borel $\sigma$-algebra of $\calX$, so $(\sigma_\calX,\zeta)$ is mixing (see \cite[Theorem 1.17]{Walters1982}), and in particular ergodic.
\end{proof}

The lemma below is a standard entropy computation, which counts the elements of the partitions $\scrQ_s$.
\begin{lemma}
\label{lem:baseentropy}
For every $\sigma_\calX$-invariant Borel probability $\eta$ on $\calX$,
$$
h_\eta(\sigma_\calX)\le\log2+\sum_{k\ge1}(4L)^{-k}(1+k\log(4L))\le2\log2 .
$$
\end{lemma}

\begin{proof}
Let $s\ge1$ and $n\ge2$, and define $N:=n+2s$. The partition $\bigvee_{i=0}^{n-1}\sigma_\calX^{-i}\scrQ_s$ is the partition of $\calX$ according to the values of $z[-s,n-1+s]$ and of $(h_{i,k})_{i\in[-s,n-1+s],\,0\le k\le s}$. We count these values on $\calX$.
\begin{itemize}[leftmargin=*]
\item The word $z[-s,n-1+s]$ takes at most $2^N$ values.
\item The row $(h_{i,0})_{i\in[-s,n-1+s]}$ takes one value, since $\calM_0(h)=\Z$.
\item Let $1\le k\le s$. By Lemma~\ref{lem:markers}\ref{mk:length}, the set $\calM_k(h)\cap[-s,n-1+s]$ is a subset of an interval of length $N$ whose distinct points are at distance greater than $(4L)^k$. It has at most $$r_k:=\lfloor N/(4L)^k\rfloor+1\le N\ \text{points},$$ so the row $(h_{i,k})_{i\in[-s,n-1+s]}$ takes at most $\sum_{i=0}^{r_k}\binom Ni$ values. Since $r_k/N\le1$,
$$
\Bigl(\frac{r_k}N\Bigr)^{r_k}\sum_{i=0}^{r_k}\binom Ni\le\sum_{i=0}^{r_k}\binom Ni\Bigl(\frac{r_k}N\Bigr)^i\le\Bigl(1+\frac{r_k}N\Bigr)^N\le e^{r_k},
$$
and, since $N/r_k< (4L)^k$, the number of values is at most
$$
\Bigl(\frac{eN}{r_k}\Bigr)^{r_k}\le e^{r_k(1+k\log(4L))}\le e^{(N/(4L)^k+1)(1+k\log(4L))}.
$$
\end{itemize}
Hence the bound $H_\mu(\xi)\le\log\#\xi$ of \eqref{eq:shannon} and the three items above give
\begin{align*}
H_\eta\Bigl(\bigvee_{i=0}^{n-1}\sigma_\calX^{-i}\scrQ_s\Bigr)&\le\log\#\Bigl(\bigvee_{i=0}^{n-1}\sigma_\calX^{-i}\scrQ_s\Bigr)\le\log\Bigl(2^N\prod_{k=1}^se^{(N/(4L)^k+1)(1+k\log(4L))}\Bigr)\\
&=N\log2+\sum_{k=1}^s\Bigl(\frac N{(4L)^k}+1\Bigr)(1+k\log(4L)).
\end{align*}
Dividing by $n$ and letting $n\to\infty$ gives $h_\eta(\sigma_\calX,\scrQ_s)\le\log2+\sum_{k=1}^s(4L)^{-k}(1+k\log(4L))$.

The partitions $\scrQ_s$ increase and generate the Borel $\sigma$-algebra of $\calX$, so $$h_\eta(\sigma_\calX)=\lim_{s\to\infty}h_\eta(\sigma_\calX,\scrQ_s)$$ (see \cite[Theorem 4.22]{Walters1982}), which gives the first inequality. For the second, $\sum_{k\ge1}(4L)^{-k}=\frac1{4L-1}$ and $\sum_{k\ge1}k(4L)^{-k}=\frac{4L}{(4L-1)^2}\le\frac1{2L}$, and $4L\ge288$ by \eqref{eq:a}, so
$$
\sum_{k\ge1}(4L)^{-k}(1+k\log(4L))\le\frac1{287}+\frac{2\log288}{288}<\frac1{20}<\log2 .
$$
\end{proof}

\section{Measures on finite products}
\label{sec:products}

This section is purely probabilistic and can be read independently of Section~\ref{sec:BCC}. In Proposition~\ref{prop:products} we construct two probabilities on a finite product. Both have the same marginals as the product measure on any set of at most $2R$ coordinates, and both give full measure to the set where a sum of centred functions of the coordinates is small. One has entropy close to that of the product measure, and the other has small support. We use these probabilities to define the laws of the blocks in Section~\ref{sec:extension}. The marginals are used to prove mixing (Section~\ref{sec:mixing}), and the entropy and support bounds are used to prove the entropy gap (Section~\ref{sec:entropy}).

Let $E_1,\dots,E_n$ be finite sets and let $\mathbb{Q}_j$ be a probability measure with full support on $E_j$ for each $j$. Set
$$
E:=E_1\times\dots\times E_n\quad\text{and}\quad \mathbb{Q}:=\mathbb{Q}_1\times\dots\times\mathbb{Q}_n .
$$
For $x\in E$ and $S\subset[1,n]$, write $x_S:=(x_j)_{j\in S}$. For $g:E\to\mathbb{R}$, let $\mathbb{E}_Sg$ denote the conditional expectation of $g$ under $\mathbb{Q}$ given the coordinates in $S$. Thus
$$
(\mathbb{E}_Sg)(x)=\frac{1}{\mathbb{Q}(\{y\in E\mid y_S=x_S\})}\sum_{\substack{y\in E\\ y_S=x_S}}g(y)\,\mathbb{Q}(\{y\}) .
$$
In particular, $\mathbb{E}_\emptyset g=\mathbb{Q}(g)$ and $\mathbb{E}_{[1,n]}g=g$.

The projection onto the terms of order at most $2R$ in the Hoeffding--Efron--Stein decomposition \cite{Hoeffding1948,EfronStein1981} (see also \cite[Section~8.3]{ODonnell2014}) is given by
$$
\Pi_{2R}:=\sum_{\substack{S\subset[1,n]\\ |S|\le 2R}}\ \sum_{T\subset S}(-1)^{|S\setminus T|}\,\mathbb{E}_T .
$$

\begin{lemma}
\label{lem:decomposition}
For all $S,U\subset[1,n]$ the following hold.
\begin{enumerate}[label=\textup{(\roman*)},leftmargin=*]
\item\label{dec:prod} $\E_S\E_U=\E_{S\cap U}$.
\item\label{dec:pi} If $|U|\le2R$, then $\E_U\Pi_{2R}=\E_U$.
\end{enumerate}
\end{lemma}

\begin{proof}
We start by proving \ref{dec:prod}. Using that $\Q$ is a product measure, $(\E_Ug)(x)$ is the integral of $g(x_U,y_{[1,n]\setminus U})$ with respect to $\prod_{j\notin U}\Q_j$, and Fubini's we obtain that $\E_S\E_U=\E_{S\cap U}$.

Finally we prove \ref{dec:pi}. By \ref{dec:prod},
$$
\E_U\Pi_{2R}=\sum_{|S|\le2R}\ \sum_{T\subset S}(-1)^{|S\setminus T|}\E_{U\cap T} .
$$
Let $|S|\le2R$. If there is $s\in S\setminus U$, then for every $T\subset S\setminus\{s\}$ the terms of $T$ and of $T\cup\{s\}$ have the same operator $\E_{U\cap T}$ and opposite signs, so the inner sum vanishes. If $S\subset U$, then $U\cap T=T$. Every $S\subset U$ satisfies $|S|\le2R$, so, exchanging the sums,
$$
\E_U\Pi_{2R}=\sum_{S\subset U}\ \sum_{T\subset S}(-1)^{|S\setminus T|}\E_T=\sum_{T\subset U}\E_T\sum_{T\subset S\subset U}(-1)^{|S\setminus T|}=\E_U ,
$$
since the inner sum equals $1$ if $T=U$ and $0$ otherwise.
\end{proof}

The following proposition is technical. It is applied in Definition~\ref{def:blocklaws} to define the laws of the blocks, and its properties are used in Lemma~\ref{lem:blocklawsint}. Its proof uses the bounded differences inequality of McDiarmid \cite{McDiarmid1989} (see also \cite[Theorem 6.2]{BoucheronLugosiMassart2013}).

\begin{proposition}
\label{prop:products}
Assume the following.
\begin{itemize}[leftmargin=*]
\item The integer $n$ satisfies
\begin{equation}
\label{eq:nk}
n\ge\max\{192R,3072\}\quad\text{and}\quad(2R+1)(2n)^{2R} e^{-n/1536}\le\frac14 .
\end{equation}
\item $\mathcal V$ is a finite set, and $\norm\cdot$ and $\langle\cdot,\cdot\rangle$ are the Euclidean norm and inner product on $\R^{\mathcal V}$.
\item $b_1,\dots,b_n\ge0$ satisfy $\max_jb_j\le\frac3n\sum_jb_j$.
\item For $1\le j\le n$, the function $X_j:E_j\to\R^{\mathcal V}$ satisfies $\sum_{w\in E_j}\Q_j(\{w\})X_j(w)=0$ and $\norm{X_j}\le b_j$.
\end{itemize}
Define
$$
G:=\Bigl\{x\in E\ ;\ \Bigl\|\sum_{j=1}^nX_j(x_j)\Bigr\|\le\frac18\sum_{j=1}^nb_j\Bigr\}.
$$
Then there are probabilities $\widetilde\Q$ and $\overline\Q$ on $E$ with the following properties.
\begin{enumerate}[label=\textup{(\roman*)},leftmargin=*]
\item\label{cp:supp} $\widetilde\Q(G)=\overline\Q(G)=1$.
\item\label{cp:marg} For every $S\subset[1,n]$ with $|S|\le2R$ and every $x\in E$,
$$
\widetilde\Q\bigl(\{y\in E\ ;\ y_S=x_S\}\bigr)=\overline\Q\bigl(\{y\in E\ ;\ y_S=x_S\}\bigr)=\prod_{j\in S}\Q_j(\{x_j\}).
$$
\item\label{cp:entropy} For $F\subset[1,n]$ let $\scrP_F$ be the partition of $E$ into the cylinders $\{y\in E\ ;\ y_F=x_F\}$, where $x\in E$ (see Section~\ref{sec:entropy-def}). Then
$$
\sum_{j=1}^nH_{\widetilde\Q}\bigl(\scrP_{\{j\}}\bigr)-H_{\widetilde\Q}\bigl(\scrP_{[1,n]}\bigr)\le\log\frac43 .
$$
\item\label{cp:support} $\#\supp\overline\Q\le2^n\bigl(\max_j\#E_j\bigr)^{2R}$.
\end{enumerate}
\end{proposition}

\begin{proof}
Define $B:=\sum_jb_j$ and $\delta:=e^{-n/1536}$. We divide the proof into four steps.

\begin{step}[1]
We show that, for every $S\subset[1,n]$ with $|S|\le2R$ and every $x\in E$,
$$
\Q\bigl(\{y\in E\ ;\ y_S=x_S,\ y\notin G\}\bigr)\le\delta\,\Q\bigl(\{y\in E\ ;\ y_S=x_S\}\bigr).
$$
\end{step}

\begin{stepproof}
If $B=0$, then $b_j=0$, and hence $X_j\equiv0$, for every $j$, so $G=E$ and there is nothing to prove. Assume $B>0$. Let $\Q^x(\cdot):=\Q(\;\cdot \mid \{y\in E\ ;\ y_S=x_S\})$, and define
$$
f(y):=\Bigl\|\sum_{j\notin S}X_j(y_j)\Bigr\|,\quad\text{where }y\in E .
$$
Since $n\ge192R$, for every $y\in E$ with $y_S=x_S$,
$$
\Bigl\|\sum_{j\in S}X_j(y_j)\Bigr\|\le2R\max_jb_j\le\frac{6R}nB\le\frac B{32},
$$
so $y\notin G$ implies $f(y)>\frac{3B}{32}$.

Under $\Q^x$, the coordinates $y_j$ with $j\notin S$ are independent with laws $\Q_j$. Hence the vectors $X_j(y_j)$, $j\notin S$, are independent and centred, so
$$
\Q^x\bigl(\langle X_j(y_j),X_{j'}(y_{j'})\rangle\bigr)=0\quad\text{for }j\ne j' .
$$
Moreover $\sum_jb_j^2\le\max_jb_j\,B\le\frac3nB^2$. By Jensen's inequality and since $n\ge3072$,
$$
\Q^x(f)\le \Q^x\bigl(f^2\bigr)^{1/2}=\Bigl(\sum_{j\notin S}\Q^x\bigl(\norm{X_j(y_j)}^2\bigr)\Bigr)^{1/2}\le\Bigl(\frac3n\Bigr)^{1/2}B\le\frac B{32} .
$$
If $y,y'\in E$ satisfy $y_S=y'_S=x_S$ and differ only in the coordinate $j\notin S$, then $|f(y)-f(y')|\le\norm{X_j(y_j)-X_j(y'_j)}\le2b_j$. Hence the bounded differences inequality (see \cite[Theorem 6.2]{BoucheronLugosiMassart2013}) gives
$$
\Q^x\bigl(\{y\ ;\ y\notin G\}\bigr)\le \Q^x\Bigl(f-\Q^x(f)\ge\frac B{16}\Bigr)\le\exp\Bigl(-\frac{2(B/16)^2}{\sum_{j\notin S}(2b_j)^2}\Bigr)\le\exp\Bigl(-\frac{B^2/128}{12B^2/n}\Bigr)=\delta .
$$
\end{stepproof}

\begin{step}[2]
We show that the operator $\Phi g:=\Pi_{2R}\bigl(\1_{E\setminus G}\,g\bigr)$, acting on functions $g:E\to\R$, satisfies $\norm{\Phi g}_\infty\le\frac14\norm g_\infty$, and that $v:=\sum_{r\ge0}\Phi^r1$ satisfies $\norm{v-1}_\infty\le\frac13$ and $\E_S(\1_Gv)=1$ for every $S\subset[1,n]$ with $|S|\le2R$.
\end{step}

\begin{stepproof}
Write $b:=\1_{E\setminus G}$. Let $T\subset[1,n]$ with $|T|\le2R$ and $x\in E$. By Step~1,
$$
\bigl|\E_T(bg)(x)\bigr|\le\norm g_\infty\,\E_T(b)(x)=\norm g_\infty\,\frac{\Q(\{y\in E\ ;\ y_T=x_T,\ y\notin G\})}{\Q(\{y\in E\ ;\ y_T=x_T\})}\le\delta\norm g_\infty .
$$
For $0\le s\le2R$ there are $\binom ns$ sets $S$ with $|S|=s$, each with $2^s$ subsets $T$, and $2^s\binom ns\le(2n)^{2R}$. Hence, by \eqref{eq:nk},
$$
\norm{\Phi g}_\infty\le\sum_{s=0}^{2R}2^s\binom ns\,\delta\norm g_\infty\le(2R+1)(2n)^{2R}\delta\norm g_\infty\le\frac14\norm g_\infty .
$$
Therefore $\norm{\Phi^r1}_\infty\le4^{-r}$, the series $v$ converges, $\norm{v-1}_\infty\le\sum_{r\ge1}4^{-r}=\frac13$, and $v=1+\Phi v$. Let $|S|\le2R$. By Lemma~\ref{lem:decomposition}\ref{dec:pi}, $\E_S(bv)=\E_S\Pi_{2R}(bv)=\E_S\Phi v$, so
$$
\E_S(\1_Gv)=\E_Sv-\E_S(bv)=\E_S(v-\Phi v)=\E_S1=1 .
$$
\end{stepproof}

\begin{step}[3]
We define $\widetilde\Q$ and prove \ref{cp:supp}--\ref{cp:entropy} for it.
\end{step}

\begin{stepproof}
Define $\widetilde\Q(\{x\}):=\1_G(x)\,v(x)\,\Q(\{x\})$, where $x\in E$. By Step~2, $v\ge\frac23$, so $\widetilde\Q(\{x\})\ge0$ and $\widetilde\Q(G)=\widetilde\Q(E)$. Let $|S|\le2R$ and $x\in E$. By Step~2,
\begin{align*}
\widetilde\Q\bigl(\{y\ ;\ y_S=x_S\}\bigr)&=\sum_{\substack{y\in E\\ y_S=x_S}}\1_G(y)v(y)\Q(\{y\})\\
&=\Q\bigl(\{y\ ;\ y_S=x_S\}\bigr)\,\E_S(\1_Gv)(x)=\prod_{j\in S}\Q_j(\{x_j\}),
\end{align*}
which is \ref{cp:marg}. For $S=\emptyset$ it gives $\widetilde\Q(E)=1$, so $\widetilde\Q$ is a probability with $\widetilde\Q(G)=1$. By \ref{cp:marg} with $|S|=1$, the marginal of $\widetilde\Q$ on the coordinate $j$ is $\Q_j$, so
$$
H_{\widetilde\Q}\bigl(\scrP_{\{j\}}\bigr)=-\sum_{w\in E_j}\Q_j(\{w\})\log\Q_j(\{w\}).
$$
The partition $\scrP_{[1,n]}$ is the partition of $E$ into points, and $\log\widetilde\Q(\{x\})=\log v(x)+\sum_j\log \Q_j(\{x_j\})$ when $\widetilde\Q(\{x\})>0$. Hence
$$
H_{\widetilde\Q}\bigl(\scrP_{[1,n]}\bigr)=-\sum_{x\in G}\widetilde\Q(\{x\})\log v(x)-\sum_{j=1}^n\sum_{w\in E_j}\Q_j(\{w\})\log \Q_j(\{w\})\ge-\log\frac43+\sum_{j=1}^nH_{\widetilde\Q}\bigl(\scrP_{\{j\}}\bigr),
$$
since $v\le\frac43$ by Step~2. This is \ref{cp:entropy}.
\end{stepproof}

\begin{step}[4]
Finally we define $\overline\Q$ and prove \ref{cp:supp}, \ref{cp:marg} and \ref{cp:support} for it.
\end{step}

\begin{stepproof}
Let $\calW$ be the set of pairs $(S,w)$ with $S\subset[1,n]$, $1\le|S|\le2R$ and $w\in\prod_{j\in S}E_j$. For $x\in E$ define $\upsilon(x)\in\R^\calW$ by $\upsilon(x)_{(S,w)}:=\1_{\{x_S=w\}}$, and define $m\in\R^\calW$ by $m_{(S,w)}:=\prod_{j\in S}\Q_j(\{w_j\})$. A probability $\Q'$ on $E$ satisfies \ref{cp:marg} exactly when $\sum_x\Q'(\{x\})\upsilon(x)=m$. By Step~3, the vector $m=\sum_{x\in G}\widetilde\Q(\{x\})\upsilon(x)$ lies in the convex hull of $\{\upsilon(x)\ ;\ x\in G\}$. By Carath\'eodory's theorem \cite[Theorem 17.1]{Rockafellar1970}, $m$ is a convex combination of at most $\#\calW+1$ vectors $\upsilon(x)$ with $x\in G$, and the probability $\overline\Q$ with these weights at these points satisfies \ref{cp:supp} and \ref{cp:marg}. Write $M:=\max_j\#E_j\ge1$. Then
$$
\#\supp\overline\Q\le\#\calW+1\le1+\sum_{s=1}^{2R}\binom nsM^s\le\Bigl(1+\sum_{s=1}^{2R}\binom ns\Bigr)M^{2R}\le2^nM^{2R} ,
$$
which is \ref{cp:support}.
\end{stepproof}
\end{proof}

\section{Construction of \texorpdfstring{$Y$}{Y}}
\label{sec:extension}

In this section we construct the subshift $(Y,\sigma)$ and the invariant measure $\mu_0$ on $Y$ of the symbolic theorem (Theorem~\ref{thm:Bformal}), which are the heart of the proof of Theorem~\ref{thm:A}. The subshift $Y$ is obtained by adding auxiliary coordinates to $\Sig$, and the first coordinate of $\mu_0$ has law $\nu$, so $\mu_0$ inherits the failure of the CLT. In Sections~\ref{sec:mixing} and~\ref{sec:entropy} we prove that $(\sigma,\mu_0)$ is exponentially mixing of order $R$ and that it is the unique equilibrium state on $Y$ of every continuous potential with small oscillation. In Section~\ref{sec:realisation} we realise $(Y,\sigma)$ as an invariant set of a smooth diffeomorphism.

Throughout Sections~\ref{sec:extension}, \ref{sec:mixing} and \ref{sec:entropy}, $q\ge1$ is a fixed integer and $L$ is given by \eqref{eq:a}. We write
$$
\calA_\aux:=\{-1,1\}^q\quad\text{and}\quad\calA:=\{-1,1\}\times\calA_\aux .
$$
A point of $\calA^\Z$ is a pair $x=(z,\omega)$ with $z\in\Sig$ and $\omega\in\calA_\aux^\Z$, so that at time $i$ the point $x$ has the base coordinate $z_i$ and the auxiliary coordinate $\omega_i$. We write
$$
\pi_1:(z,\omega)\in \calA^\Z\mapsto z\in \Sig
$$
for the projection on the first coordinate, and
$$
\widehat\sigma:(z,h,\omega)\in \calX\times\calA_\aux^\Z\mapsto \bigl(\sigma_\calX(z,h),\sigma\omega\bigr)\in \calX\times\calA_\aux^\Z
$$
for the shift of $\calX\times\calA_\aux^\Z$, as in Section~\ref{sec:entropy-def}. To construct $Y$ we work on the larger space $\calX\times\calA_\aux^\Z$, which also records the marker sets, and we obtain $Y$ as the image of a closed $\widehat\sigma$-invariant set $\widehat Y\subset\calX\times\calA_\aux^\Z$ under the projection $\widehat\pi(z,h,\omega)=(z,\omega)$.

We now sketch the construction of $Y$ from the auxiliary shift $(\calX,\sigma_\calX,\zeta)$ of Section~\ref{sec:base}. The precise statements are in the sections indicated in each step.
\begin{enumerate}[label=\textup{(\arabic*)},leftmargin=*]
\item \emph{Recursive construction} (Sections~\ref{sec:blocklaws} and~\ref{sec:construction}). For $(z,h)\in\calX$ and every block $B\in\calP_k(h)$ we define a probability $\P^{z,h}_B$ on $\calA_\aux^B$, obtained from the probabilities of the blocks of level $k-1$ contained in $B$ by Proposition~\ref{prop:products}. Under $\P^{z,h}_B$ any $2R$ of these blocks are independent, the local statistics of every $w\in\supp\P^{z,h}_B$ are close to their means, and $\P^{z,h}_B$ has almost maximal entropy if $(z,h)$ is typical for $\zeta$ along $B$ and small support otherwise.
\item \emph{The extension} (Section~\ref{sec:sample}). The sample measures are the limits
$$
\mu_{z,h}:=\lim_{k\to\infty}\ \bigotimes_{B\in\calP_k(h)}\P^{z,h}_B ,
$$
and we define
$$
\widehat Y:=\bigl\{(z,h,\omega)\ ;\ \omega[B]\in\supp\P^{z,h}_B\text{ for all }B\bigr\}\quad\text{and}\quad Y:=\widehat\pi(\widehat Y),
$$
$$
\widehat\mu_0(dz,dh,d\omega):=\mu_{z,h}(d\omega)\,\zeta(dz,dh)\quad\text{and}\quad\mu_0:=\widehat\pi_*\widehat\mu_0 .
$$
\item \emph{Uniqueness of the lift} (Section~\ref{sec:faithful}). Every $\sigma_\calX$-invariant probability $\eta$ on $\calX$ has exactly one $\widehat\sigma$-invariant lift to $\widehat Y$, namely $\mu_{z,h}(d\omega)\,\eta(dz,dh)$. Moreover $Y$ has no fixed point.
\item \emph{Mixing} (Section~\ref{sec:mixing}). Exact independence of the laws of $2R$ blocks, the finite dependence of the marker sets and Theorem~\ref{thm:BCC}\ref{BCC:ii} give exact independence of $R$ distant blocks of coordinates under $\mu_0$, and hence, for $0<\theta<1$ and $F_1,\dots,F_R\in\calF_\theta(\calA^\Z)$,
$$
\Bigl|\mu_0\Bigl(\prod_{s=1}^RF_s\circ\sigma^{t_s}\Bigr)-\prod_{s=1}^R\mu_0(F_s)\Bigr|\le C_\theta e^{-c_\theta\Delta(t_1,\dots,t_R)}\prod_{s=1}^R\norm{F_s}_\theta .
$$
\item \emph{Entropy gap} (Section~\ref{sec:entropy}). For every ergodic $\sigma$-invariant probability $\rho\ne\mu_0$ on $Y$,
$$
h_{\mu_0}(\sigma)\ge\Bigl(\frac78q-1\Bigr)\log2\quad\text{and}\quad h_\rho(\sigma)\le\Bigl(2+\frac q8\Bigr)\log2 ,
$$
so $h_{\mu_0}(\sigma)-h_\rho(\sigma)\ge\frac q2\log2$ when $q\ge12$. Hence, for every continuous $\psi$ with $\sup_Y\psi-\inf_Y\psi<\frac q2\log2$ and every such $\rho$,
$$
h_\rho(\sigma)+\rho(\psi)\le h_{\mu_0}(\sigma)-\frac q2\log2+\sup_Y\psi<h_{\mu_0}(\sigma)+\inf_Y\psi\le h_{\mu_0}(\sigma)+\mu_0(\psi).
$$
Since almost every ergodic component of an equilibrium state is an equilibrium state (see Section~\ref{sec:entropy-def}), $\mu_0$ is the unique equilibrium state of $\psi|_Y$ (Theorem~\ref{thm:Bformal}\ref{B:unique}).
\item \emph{Smooth realisation} (Section~\ref{sec:realisation}). We choose $q$ with $$\frac q2\log2>\sup_M\varphi-\inf_M\varphi.$$ Using a horseshoe and the neutral isolation of Lin and Tian \cite{LinTian2026}, we realise $\sigma|_Y$ on an invariant set of a diffeomorphism $T$ of $M$ whose other invariant probabilities are supported on fixed points. By (5), the image of $\mu_0$ is the unique equilibrium state of $(T,\varphi)$.
\end{enumerate}

\subsection{Definition of \texorpdfstring{$Y$}{Y}}
\label{sec:blocklaws}
Recall that, for $h\in\Sigma_\calH$ and $k\ge0$, $\calP_k(h)$ is the partition of $\Z$ of Lemma~\ref{lem:markers}\ref{mk:length}, given by
$$
\calP_k(h)=\bigl\{[u,v-1]\ ;\ u<v\text{ are consecutive points of }\calM_k(h)\bigr\},
$$
where $\calM_k(h)=\{i\in\Z\ ;\ h_{i,k}=1\}$. We write
$$
J_k(h):=\text{the unique element of }\calP_k(h)\text{ which contains }0 .
$$

In this subsection we define the subshift $Y$ (Definition~\ref{def:Y}). It is built from a family of probabilities $\P^{z,h}_B$ on $\calA_\aux^B$, one for each $(z,h)\in\calX$, $k\ge0$ and $B\in\calP_k(h)$. The properties of these probabilities which we use are stated in Proposition~\ref{prop:blocklaws}, which we prove in Section~\ref{sec:construction}. 

To state property \ref{P:freq} of Proposition~\ref{prop:blocklaws} we need the following notation. Let $s\ge0$.
\begin{itemize}[leftmargin=*]
\item $\mathcal K_s$ is the partition of $\calX\times\calA_\aux^\Z$ into the nonempty sets
$$
\bigl\{(z,h,\omega)\ ;\ z[-s,s]=u,\ (h_{i,k})_{i\in[-s,s],\,0\le k\le s}=v,\ \omega[-s,s]=w\bigr\},
$$
where $u\in\{-1,1\}^{[-s,s]}$, $v\in\{0,1\}^{[-s,s]\times[0,s]}$ and $w\in\calA_\aux^{[-s,s]}$. The partitions $\mathcal K_s$ are finite and increase with $s$.

\item For an interval $I=[a,b]=\{a,a+1,\dots,b\}$, recall from \ref{not:intervals} that
$$
I^{(s)}=\bigl\{i\in I\ ;\ [i-s,i+s]\subset I\bigr\}
$$
is the set of points of $I$ at distance greater than $s$ from $\Z\setminus I$. Thus $I^{(s)}=[a+s,b-s]$ if $b-a\ge2s$ and $I^{(s)}=\emptyset$ otherwise.

\item For an interval $I$ and $\Gamma\in\mathcal K_s$, the number of occurrences of $\Gamma$ along $I$ is
$$
S_{\Gamma,I}(z,h,\omega):=\sum_{i\in I^{(s)}}\1_\Gamma\bigl(\widehat\sigma^i(z,h,\omega)\bigr),\quad\text{where }(z,h,\omega)\in\calX\times\calA_\aux^\Z ,
$$
and $\widehat\sigma(z,h,\omega)=(\sigma_\calX(z,h),\sigma\omega)$ is the shift of $\calX\times\calA_\aux^\Z$. The number $S_{\Gamma,I}(z,h,\omega)$ depends only on $z[I]$, on $(h_{i,k})_{i\in I,\,0\le k\le s}$ and on $\omega[I]$. For fixed $(z,h)$ we regard it as a function on $\calA_\aux^I$.
\end{itemize}

The next proposition gives, for every base point $(z,h)\in\calX$ and every block $B$, a law $\P^{z,h}_B$ for the auxiliary coordinates along $B$, and it collects the properties of these laws which we use to construct $\mu_0$ (Definition~\ref{def:mu0}) and to prove its properties.

\begin{proposition}
\label{prop:blocklaws}
There exist a sequence $(\calX_k)_{k\ge0}$ of Borel subsets of $\calX$ and a family of probabilities $\P^{z,h}_B$ on $\calA_\aux^B$, one for each $(z,h)\in\calX$, $k\ge0$ and $B\in\calP_k(h)$, with the following properties:
\begin{enumerate}[label=\textup{(P\arabic*)},leftmargin=*]
\item\label{P:unif} If $k=0$, then $B=\{i\}$ for some $i\in\Z$, and $\P^{z,h}_B$ is the uniform probability on $\calA_\aux^{\{i\}}$.
\item\label{P:freq} For $0\le s<k$, $\Gamma\in\mathcal K_s$ and $w\in\supp\P^{z,h}_B$,
$$
\bigl|S_{\Gamma,B}(z,h,w)-\P^{z,h}_B\bigl(S_{\Gamma,B}(z,h,\cdot)\bigr)\bigr|\le3\cdot2^{s+1-k}|B| .
$$
\item\label{P:marg} Let $k\ge1$, and let $B_1,\dots,B_n$ be the elements of $\calP_{k-1}(h)$ contained in $B$, ordered from left to right. By Lemma~\ref{lem:markers}\ref{mk:children}, $B=B_1\cup\dots\cup B_n$. For every $T\subset[1,n]$ with $\#T\le2R$ and all $u_j\in\calA_\aux^{B_j}$, where $j\in T$,
$$
\P^{z,h}_B\bigl(\bigl\{w\in\calA_\aux^B\ ;\ w[B_j]=u_j\text{ for all }j\in T\bigr\}\bigr)=\prod_{j\in T}\P^{z,h}_{B_j}\bigl(\{u_j\}\bigr).
$$
\item\label{P:entropy} Let $k\ge1$ and $B=B_1\cup\dots\cup B_n$ as in \ref{P:marg}, and for $F\subset B$ let $\scrP_F$ be the partition of $\calA_\aux^B$ into the cylinders over $F$ (see Section~\ref{sec:entropy-def}). If $\sigma_\calX^{\min B}(z,h)\in\calX_k$, then
$$
\sum_{j=1}^nH_{\P^{z,h}_B}\bigl(\scrP_{B_j}\bigr)-H_{\P^{z,h}_B}\bigl(\scrP_B\bigr)\le\log\frac43 .
$$
If $\sigma_\calX^{\min B}(z,h)\notin\calX_k$, then $\#\supp\P^{z,h}_B\le2^{9Rq|B|/(2L)+12L}$.
\item\label{P:cov} Let $(z',h')\in\calX$ and $c\in\Z$ with $B+c\in\calP_k(h')$, and assume that $z'_{i+c}=z_i$ and $h'_{i+c,k'}=h_{i,k'}$ for all $i\in B$ and $0\le k'<k$. Then $$\P^{z',h'}_{B+c}\ \text{is the pushforward of }\P^{z,h}_B\text{ under the map }w\mapsto(w_{i-c})_{i\in B+c}.$$
\item\label{P:typ} For $k\ge0$, $B\in\calP_k(h)$ and $i\in B$, we have $\sigma_\calX^i(z,h)\in\calX_k$ exactly when $\sigma_\calX^{\min B}(z,h)\in\calX_k$. Moreover, for the probability $\zeta=\nu\times\P_\calH$ on $\calX$ of Section~\ref{sec:base},
$$
\sum_{k\ge1}\zeta(\calX\setminus\calX_k)\le\frac18,\ \text{and }\liminf_{k\to\infty}\eta(\calX_k)=0$$ 
for every ergodic $\sigma_\calX$-invariant Borel probability $\eta$ on $\calX$ with $\eta\ne\zeta$.
\end{enumerate}
\end{proposition}

In particular, by \ref{P:cov} with $c=0$, the probability $\P^{z,h}_B$ depends only on $B$, on $z[B]$ and on $(h_{i,k'})_{i\in B,\,0\le k'<k}$. Property \ref{P:marg} gives the mixing of order $R$ (Theorem~\ref{thm:mixing}, through Lemma~\ref{lem:sample}), property \ref{P:freq} gives the uniqueness of the lift (Theorem~\ref{thm:faithful}), and properties \ref{P:entropy} and \ref{P:typ} give the entropy gap (Propositions~\ref{prop:lower} and~\ref{prop:upper}).

From now on we fix the sets $\calX_k$ and the probabilities $\P^{z,h}_B$ given by Proposition~\ref{prop:blocklaws}.

\begin{definition}
\label{def:Y}
Let
$$
\widehat\pi:(z,h,\omega)\in \calX\times\calA_\aux^\Z\mapsto (z,\omega)\in \calA^\Z .
$$
We define
$$
\widehat Y:=\bigl\{(z,h,\omega)\in\calX\times\calA_\aux^\Z\ ;\ \omega[B]\in\supp\P^{z,h}_B\text{ for all }k\ge0\text{ and }B\in\calP_k(h)\bigr\}
$$
and $Y:=\widehat\pi(\widehat Y)$.
\end{definition}

In other words, $(z,\omega)\in\calA^\Z$ lies in $Y$ exactly when there is $h\in\Sigma_\calH$ such that, for every $k\ge0$ and every $B\in\calP_k(h)$, the word $\omega[B]$ belongs to $\supp\P^{z,h}_B$. In Lemma~\ref{lem:Yhat} we show that $Y$ is a subshift.

\subsection{The sample measures and the measure \texorpdfstring{$\mu_0$}{mu0}}
\label{sec:sample}
We recall that Proposition \ref{prop:blocklaws} is being assumed in this subsection. For $(z,h)\in\calX$ and $k\ge0$ we define the product probability
$$
\mu^{(k)}_{z,h}:=\bigotimes_{B\in\calP_k(h)}\P^{z,h}_B\quad\text{on}\quad\calA_\aux^\Z=\prod_{B\in\calP_k(h)}\calA_\aux^B ,
$$
under which the restrictions $\omega[B]$, $B\in\calP_k(h)$, are independent with laws $\P^{z,h}_B$. This is the analogue of the measure $\nu^{(k)}_h$ of \cite[(34)]{BassolsCastro2026}, with the auxiliary coordinates in place of the coordinates of $\Sig$ and with the laws $\P^{z,h}_B$ in place of the laws $\rho_B$. The next lemma shows that the joint law of at most $2R$ blocks does not change when $k$ increases.

\begin{lemma}
\label{lem:marginals}
Let $(z,h)\in\calX$, $k\ge0$, and let $\calF\subset\calP_k(h)$ with $\#\calF\le2R$. Then for every $k'\ge k$ and all $u_B\in\calA_\aux^B$, where $B\in\calF$,
$$
\mu^{(k')}_{z,h}\bigl(\bigl\{\omega\in\calA_\aux^\Z\ ;\ \omega[B]=u_B\text{ for all }B\in\calF\bigr\}\bigr)=\prod_{B\in\calF}\P^{z,h}_B\bigl(\{u_B\}\bigr).
$$
\end{lemma}

\begin{proof}
For $k'\ge k$ let
$$
\calF_{k'}:=\bigl\{B'\in\calP_{k'}(h)\ ;\ B\subset B'\text{ for some }B\in\calF\bigr\}.
$$
By Lemma~\ref{lem:markers}\ref{mk:length} and~\ref{mk:children}, $\calP_{k'}(h)$ is a partition of $\Z$ and every element of $\calP_{k'}(h)$ is a union of elements of $\calP_k(h)$, so every element of $\calP_k(h)$ is contained in exactly one element of $\calP_{k'}(h)$.  Hence $\calF_k=\calF$ and $\#\calF_{k'}\le\#\calF\le2R$. By Lemma~\ref{lem:markers}\ref{mk:children}, every $B'\in\calF_{k'}$ is contained in exactly one $J\in\calF_{k'+1}$, and for $J\in\calF_{k'+1}$ there are at most $2R$ elements $B'\in\calF_{k'}$ with $B'\subset J$.

Let $k'\ge k$ and $u_{B'}\in\calA_\aux^{B'}$, where $B'\in\calF_{k'}$, and let
$$
C:=\bigl\{\omega\in\calA_\aux^\Z\ ;\ \omega[B']=u_{B'}\text{ for all }B'\in\calF_{k'}\bigr\}.
$$
For $J\in\calF_{k'+1}$ let
$$
C_J:=\bigl\{w\in\calA_\aux^J\ ;\ w[B']=u_{B'}\text{ for all }B'\in\calF_{k'}\text{ with }B'\subset J\bigr\}.
$$
Under $\mu^{(k'+1)}_{z,h}$ the restrictions $\omega[J]$, $J\in\calP_{k'+1}(h)$, are independent with laws $\P^{z,h}_J$. Hence, by Proposition~\ref{prop:blocklaws}\ref{P:marg} and the definition of $\mu^{(k')}_{z,h}$,
$$
\mu^{(k'+1)}_{z,h}(C)=\prod_{J\in\calF_{k'+1}}\P^{z,h}_J(C_J)=\prod_{J\in\calF_{k'+1}}\ \prod_{\substack{B'\in\calF_{k'}\\ B'\subset J}}\P^{z,h}_{B'}\bigl(\{u_{B'}\}\bigr)=\prod_{B'\in\calF_{k'}}\P^{z,h}_{B'}\bigl(\{u_{B'}\}\bigr)=\mu^{(k')}_{z,h}(C).
$$
Since every element of $\calF$ is contained in an element of $\calF_{k'}$, the set $\{\omega\ ;\ \omega[B]=u_B\text{ for all }B\in\calF\}$ is a finite disjoint union of sets of the form $\{\omega\ ;\ \omega[B']=u_{B'}\text{ for all }B'\in\calF_{k'}\}$. Hence $\mu^{(k'+1)}_{z,h}$ and $\mu^{(k')}_{z,h}$ give the same probability to it. For $k'=k$ this probability is $\prod_{B\in\calF}\P^{z,h}_B(\{u_B\})$ by the definition of $\mu^{(k)}_{z,h}$, and induction on $k'$ gives the lemma.
\end{proof}

By Lemma~\ref{lem:marginals}, the probability under $\mu^{(k)}_{z,h}$ of a cylinder over a fixed interval does not depend on $k$ for $k$ large. The next lemma defines the sample measure $\mu_{z,h}$ as the limit of the probabilities $\mu^{(k)}_{z,h}$. 

\begin{lemma}
\label{lem:sample}
Let $(z,h)\in\calX$. There is a unique probability $\mu_{z,h}$ on $\calA_\aux^\Z$ such that, for every $k\ge0$, every $\calF\subset\calP_k(h)$ with $\#\calF\le2R$ and all $u_B\in\calA_\aux^B$, where $B\in\calF$,
\begin{equation}
\label{eq:sample}
\mu_{z,h}\bigl(\bigl\{\omega\in\calA_\aux^\Z\ ;\ \omega[B]=u_B\text{ for all }B\in\calF\bigr\}\bigr)=\prod_{B\in\calF}\P^{z,h}_B\bigl(\{u_B\}\bigr).
\end{equation}
Moreover the following hold.
\begin{enumerate}[label=\textup{(\roman*)},leftmargin=*]
\item\label{sm:stab} For $n\ge0$ and $k\ge0$ with $(4L)^k\ge2n+1$, and for every $w\in\calA_\aux^{[-n,n]}$,
$$
\mu_{z,h}\bigl(\{\omega\ ;\ \omega[-n,n]=w\}\bigr)=\mu^{(k)}_{z,h}\bigl(\{\omega\ ;\ \omega[-n,n]=w\}\bigr).
$$
\item\label{sm:cont} For $n,k\ge0$ and $w\in\calA_\aux^{[-n,n]}$, the maps $$(z,h)\mapsto\mu_{z,h}(\omega[-n,n]=w)\ \text{and}\ (z,h)\mapsto\mu^{(k)}_{z,h}(\omega[-n,n]=w)$$ 
are locally constant on $\calX$.
\item\label{sm:cov} $\mu_{\sigma_\calX(z,h)}=\sigma_*\mu_{z,h}$, and $\mu^{(k)}_{\sigma_\calX(z,h)}=\sigma_*\mu^{(k)}_{z,h}$ for every $k\ge0$.
\item\label{sm:supp} $\mu_{z,h}(\widehat Y_{z,h})=1$, where $\widehat Y_{z,h}:=\bigl\{\omega\in\calA_\aux^\Z\ ;\ (z,h,\omega)\in\widehat Y\bigr\}$ and $\widehat Y$ is given by Definition~\ref{def:Y}.
\end{enumerate}
\end{lemma}

\begin{proof}
We start with the existence and \ref{sm:stab}. Let $n\ge0$ and $k\ge0$ with $(4L)^k\ge2n+1$. If $k\ge1$, the interval $[-n,n]$ is shorter than every element of $\calP_k(h)$ by Lemma~\ref{lem:markers}\ref{mk:length}, so it intersects at most two of them, and if $k=0$, then $n=0$ and $[-n,n]$ intersects one element. 

Let $\calF$ be the set of elements of $\calP_k(h)$ which intersect $[-n,n]$. For $w\in\calA_\aux^{[-n,n]}$, the set $\{\omega\ ;\ \omega[-n,n]=w\}$ is a finite disjoint union of sets $\{\omega\ ;\ \omega[B]=u_B\text{ for all }B\in\calF\}$, so by Lemma~\ref{lem:marginals} the number $\mu^{(k')}_{z,h}(\{\omega\ ;\ \omega[-n,n]=w\})$ is the same for all $k'\ge k$. These numbers, for $n\ge0$ and $w\in\calA_\aux^{[-n,n]}$, form a consistent family, and the Kolmogorov extension theorem gives a probability $\mu_{z,h}$ on $\calA_\aux^\Z$ which satisfies \ref{sm:stab}.

Let $k\ge0$, $\calF\subset\calP_k(h)$ with $\#\calF\le2R$ and $u_B\in\calA_\aux^B$ for $B\in\calF$, choose $n$ with $$\bigcup\calF\subset[-n,n]\text{ and }k'\ge k\text{ with }(4L)^{k'}\ge2n+1.$$ The set $\{\omega\ ;\ \omega[B]=u_B\text{ for all }B\in\calF\}$ is a finite disjoint union of sets $\{\omega\ ;\ \omega[-n,n]=w\}$, so \eqref{eq:sample} follows from \ref{sm:stab} for $k'$ and Lemma~\ref{lem:marginals}.

Uniqueness holds because, for $(4L)^k\ge2n+1$, the set $\{\omega\ ;\ \omega[-n,n]=w\}$ is a finite disjoint union of sets $\{\omega\ ;\ \omega[B]=u_B\text{ for all }B\in\calF\}$, where $\calF$ is the set of the at most two elements of $\calP_k(h)$ which intersect $[-n,n]$, so its probability is determined by \eqref{eq:sample}.

We now prove \ref{sm:cont}. Let $J_n:=[-n-3(4L)^k,n+3(4L)^k]$. By Lemma~\ref{lem:markers}\ref{mk:length}, the elements of $\calP_k(h)$ which intersect $[-n,n]$ are contained in $J_n$, and they are determined by $(h_{i,k})_{i\in J_n}$, since the largest point of $\calM_k(h)$ in $(-\infty,-n]$ and the smallest point of $\calM_k(h)$ in $(n,\infty)$ lie in $J_n$. 

Let $(z',h')\in\calX$ with $z'[J_n]=z[J_n]$ and $h'_{i,k'}=h_{i,k'}$ for $i\in J_n$ and $k'\le k$. Then $\calP_k(h)$ and $\calP_k(h')$ have the same elements intersecting $[-n,n]$, and by Proposition~\ref{prop:blocklaws}\ref{P:cov} with $c=0$ these blocks have the same laws. Hence $$\mu^{(k)}_{z',h'}(\omega[-n,n]=w)=\mu^{(k)}_{z,h}(\omega[-n,n]=w)$$ for every $(z',h')$ in the element of $\scrQ_{n+3(4L)^k}$ which contains $(z,h)$. The statement for $\mu_{z,h}$ follows from this and \ref{sm:stab}.

We now prove \ref{sm:cov}. Since $\calM_k(\sigma_\calH h)=\calM_k(h)-1$, the elements of $\calP_k(\sigma_\calH h)$ are the intervals $B-1$, $B\in\calP_k(h)$. By Proposition~\ref{prop:blocklaws}\ref{P:cov} with $(z',h')=\sigma_\calX(z,h)$ and $c=-1$, the probability $$\P^{\sigma_\calX(z,h)}_{B-1}\text{ is the pushforward of }\P^{z,h}_B\text{ under the map }\tau_{B,-1},$$ that is, under the map $w\mapsto w'$, where $w'_{i-1}:=w_i$. Since $(\sigma\omega)_{i-1}=\omega_i$, this gives $$\mu^{(k)}_{\sigma_\calX(z,h)}=\sigma_*\mu^{(k)}_{z,h}.$$ For $n\ge0$ and $k$ with $(4L)^k\ge2n+3$, \ref{sm:stab} shows that $\mu_{\sigma_\calX(z,h)}$ and $\sigma_*\mu_{z,h}$ give the same probabilities to the sets $\{\omega\ ;\ \omega[-n,n]=w\}$ as $\mu^{(k)}_{\sigma_\calX(z,h)}$ and $\sigma_*\mu^{(k)}_{z,h}$, respectively, so $\mu_{\sigma_\calX(z,h)}=\sigma_*\mu_{z,h}$.

Finally we prove \ref{sm:supp}. For $k\ge0$ and $B\in\calP_k(h)$, \eqref{eq:sample} with $\calF=\{B\}$ gives $$\mu_{z,h}(\{\omega\ ;\ \omega[B]=u\})=\P^{z,h}_B(\{u\})$$ for every $u\in\calA_\aux^B$, so $\mu_{z,h}(\{\omega\ ;\ \omega[B]\in\supp\P^{z,h}_B\})=1$. There are countably many pairs $(k,B)$, which gives \ref{sm:supp}.
\end{proof}

We call the probabilities $\mu_{z,h}$ the \emph{sample measures}. We now define the distinguished measure.

\begin{definition}
\label{def:mu0}
Let
$$
\pi_\calX:(z,h,\omega)\in \calX\times\calA_\aux^\Z\mapsto (z,h)\in \calX .
$$
We define the probabilities $\widehat\mu_0$ on $\calX\times\calA_\aux^\Z$ and $\mu_0$ on $\calA^\Z$ by
\begin{equation}
\label{eq:mu0}
\widehat\mu_0(dz,dh,d\omega):=\mu_{z,h}(d\omega)\,\zeta(dz,dh)\quad\text{and}\quad\mu_0:=\widehat\pi_*\widehat\mu_0 ,
\end{equation}
where $\zeta=\nu\times\P_\calH$ is the probability on $\calX$ of Section~\ref{sec:base}.
\end{definition}

Explicitly, $$\widehat\mu_0(F)=\int_\calX\mu_{z,h}\bigl(F(z,h,\cdot)\bigr)\,d\zeta(z,h)$$ for every bounded Borel $F$. The integrand is a Borel function of $(z,h)$, by Lemma~\ref{lem:sample}\ref{sm:cont} for indicators of the elements of the partitions $\mathcal K_s$ and by the monotone class theorem in general.

The next lemma shows that Definitions~\ref{def:Y} and~\ref{def:mu0} give a subshift $Y$ and an invariant probability $\mu_0$ on $Y$ whose first coordinate has law $\nu$. 
\begin{lemma}
\label{lem:Yhat}
Let $\widehat Y$ and $Y$ be given by Definition~\ref{def:Y}, and $\widehat\mu_0$ and $\mu_0$ by \eqref{eq:mu0}. The following hold.
\begin{enumerate}[label=\textup{(\roman*)},leftmargin=*]
\item\label{Y:compact} $\widehat Y$ is compact, $\widehat\sigma(\widehat Y)=\widehat Y$ and $\pi_\calX(\widehat Y)=\calX$.
\item\label{Y:sub} $Y$ is a subshift of $\calA^\Z$, and $\pi_1(Y)=\Sig$.
\item\label{Y:mu} $\widehat\mu_0$ is $\widehat\sigma$-invariant, $\widehat\mu_0(\widehat Y)=1$ and $(\pi_\calX)_*\widehat\mu_0=\zeta$. Hence, $\mu_0$ is $\sigma$-invariant, $\mu_0(Y)=1$ and $(\pi_1)_*\mu_0=\nu$.
\end{enumerate}
\end{lemma}

\begin{proof}
We start by proving \ref{Y:compact}. Let $(z,h,\omega)\notin\widehat Y$. Then there are $k\ge0$ and $B=[u,v-1]\in\calP_k(h)$ with $\omega[B]\notin\supp\P^{z,h}_B$. Let
$$
U:=\left\{(z',h',\omega')\in\calX\times\calA_\aux^\Z\ ;\ \begin{array}{l}z'[B]=z[B],\ \omega'[B]=\omega[B],\\ h'_{i,k'}=h_{i,k'}\text{ for }(i,k')\in[u,v]\times[0,k]\end{array}\right\}.
$$
The set $U$ is open, since it is defined by finitely many coordinates, and $(z,h,\omega)\in U$. Let $(z',h',\omega')\in U$. Since $h'_{i,k}=h_{i,k}$ for $i\in[u,v]$, the points $u<v$ are consecutive points of $\calM_k(h')$, so $B\in\calP_k(h')$. Since moreover $z'[B]=z[B]$ and $h'_{i,k'}=h_{i,k'}$ for $i\in B$ and $0\le k'<k$, Proposition~\ref{prop:blocklaws}\ref{P:cov} with $c=0$ gives $\P^{z',h'}_B=\P^{z,h}_B$. Hence $\omega'[B]=\omega[B]\notin\supp\P^{z',h'}_B$, so $(z',h',\omega')\notin\widehat Y$. Thus every point outside $\widehat Y$ has an open neighbourhood disjoint from $\widehat Y$, so $\widehat Y$ is closed in the compact space $\calX\times\calA_\aux^\Z$.  

We now prove \ref{Y:sub}. The map $\widehat\pi$ is continuous and $\widehat\pi\circ\widehat\sigma=\sigma\circ\widehat\pi$, so $Y$ is compact and $\sigma(Y)=Y$ by \ref{Y:compact}. Since $\Sigma_\calH\ne\emptyset$ and $\pi_\calX(\widehat Y)=\calX$, we have $\pi_1(Y)=\Sig$.

Finally we prove \ref{Y:mu}. For bounded Borel $F$, by Lemma~\ref{lem:sample}\ref{sm:cov} and the invariance of $\zeta$ (see Lemma~\ref{lem:base}\ref{base:mixing}),
$$
\widehat\mu_0(F\circ\widehat\sigma)=\int_\calX(\sigma_*\mu_{z,h})\bigl(F(\sigma_\calX(z,h),\cdot)\bigr)\,d\zeta=\int_\calX\mu_{\sigma_\calX(z,h)}\bigl(F(\sigma_\calX(z,h),\cdot)\bigr)\,d\zeta=\widehat\mu_0(F).
$$
By Lemma~\ref{lem:sample}\ref{sm:supp}, $\widehat\mu_0(\widehat Y)=\int\mu_{z,h}(\widehat Y_{z,h})\,d\zeta=1$, and $(\pi_\calX)_*\widehat\mu_0=\zeta$ by \eqref{eq:mu0}. The statements on $\mu_0$ follow, since $\widehat\pi\circ\widehat\sigma=\sigma\circ\widehat\pi$, $\widehat\pi(\widehat Y)=Y$, and $\pi_1\circ\widehat\pi(z,h,\omega)=z$ while the pushforward of $\zeta$ under the map $(z,h)\mapsto z$ is $\nu$.
\end{proof}

\label{sec:faithful}
We recall that Proposition~\ref{prop:blocklaws} is being assumed in this subsection. We now prove that every invariant probability of $\calX$ has exactly one invariant lift to $\widehat Y$, that is, the extension $\pi_\calX:\widehat Y\to\calX$ is faithful in the sense of Downarowicz and Huczek \cite{DownarowiczHuczek2012}.

\begin{theorem}
\label{thm:faithful}
Let $\eta$ be a $\sigma_\calX$-invariant Borel probability on $\calX$, and let $\widehat\rho$ be a $\widehat\sigma$-invariant Borel probability on $\calX\times\calA_\aux^\Z$ with $\widehat\rho(\widehat Y)=1$ and $(\pi_\calX)_*\widehat\rho=\eta$. Then
$$
\widehat\rho(dz,dh,d\omega)=\mu_{z,h}(d\omega)\,\eta(dz,dh).
$$
\end{theorem}

\begin{proof}
For $s\ge0$ and $\Gamma\in\mathcal K_s$ define $g_\Gamma(z,h):=\mu_{z,h}\bigl(\{\omega\ ;\ (z,h,\omega)\in\Gamma\}\bigr)$, a locally constant function on $\calX$ by Lemma~\ref{lem:sample}\ref{sm:cont}. We divide the proof into three steps.

\begin{step}[1]
We show that, for all $s\ge0$, $\Gamma\in\mathcal K_s$, $k>s$, $(z,h,\omega)\in\widehat Y$ and $B\in\calP_k(h)$,
$$
\Bigl|S_{\Gamma,B}(z,h,\omega)-\sum_{i\in B^{(s)}}g_\Gamma\bigl(\sigma_\calX^i(z,h)\bigr)\Bigr|\le3\cdot2^{s+1-k}|B| .
$$
\end{step}

\begin{stepproof}
By \eqref{eq:sample} with $\calF=\{B\}$, we have $\mu_{z,h}(\{\omega\ ;\ \omega[B]=u\})=\P^{z,h}_B(\{u\})$ for every $u\in\calA_\aux^B$. For $i\in B^{(s)}$ the event $\{\widehat\sigma^i(z,h,\omega)\in\Gamma\}$ depends only on $\omega[i-s,i+s]$, and $[i-s,i+s]\subset B$, so by Lemma~\ref{lem:sample}\ref{sm:cov},
$$
\begin{aligned}
\P^{z,h}_B\bigl(S_{\Gamma,B}(z,h,\cdot)\bigr)&=\sum_{i\in B^{(s)}}\mu_{z,h}\bigl(\{\omega\ ;\ (\sigma_\calX^i(z,h),\sigma^i\omega)\in\Gamma\}\bigr)\\
&=\sum_{i\in B^{(s)}}\mu_{\sigma_\calX^i(z,h)}\bigl(\{\omega\ ;\ (\sigma_\calX^i(z,h),\omega)\in\Gamma\}\bigr)=\sum_{i\in B^{(s)}}g_\Gamma\bigl(\sigma_\calX^i(z,h)\bigr).
\end{aligned}
$$
Since $\omega\in\widehat Y_{z,h}$, we have $\omega[B]\in\supp\P^{z,h}_B$, and Proposition~\ref{prop:blocklaws}\ref{P:freq} gives the claim.
\end{stepproof}

\begin{step}[2]
We show that, for all $s\ge0$, $\Gamma\in\mathcal K_s$, $k>s$, $n\ge1$ and $(z,h,\omega)\in\widehat Y$,
$$
\Bigl|\sum_{i=0}^{n-1}\Bigl(\1_\Gamma\bigl(\widehat\sigma^i(z,h,\omega)\bigr)-g_\Gamma\bigl(\sigma_\calX^i(z,h)\bigr)\Bigr)\Bigr|\le3\cdot2^{s+1-k}n+2s\,(4L)^{-k}n+6(4L)^k .
$$
\end{step}

\begin{stepproof}
Let
$$
\calF:=\bigl\{B\in\calP_k(h)\ ;\ B\subset[0,n-1]\bigr\}.
$$
Since the elements of $\calF$ are disjoint subsets of $[0,n-1]$, the interval $[0,n-1]$ is the disjoint union of the sets $B^{(s)}$ and $B\setminus B^{(s)}$, where $B\in\calF$, and of the set $[0,n-1]\setminus\bigcup\calF$. We bound the sum in the claim over each of these sets, using that every term $\1_\Gamma(\widehat\sigma^i(z,h,\omega))-g_\Gamma(\sigma_\calX^i(z,h))$ has absolute value at most $1$.

Let $B\in\calF$. Since $S_{\Gamma,B}(z,h,\omega)=\sum_{i\in B^{(s)}}\1_\Gamma(\widehat\sigma^i(z,h,\omega))$, Step~1 gives
$$
\Bigl|\sum_{i\in B^{(s)}}\Bigl(\1_\Gamma\bigl(\widehat\sigma^i(z,h,\omega)\bigr)-g_\Gamma\bigl(\sigma_\calX^i(z,h)\bigr)\Bigr)\Bigr|\le3\cdot2^{s+1-k}|B| .
$$
Moreover, $B\setminus B^{(s)}$ has at most $2s$ elements by \eqref{eq:interior}.

Let $i\in[0,n-1]\setminus\bigcup\calF$, and let $B$ be the element of $\calP_k(h)$ which contains $i$. The interval $B$ contains $i$ and is not contained in $[0,n-1]$, so it contains $0$ or $n-1$. By Lemma~\ref{lem:markers}\ref{mk:length} the elements of $\calP_k(h)$ which contain $0$ or $n-1$ have length at most $3(4L)^k$, so
$$
\#\Bigl([0,n-1]\setminus\bigcup\calF\Bigr)\le6(4L)^k .
$$

Finally, every element of $\calP_k(h)$ has length at least $(4L)^k$ by Lemma~\ref{lem:markers}\ref{mk:length}, so $\#\calF\le n/(4L)^k$, and $\sum_{B\in\calF}|B|\le n$. Hence, by the triangle inequality and the bounds above,
$$
\begin{aligned}
&\Bigl|\sum_{i=0}^{n-1}\Bigl(\1_\Gamma\bigl(\widehat\sigma^i(z,h,\omega)\bigr)-g_\Gamma\bigl(\sigma_\calX^i(z,h)\bigr)\Bigr)\Bigr|\\
\le&\sum_{B\in\calF}\Bigl|\sum_{i\in B^{(s)}}\Bigl(\1_\Gamma\bigl(\widehat\sigma^i(z,h,\omega)\bigr)-g_\Gamma\bigl(\sigma_\calX^i(z,h)\bigr)\Bigr)\Bigr|+\sum_{B\in\calF}\#\bigl(B\setminus B^{(s)}\bigr)+\#\Bigl([0,n-1]\setminus\bigcup\calF\Bigr)\\
\le & \sum_{B\in\calF}3\cdot2^{s+1-k}|B|+2s\,\#\calF+6(4L)^k\le  3\cdot2^{s+1-k}n+2s\,(4L)^{-k}n+6(4L)^k .
\end{aligned}
$$
\end{stepproof}

\begin{step}[3]
Finally we prove that $\widehat\rho(dz,dh,d\omega)=\mu_{z,h}(d\omega)\,\eta(dz,dh)$.
\end{step}

\begin{stepproof}
Let $s\ge0$ and $\Gamma\in\mathcal K_s$. Since $\widehat\rho$ is $\widehat\sigma$-invariant, $\widehat\rho(\widehat Y)=1$ and $(\pi_\calX)_*\widehat\rho=\eta$ is $\sigma_\calX$-invariant, Step~2 gives, for all $k>s$ and $n\ge1$,
\begin{align*}
\bigl|\widehat\rho(\Gamma)-\eta(g_\Gamma)\bigr|&=\frac1n\Bigl|\int\sum_{i=0}^{n-1}\Bigl(\1_\Gamma\bigl(\widehat\sigma^i(z,h,\omega)\bigr)-g_\Gamma\bigl(\sigma_\calX^i(z,h)\bigr)\Bigr)\,d\widehat\rho(z,h,\omega)\Bigr|\\
&\le3\cdot2^{s+1-k}+2s\,(4L)^{-k}+\frac{6(4L)^k}n .
\end{align*}
Letting $n\to\infty$ and then $k\to\infty$ gives
$$
\widehat\rho(\Gamma)=\eta(g_\Gamma)=\int_\calX\mu_{z,h}\bigl(\{\omega\ ;\ (z,h,\omega)\in\Gamma\}\bigr)\,d\eta(z,h).
$$
Since the partitions $\mathcal K_s$ increase with $s$, the family $\{\emptyset\}\cup\bigcup_{s\ge0}\mathcal K_s$ is closed under intersections, which generates the Borel $\sigma$-algebra of $\calX\times\calA_\aux^\Z$, and the probabilities $\widehat\rho$ and $\mu_{z,h}(d\omega)\,\eta(dz,dh)$ coincide.
\end{stepproof}
\end{proof}

\begin{corollary}
\label{cor:faithful}
The measure $\widehat\mu_0$ is the only $\widehat\sigma$-invariant Borel probability $\widehat\rho$ on $\calX\times\calA_\aux^\Z$ with $\widehat\rho(\widehat Y)=1$ and $(\pi_\calX)_*\widehat\rho=\zeta$, where $\zeta=\nu\times\P_\calH$ (see Section~\ref{sec:base}). Moreover $(\widehat\sigma,\widehat\mu_0)$ and $(\sigma,\mu_0)$ are ergodic.
\end{corollary}

\begin{proof}
We start by proving the first statement. By Lemma~\ref{lem:Yhat}\ref{Y:mu}, the probability $\widehat\mu_0$ is $\widehat\sigma$-invariant, $\widehat\mu_0(\widehat Y)=1$ and $(\pi_\calX)_*\widehat\mu_0=\zeta$, so $\widehat\mu_0$ has the stated properties. Conversely, let $\widehat\rho$ be a $\widehat\sigma$-invariant Borel probability on $\calX\times\calA_\aux^\Z$ with $\widehat\rho(\widehat Y)=1$ and $(\pi_\calX)_*\widehat\rho=\zeta$. The probability $\zeta$ is $\sigma_\calX$-invariant by Lemma~\ref{lem:base}\ref{base:mixing}, so Theorem~\ref{thm:faithful} with $\eta=\zeta$ and \eqref{eq:mu0} give
$$
\widehat\rho(dz,dh,d\omega)=\mu_{z,h}(d\omega)\,\zeta(dz,dh)=\widehat\mu_0(dz,dh,d\omega).
$$

We now prove that $(\widehat\sigma,\widehat\mu_0)$ and $(\sigma,\mu_0)$ are ergodic. Let $F\subset\calX\times\calA_\aux^\Z$ be a $\widehat\sigma$-invariant Borel set with $\widehat\mu_0(F)>0$, and define $\widehat\rho:=\widehat\mu_0(\,\cdot\cap F)/\widehat\mu_0(F)$. Then $\widehat\rho$ is $\widehat\sigma$-invariant, $\widehat\rho(\widehat Y)=1$, and $\eta:=(\pi_\calX)_*\widehat\rho$ is $\sigma_\calX$-invariant with $\eta\le\zeta/\widehat\mu_0(F)$. The density $d\eta/d\zeta$ is $\sigma_\calX$-invariant $\zeta$-almost everywhere, so it is constant since $(\sigma_\calX,\zeta)$ is ergodic (see Lemma~\ref{lem:base}\ref{base:mixing}), and $\eta=\zeta$. By the first statement $\widehat\rho=\widehat\mu_0$, so $\widehat\mu_0(F)=\widehat\rho(F)=1$. Hence $(\widehat\sigma,\widehat\mu_0)$ is ergodic, and so is its factor $(\sigma,\mu_0)$, since $\mu_0=\widehat\pi_*\widehat\mu_0$.
\end{proof}

The next lemma shows that the shift $(Y,\sigma)$ has no fixed points. This allows us to apply the neutral isolation construction of Lin and Tian \cite[Proposition 3.3]{LinTian2026} (see Proposition~\ref{prop:LT}), which is the key step in the smooth realisation of $(Y,\sigma)$ in Section~\ref{sec:realisation}.

\begin{lemma}
\label{lem:nofixed}
The subshift $Y$ has no fixed point, that is, $\sigma x\ne x$ for every $x\in Y$.
\end{lemma}

\begin{proof}
Assume for a contraction that $x=(z,\omega)\in Y$ satisfies $\sigma x=x$, and let $c:=\omega_0$, so that $\omega_i=c$ for every $i\in\Z$. The set $F:=\{y\in\widehat Y\ ;\ \widehat\pi(y)=x\}$ is nonempty and compact, and $\widehat\sigma(F)=F$ since $\widehat\sigma(\widehat Y)=\widehat Y$ and $\widehat\pi\circ\widehat\sigma=\sigma\circ\widehat\pi$. Hence there is a $\widehat\sigma$-invariant Borel probability $\widehat\rho$ with $\widehat\rho(F)=1$ (see \cite[Corollary 6.9.1]{Walters1982}). Let $\eta:=(\pi_\calX)_*\widehat\rho$. Since $F\subset\widehat Y$ and every point of $F$ has auxiliary coordinate $c$ at time $0$, Theorem~\ref{thm:faithful} gives
$$
1=\widehat\rho\bigl(\{(z',h',\omega')\ ;\ \omega'_0=c\}\bigr)=\int_\calX\mu_{z',h'}\bigl(\{\omega'\ ;\ \omega'_0=c\}\bigr)\,d\eta(z',h') .
$$
By Lemma~\ref{lem:sample} with $\calF=\{\{0\}\}\subset\calP_0(h')$ and by Proposition~\ref{prop:blocklaws}\ref{P:unif}, $$\mu_{z',h'}(\{\omega'\ ;\ \omega'_0=c\})=\P^{z',h'}_{\{0\}}(\{c\})=2^{-q},$$ where $\calA_\aux=\{-1,1\}^q$, so the integrand equals $2^{-q}<1$. This is a contradiction.
\end{proof}

\subsection{Proof of Proposition~\ref{prop:blocklaws}}
\label{sec:construction}
In this subsection we construct the sets $\calX_k$ and the probabilities $\P^{z,h}_B$ of Proposition~\ref{prop:blocklaws}, recursively in $k$, and we prove their properties.

For $k=0$, the set $B$ is a singleton and $\P^{z,h}_B$ is the uniform probability on $\calA_\aux^B$. For $k\ge1$, we obtain $\P^{z,h}_B$ by applying Proposition~\ref{prop:products} to the probabilities $\P^{z,h}_{B_1},\dots,\P^{z,h}_{B_n}$, where $B_1,\dots,B_n\in\calP_{k-1}(h)$ are the elements contained in $B$, ordered from left to right.

To apply Proposition~\ref{prop:products} we need the functions $X_j$, which we build from the vectors $V_{k,I}$ of Definition~\ref{def:V}, and a choice between $\widetilde\Q$ and $\overline\Q$, which is made by the sets $\calX_k$ of Definition~\ref{def:Xk}.

\begin{definition}
\label{def:V}
Let $\mathcal K_{<k}:=\mathcal K_0\cup\dots\cup\mathcal K_{k-1}$ for $k\ge1$ and $\mathcal K_{<0}:=\emptyset$, with $\R^\emptyset:=\{0\}$, and let $\norm\cdot$ be the Euclidean norm on $\R^{\mathcal K_{<k}}$. For $k\ge0$, an interval $I$ and $(z,h,\omega)\in\calX\times\calA_\aux^\Z$, let $V_{k,I}(z,h,\omega)\in\R^{\mathcal K_{<k}}$ be the vector with coordinates
$$
V_{k,I}(z,h,\omega)_\Gamma:=2^{-s-1}S_{\Gamma,I}(z,h,\omega),\quad\text{where }\Gamma\in\mathcal K_{s}\text{ and }0\le s<k ,
$$
and $S_{\Gamma,I}$ is the number of occurrences of $\Gamma$ along $I$ of Section~\ref{sec:blocklaws}.
\end{definition}

Since $\mathcal K_s$ is a partition,
$$
\sum_{\Gamma\in\mathcal K_{s}}S_{\Gamma,I}=\#I^{(s)}\le|I|,
$$
and since the numbers $S_{\Gamma,I}$ are nonnegative,
\begin{equation}
\label{eq:Vnorm}
\Bigl(\sum_{\Gamma\in\mathcal K_{s}}V_{k,I}(z,h,\omega)_\Gamma^2\Bigr)^{1/2}\le\sum_{\Gamma\in\mathcal K_{s}}V_{k,I}(z,h,\omega)_\Gamma\le2^{-s-1}|I|\quad\text{for every }0\le s<k .
\end{equation}

\begin{definition}
\label{def:test}
Let $m\ge1$, and recall from Section~\ref{sec:base} that $\scrQ_s$ is the partition of $\calX$ according to the values of $z[-s,s]$ and of $(h_{i,k})_{i\in[-s,s],\,0\le k\le s}$, and that $\zeta=\nu\times\P_\calH$.
\begin{itemize}[leftmargin=*]
\item For $s\ge0$, $C\in\scrQ_s$ and an interval $I$, the number of occurrences of $C$ along $I$ is
$$
S_{C,I}(z,h):=\sum_{i\in I^{(s)}}\1_C\bigl(\sigma_\calX^i(z,h)\bigr),\quad\text{where }(z,h)\in\calX ,
$$
and $\sigma_\calX(z,h)=(\sigma z,\sigma_\calH h)$ is the shift of $\calX$.
\item For an interval $I$,
\begin{equation}
\label{eq:test}
\mathcal T_m(I):=\left\{(z,h)\in\calX\ ;\ \begin{array}{l}\bigl|S_{C,I}(z,h)-|I|\,\zeta(C)\bigr|\leq {|I|}/m\\ \text{for every }0\le s<m\text{ and every }C\in\scrQ_s\end{array}\right\}
\end{equation}
is the set of points of $\calX$ at which the frequencies along $I$ of the elements of $\scrQ_0,\dots,\scrQ_{m-1}$ are $1/m$-close to their probabilities under $\zeta$.
\item For $k\ge0$, with $J_k(h)$ the element of $\calP_k(h)$ which contains $0$,
$$
T_{m,k}:=\bigl\{(z,h)\in\calX\ ;\ (z,h)\in\mathcal T_m\bigl(J_k(h)\bigr)\bigr\}.
$$
\end{itemize}
\end{definition}

The next lemma is a standard consequence of the ergodic theorem for $\sigma_\calX$ and for $\sigma_\calX^{-1}$, for averages over intervals $[-n_0,n_1]$ which contain $0$ but need not start at $0$. We use it for the intervals $J_k(h)$.

\begin{lemma}
\label{lem:birkhoff}
Let $\eta$ be an ergodic $\sigma_\calX$-invariant Borel probability on $\calX$, and let $C\in\scrQ_s$ for some $s\ge0$. Then there is a Borel set $G\subset\calX$ with $\eta(G)=1$ such that for every $(z,h)\in G$ and every $\varepsilon>0$ there is $N\ge1$ with
$$
\Bigl|\frac1{n_0+n_1+1}\sum_{i=-n_0}^{n_1}\1_C\bigl(\sigma_\calX^i(z,h)\bigr)-\eta(C)\Bigr|\le\varepsilon\quad\text{for all integers }n_0,n_1\ge0\text{ with }n_0+n_1\ge N .
$$
\end{lemma}

\begin{proof}
By Birkhoff's ergodic theorem applied to the maps $\sigma_\calX$ and $\sigma_\calX^{-1}$, there is a Borel set $G\subset\calX$ with $\eta(G)=1$ such that for every $(z,h)\in G$ and every $\varepsilon>0$ there is $n_\varepsilon$ with
$$
\Bigl|\sum_{i=0}^{n-1}\Bigl(\1_C\bigl(\sigma_\calX^i(z,h)\bigr)-\eta(C)\Bigr)\Bigr|\le\frac\varepsilon2\,n\quad\text{and}\quad\Bigl|\sum_{i=1}^{n}\Bigl(\1_C\bigl(\sigma_\calX^{-i}(z,h)\bigr)-\eta(C)\Bigr)\Bigr|\le\frac\varepsilon2\,n\quad\text{for }n\ge n_\varepsilon .
$$
Fix such $(z,h)$ and $\varepsilon$, and let $n_0,n_1\ge0$. We split the sum at $0$,
$$
\sum_{i=-n_0}^{n_1}\Bigl(\1_C\bigl(\sigma_\calX^i(z,h)\bigr)-\eta(C)\Bigr)=\sum_{i=0}^{n_1}\Bigl(\1_C\bigl(\sigma_\calX^i(z,h)\bigr)-\eta(C)\Bigr)+\sum_{i=1}^{n_0}\Bigl(\1_C\bigl(\sigma_\calX^{-i}(z,h)\bigr)-\eta(C)\Bigr).
$$
Every term has absolute value at most $1$, so each of the two sums on the right is bounded by its number of terms. If $n_1+1\ge n_\varepsilon$ we use the first estimate with $n=n_1+1$, and otherwise we bound the first sum by its number of terms, so
$$
\Bigl|\sum_{i=0}^{n_1}\Bigl(\1_C\bigl(\sigma_\calX^i(z,h)\bigr)-\eta(C)\Bigr)\Bigr|\le\begin{cases}\frac\varepsilon2(n_1+1)&\text{if }n_1+1\ge n_\varepsilon,\\ n_1+1<n_\varepsilon&\text{if }n_1+1<n_\varepsilon,\end{cases}
$$
and in both cases the left side is at most $\frac\varepsilon2(n_1+1)+n_\varepsilon$. In the same way, using the second estimate with $n=n_0$,
$$
\Bigl|\sum_{i=1}^{n_0}\Bigl(\1_C\bigl(\sigma_\calX^{-i}(z,h)\bigr)-\eta(C)\Bigr)\Bigr|\le\begin{cases}\frac\varepsilon2n_0&\text{if }n_0\ge n_\varepsilon,\\ n_0<n_\varepsilon&\text{if }n_0<n_\varepsilon,\end{cases}
$$
so the left side is at most $\frac\varepsilon2n_0+n_\varepsilon$. Adding the two bounds,
$$
\Bigl|\sum_{i=-n_0}^{n_1}\Bigl(\1_C\bigl(\sigma_\calX^i(z,h)\bigr)-\eta(C)\Bigr)\Bigr|\le\frac\varepsilon2\,(n_0+n_1+1)+2n_\varepsilon .
$$
If $n_0+n_1\ge N:=\lceil4n_\varepsilon/\varepsilon\rceil$, then $2n_\varepsilon\le\frac\varepsilon2(n_0+n_1+1)$, and dividing by $n_0+n_1+1$ gives
$$
\Bigl|\frac1{n_0+n_1+1}\sum_{i=-n_0}^{n_1}\1_C\bigl(\sigma_\calX^i(z,h)\bigr)-\eta(C)\Bigr|\le\varepsilon .
$$
\end{proof}
The next lemma shows that the test sets $T_{m,k}$ depend on finitely many coordinates, and that for fixed $m$ the test at level $k$ is passed with $\zeta$-probability close to $1$ when $k$ is large. The second statement ensures that the integers $k_m$ of Definition~\ref{def:Xk} are well defined.

\begin{lemma}
\label{lem:Tmk}
Let $m\ge1$. For every $k\ge0$, the set $T_{m,k}$ is a union of elements of $\scrQ_{3(4L)^k+m}$, and in particular it is clopen. Moreover, $\zeta(T_{m,k})\to1$ as $k\to\infty$.
\end{lemma}

\begin{proof}
We start with the first statement. Let $k\ge0$ and $(z,h)\in\calX$, and write $J_k(h)=[u,v-1]$. By Lemma~\ref{lem:markers}\ref{mk:length}, $u\le0<v$ and $v-u\le3(4L)^k$, so $u$ and $v$ lie in $[-3(4L)^k,3(4L)^k]$. Since $u$ is the largest point of $\calM_k(h)$ in $(-\infty,0]$ and $v$ the smallest in $(0,\infty)$, the block $J_k(h)$ is determined by $(h_{i,k})_{i\in[-3(4L)^k,3(4L)^k]}$.

Let $C\in\scrQ_s$ with $s<m$ and $i\in J_k(h)^{(s)}$. Then $[i-s,i+s]\subset J_k(h)$, so $\1_C(\sigma_\calX^i(z,h))$ depends only on $z$ and $(h_{\cdot,k'})_{0\le k'\le s}$ on $J_k(h)$. Hence whether $(z,h)\in T_{m,k}$ depends only on $z[-3(4L)^k,3(4L)^k]$ and on $h_{i,k'}$ for $i\in[-3(4L)^k,3(4L)^k]$ and $0\le k'\le\max\{k,m-1\}$. Since $\max\{k,m-1\}\le3(4L)^k+m$, the set $T_{m,k}$ is a union of elements of $\scrQ_{3(4L)^k+m}$.

We now prove the second statement. The system $(\sigma_\calX,\zeta)$ is ergodic by Lemma~\ref{lem:base}\ref{base:mixing}. Applying Lemma~\ref{lem:birkhoff} with $\varepsilon=1/(2m)$ to each of the finitely many elements of $\scrQ_0\cup\dots\cup\scrQ_{m-1}$, we obtain a Borel set $G\subset\calX$ with $\zeta(G)=1$ with the following property. For every $(z,h)\in G$ there is $N\ge1$ such that, for every $C\in\scrQ_s$ with $s<m$ and every interval $I\ni0$ with $|I|>N$,
$$
\Bigl|\sum_{i\in I}\1_C\bigl(\sigma_\calX^i(z,h)\bigr)-|I|\,\zeta(C)\Bigr|\le\frac{|I|}{2m} .
$$
By \eqref{eq:interior}, the sum on the left differs from $S_{C,I}(z,h)$ by at most $2s<2m$, so
$$
\bigl|S_{C,I}(z,h)-|I|\,\zeta(C)\bigr|\le\frac{|I|}{2m}+2m\le\frac{|I|}m\quad\text{if }|I|\ge\max\{N+1,4m^2\}.
$$
Since $0\in J_k(h)$ and $|J_k(h)|>(4L)^k$ by Lemma~\ref{lem:markers}\ref{mk:length}, every $(z,h)\in G$ belongs to $T_{m,k}$ for all $k$ with $(4L)^k\ge\max\{N+1,4m^2\}$. Hence $\1_{T_{m,k}}\to1$ $\zeta$-almost everywhere as $k\to\infty$, and dominated convergence gives $\zeta(T_{m,k})\to1$.
\end{proof}

We now define the sets $\calX_k$ of typical points. For each $m\ge1$ we choose a level $k_m$ such that $\zeta(\calX\setminus T_{m,k_m})\le2^{-m-3}$, and we set $\calX_{k_m}:=T_{m,k_m}$. At the other levels every point is typical. This has two consequences, proved in Lemma~\ref{lem:tests}. The atypical points have small total $\zeta$-probability, and every ergodic $\sigma_\calX$-invariant probability $\eta\ne\zeta$ fails the tests at large levels.

\begin{definition}
\label{def:Xk}
Let $T_{m,k}$ be the sets of Definition~\ref{def:test}. By Lemma~\ref{lem:Tmk}, $\zeta(T_{m,k})\to1$ as $k\to\infty$ for every $m\ge1$, so the following integers $k_m$ are well defined, and we use them to define the sets $\calX_k$.
\begin{itemize}[leftmargin=*]
\item Let $k_0:=0$ and, recursively for $m\ge1$,
\begin{equation}
\label{eq:Jm}
k_m:=\min\Bigl\{k\in\Z\ ;\ k>k_{m-1},\ k\ge m\text{ and }\zeta\bigl(\calX\setminus T_{m,k}\bigr)\le2^{-m-3}\Bigr\} .
\end{equation}
\item For $k\ge0$,
$$
\calX_k:=T_{m,k}\quad\text{if }k=k_m\text{ for some }m\ge1,\quad\text{and}\quad\calX_k:=\calX\quad\text{otherwise}.
$$
\end{itemize}
\end{definition}

Let $k=k_m$, $B\in\calP_k(h)$ and $i\in B$. Then $$J_k(\sigma_\calH^ih)=B-i\ \text{and}\ S_{C,B-i}(\sigma_\calX^i(z,h))=S_{C,B}(z,h),$$ so
\begin{equation}
\label{eq:Xj}
\sigma_\calX^i(z,h)\in\calX_k\quad\text{exactly when}\quad(z,h)\in\mathcal T_m(B).
\end{equation}
For $s<m$ we have $s\le m-1\le k-1$ by \eqref{eq:Jm}, so this condition depends only on $z[B]$ and on $(h_{i',k'})_{i'\in B,\,0\le k'<k}$. For $k\notin\{k_m\ ;\ m\ge1\}$ we have $\sigma_\calX^i(z,h)\in\calX_k$ for every $i$.

Finally we define the laws of the blocks. For an interval $I$ with subintervals $I_1,\dots,I_n$, ordered from left to right, whose union is $I$, we identify $\calA_\aux^I$ with $\calA_\aux^{I_1}\times\dots\times\calA_\aux^{I_n}$. For an interval $I$ and $c\in\Z$ we define the translation
$$
\tau_{I,c}:w\in \calA_\aux^I\mapsto (w_{i-c})_{i\in I+c}\in \calA_\aux^{I+c} .
$$
\begin{lemma}
\label{lem:choice}
For all intervals $I_1,\dots,I_n$, ordered from left to right, whose union is an interval $I$, and all probabilities $\Q_j$ on sets $E_j\subset\calA_\aux^{I_j}$, functions $X_j$ on $E_j$ and numbers $b_j$ satisfying the hypotheses of Proposition~\ref{prop:products}, one can choose probabilities $\widetilde\Q$ and $\overline\Q$ with the properties of Proposition~\ref{prop:products} such that the following holds. For every $c\in\Z$, the probabilities chosen for
$$
(\tau_{I_j,c})_*\Q_j,\quad X_j\circ\tau_{I_j+c,-c}\quad\text{and}\quad b_j
$$
are $(\tau_{I,c})_*\widetilde\Q$ and $(\tau_{I,c})_*\overline\Q$.
\end{lemma}

\begin{proof}
Let $c:=\min I_1$. If $c=0$, let $\widetilde\Q$ and $\overline\Q$ be given by Proposition~\ref{prop:products} for $(\Q_j)_j$, $(X_j)_j$ and $(b_j)_j$. If $c\ne0$, the probabilities $\Q'_j:=(\tau_{I_j,-c})_*\Q_j$, the functions $X'_j:=X_j\circ\tau_{I_j-c,c}$ and the numbers $b_j$ are defined on the intervals $I_j-c$, with $\min(I_1-c)=0$, so probabilities $\widetilde\Q'$ and $\overline\Q'$ are already chosen for them, and we take
$$
\widetilde\Q:=(\tau_{I-c,c})_*\widetilde\Q'\quad\text{and}\quad\overline\Q:=(\tau_{I-c,c})_*\overline\Q' .
$$
The map $\tau_{I-c,c}$ is a bijection which maps cylinders to cylinders, and $X'_j(w)=X_j(\tau_{I_j-c,c}w)$, so it maps the set $G$ of Proposition~\ref{prop:products} for $(\Q'_j)_j$, $(X'_j)_j$, $(b_j)_j$ onto the set $G$ for $(\Q_j)_j$, $(X_j)_j$, $(b_j)_j$. Hence $\widetilde\Q$ and $\overline\Q$ have the properties of Proposition~\ref{prop:products}, and by construction the choice commutes with translations.
\end{proof}

\begin{definition}
\label{def:blocklaws}
Let $(z,h)\in\calX$. For $k\ge0$ and $B\in\calP_k(h)$ we define a probability $\P^{z,h}_B$ on $\calA_\aux^B$, recursively in $k$, as follows.
\begin{itemize}[leftmargin=*]
\item If $k=0$, then $B=\{i\}$ for some $i\in\Z$, and $\P^{z,h}_B$ is the uniform probability on $\calA_\aux^{\{i\}}$.
\item Let $k\ge1$, let $B\in\calP_k(h)$ with $B=B_1\cup\dots\cup B_n$, where $B_1,\dots,B_n\in\calP_{k-1}(h)$ are the elements contained in $B$, ordered from left to right, and assume that the probabilities $\P^{z,h}_{B_j}$ are defined. Let $E_j:=\supp\P^{z,h}_{B_j}$ and let $\Q_j$ be the restriction of $\P^{z,h}_{B_j}$ to $E_j$, so that $E:=E_1\times\dots\times E_n\subset\calA_\aux^B$. Define $X_j:E_j\to\R^{\mathcal K_{<k}}$ and $b_j\ge0$ by
$$
X_j(w):=V_{k,B_j}(z,h,w)-\P^{z,h}_{B_j}\bigl(V_{k,B_j}(z,h,\cdot)\bigr)\quad\text{and}\quad b_j:=2^{3-k}|B_j| ,
$$
so that the set $G$ of Proposition~\ref{prop:products} is
$$
G^{z,h}_B:=\Bigl\{w\in E\ ;\ \Bigl\|\sum_{j=1}^nX_j(w[B_j])\Bigr\|\le2^{-k}|B|\Bigr\}.
$$
Let $\widetilde\Q$ and $\overline\Q$ be the probabilities chosen in Lemma~\ref{lem:choice} for $(\Q_j)_j$, $(X_j)_j$ and $(b_j)_j$. Then, with the sets $\calX_k$ of Definition~\ref{def:Xk},
$$
\P^{z,h}_B:=\begin{cases}\widetilde\Q&\text{if }\sigma_\calX^{\min B}(z,h)\in\calX_k,\\ \overline\Q&\text{if }\sigma_\calX^{\min B}(z,h)\notin\calX_k.\end{cases}
$$
In both cases $\P^{z,h}_B$ is a probability on $\calA_\aux^B$ supported on $E$.
\end{itemize}
\end{definition}

Thus $\P^{z,h}_B$ keeps the joint law of any $2R$ of the restrictions $w[B_j]$, and it is supported on the $w\in\calA_\aux^B$ along which the numbers of occurrences of all elements of $\mathcal K_{<k}$, weighted as in $V_{k,B}$, are close to their means, up to the occurrences which intersect two of the intervals $B_j$. The following lemma shows that the definition makes sense and collects the properties of these laws.

\begin{lemma}
\label{lem:blocklawsint}
For every $k\ge0$, Definition~\ref{def:blocklaws} defines $\P^{z,h}_B$ for every $(z,h)\in\calX$ and every $B\in\calP_k(h)$, and the following hold.
\begin{enumerate}[label=\textup{(\roman*)},leftmargin=*]
\item\label{bl:freq} For all $w,w'\in\supp\P^{z,h}_B$,
$$
\bigl\|V_{k,B}(z,h,w)-V_{k,B}(z,h,w')\bigr\|\le3\cdot2^{-k}|B| .
$$
Hence, for every $0\le s<k$, every $\Gamma\in\mathcal K_s$ and every $w\in\supp\P^{z,h}_B$,
$$
\bigl|S_{\Gamma,B}(z,h,w)-\P^{z,h}_B\bigl(S_{\Gamma,B}(z,h,\cdot)\bigr)\bigr|\le3\cdot2^{s+1-k}|B| .
$$
\item\label{bl:marg} Let $k\ge1$ and let $B_1,\dots,B_n$ be the children of $B$. For every $T\subset[1,n]$ with $\#T\le2R$ and all $u_j\in\calA_\aux^{B_j}$, where $j\in T$,
$$
\P^{z,h}_B\bigl(\bigl\{w\in\calA_\aux^B\ ;\ w[B_j]=u_j\text{ for all }j\in T\bigr\}\bigr)=\prod_{j\in T}\P^{z,h}_{B_j}\bigl(\{u_j\}\bigr).
$$
\item\label{bl:entropy} Let $k\ge1$, let $B_1,\dots,B_n$ be the children of $B$ (see Lemma~\ref{lem:markers}\ref{mk:children}), let $u:=\min B$, and for $F\subset B$ let $\scrP_F$ be the partition of $\calA_\aux^B$ into the cylinders over $F$. If $\sigma_\calX^u(z,h)\in\calX_k$, then
$$
\sum_{j=1}^nH_{\P^{z,h}_B}\bigl(\scrP_{B_j}\bigr)-H_{\P^{z,h}_B}\bigl(\scrP_B\bigr)\le\log\frac43 .
$$
If $\sigma_\calX^u(z,h)\notin\calX_k$, then $\#\supp\P^{z,h}_B\le2^{9Rq|B|/(2L)+12L}$.
\item\label{bl:cov} Let $(z',h')\in\calX$ and $c\in\Z$ with $B+c\in\calP_k(h')$, and assume that $z'_{i+c}=z_i$ and $h'_{i+c,k'}=h_{i,k'}$ for all $i\in B$ and $0\le k'<k$. Then $\P^{z',h'}_{B+c}=(\tau_{B,c})_*\P^{z,h}_B$.
\end{enumerate}
\end{lemma}

\begin{proof}
We argue by induction on $k$. Let $k=0$ and $B=\{i\}$. Then $\R^{\mathcal K_{<0}}=\{0\}$, so \ref{bl:freq} holds. Property \ref{bl:cov} holds since $\tau_{B,c}$ maps the uniform probability on $\calA_\aux^B$ to the uniform probability on $\calA_\aux^{B+c}$, and \ref{bl:marg} and \ref{bl:entropy} concern $k\ge1$.

Let $k\ge1$ and assume that the lemma holds at level $k-1$. Fix $(z,h)\in\calX$ and $B\in\calP_k(h)$ with children $B_1,\dots,B_n$, let $u:=\min B$, and write $\P_B:=\P^{z,h}_B$ and $\P_{B_j}:=\P^{z,h}_{B_j}$. We divide the proof into four steps.

\begin{step}[1]
We show that $(\Q_j)_j$, $(X_j)_j$ and $(b_j)_j$ satisfy the hypotheses of Proposition~\ref{prop:products}, so that Definition~\ref{def:blocklaws} makes sense.
\end{step}

\begin{stepproof}
By Lemma~\ref{lem:markers}\ref{mk:children} and \eqref{eq:a}, $n>4L/3\ge\max\{192R,3072\}$, and since $n\le12L$,
$$
(2R+1)(2n)^{2R}e^{-n/1536}\le(2R+1)(24L)^{2R}e^{-L/1152}\le\frac14 ,
$$
so \eqref{eq:nk} holds. We have $$\sum_jb_j=2^{3-k}|B|,\quad\text{and}\quad\max_jb_j\le\frac3n\sum_jb_j$$ since, by Lemma~\ref{lem:markers}\ref{mk:children}, $|B_j|\le {3|B|}/n$. The function $X_j$ has mean zero under $\Q_j$.

Let $w,w'\in E_j$. The coordinates of $V_{k,B_j}(z,h,\cdot)$ indexed by $\mathcal K_{<k-1}$ are those of $V_{k-1,B_j}(z,h,\cdot)$, so by \ref{bl:freq} at level $k-1$ their differences at $w$ and $w'$ form a vector of norm at most $3\cdot2^{1-k}|B_j|$. The coordinates indexed by $\mathcal K_{k-1}$ form a vector of norm at most $2^{-k}|B_j|$ at $w$ and at $w'$, by \eqref{eq:Vnorm} with $s=k-1$. Hence
$$
\bigl\|V_{k,B_j}(z,h,w)-V_{k,B_j}(z,h,w')\bigr\|\le3\cdot2^{1-k}|B_j|+2\cdot2^{-k}|B_j|=b_j .
$$
The vector $\P_{B_j}\bigl(V_{k,B_j}(z,h,\cdot)\bigr)$ is a convex combination of the vectors $V_{k,B_j}(z,h,w')$, $w'\in E_j$, so $\norm{X_j(w)}\le b_j$.
\end{stepproof}

\begin{step}[2]
We prove \ref{bl:marg} and \ref{bl:entropy}.
\end{step}

\begin{stepproof}
By Proposition~\ref{prop:products}\ref{cp:supp}, in both cases of Definition~\ref{def:blocklaws}, $\P_B(G^{z,h}_B)=1$, where $$G^{z,h}_B\subset E=E_1\times\dots\times E_n.$$ Let $T\subset[1,n]$ with $\#T\le2R$, and let $u_j\in\calA_\aux^{B_j}$ for $j\in T$. If $u_j\in E_j$ for every $j\in T$, choose $x\in E$ with $x_j=u_j$ for $j\in T$. Since $\Q_j$ is the restriction of $\P_{B_j}$ to $E_j$, Proposition~\ref{prop:products}\ref{cp:marg} with $S=T$ gives
$$
\P_B\bigl(\bigl\{w\in\calA_\aux^B\ ;\ w[B_j]=u_j\text{ for all }j\in T\bigr\}\bigr)=\prod_{j\in T}\Q_j\bigl(\{u_j\}\bigr)=\prod_{j\in T}\P_{B_j}\bigl(\{u_j\}\bigr).
$$
If $u_j\notin E_j=\supp\P_{B_j}$ for some $j\in T$, then both sides vanish, since $\P_B(E)=1$. This proves \ref{bl:marg}. If $\sigma_\calX^u(z,h)\in\calX_k$, then $\P_B=\widetilde\Q$ by Definition~\ref{def:blocklaws}. Under the identification of $\calA_\aux^B$ with $\calA_\aux^{B_1}\times\dots\times\calA_\aux^{B_n}$, the intersections with $E$ of the elements of the partitions $\scrP_{B_j}$ and $\scrP_B$ of $\calA_\aux^B$ are the elements of the partitions $\scrP_{\{j\}}$ and $\scrP_{[1,n]}$ of $E$ in Proposition~\ref{prop:products}\ref{cp:entropy}. Since $\P_B(E)=1$, \ref{bl:entropy} is Proposition~\ref{prop:products}\ref{cp:entropy}. If $\sigma_\calX^u(z,h)\notin\calX_k$, then by Proposition~\ref{prop:products}\ref{cp:support}, $\#E_j\le\#\calA_\aux^{B_j}=2^{q|B_j|}$ and Lemma~\ref{lem:markers}\ref{mk:children},
$$
\#\supp\P_B\le2^n\Bigl(\max_j2^{q|B_j|}\Bigr)^{2R}\le2^{n+6qR|B|/n}\le2^{12L+9qR|B|/(2L)},
$$
since $n\le12L$ and $3/n< 9/(4L)$.
\end{stepproof}

\begin{step}[3]
We prove \ref{bl:freq}.
\end{step}

\begin{stepproof}
Let $0\le s<k$, and let $i\in B^{(s)}$ with $i\notin B_j^{(s)}$ for every $1\le j\le n$. Then $[i-s,i+s]\subset B$, but $[i-s,i+s]$ is not contained in any $B_j$. Hence $[i-s,i+s]$ contains both $\min B_j-1$ and $\min B_j$ for some $2\le j\le n$, that is,
$$
i\in\bigl[\min B_j-s,\ \min B_j+s-1\bigr]\quad\text{for some }2\le j\le n .
$$
Each of these $n-1$ intervals has $2s$ points, so there are at most $2s(n-1)$ such points $i$. Moreover, for every $i$ the point $\widehat\sigma^i(z,h,w)$ lies in exactly one element of $\mathcal K_s$. Hence, for $w\in\calA_\aux^B$,
$$
V_{k,B}(z,h,w)=\sum_{j=1}^nV_{k,B_j}(z,h,w[B_j])+V^\partial(w),
$$
where the coordinates of $V^\partial(w)\in\R^{\mathcal K_{<k}}$ are nonnegative and, for every $0\le s<k$, those indexed by $\mathcal K_s$ have sum at most $2^{-s-1}\,2s(n-1)=s2^{-s}(n-1)$. Since $\sum_{s\ge1}s2^{-s}=2$ and $n-1<|B|(4L)^{1-k}$ by Lemma~\ref{lem:markers}\ref{mk:children},
$$
\norm{V^\partial(w)}\le\sum_{s=1}^{k-1}s2^{-s}(n-1)\le2(4L)^{1-k}|B| ,
$$
and $V^\partial(w)=0$ if $k=1$. Define $m_B:=\sum_j\P_{B_j}\bigl(V_{k,B_j}(z,h,\cdot)\bigr)$, which does not depend on $w$. Then
$$
V_{k,B}(z,h,w)=m_B+\sum_{j=1}^nX_j(w[B_j])+V^\partial(w)\quad\text{for }w\in E .
$$
Let $w,w'\in\supp\P_B$. Since $\P_B(G^{z,h}_B)=1$, we have $w,w'\in G^{z,h}_B$, and therefore
$$
\bigl\|V_{k,B}(z,h,w)-V_{k,B}(z,h,w')\bigr\|\le2\cdot2^{-k}|B|+\norm{V^\partial(w)}+\norm{V^\partial(w')} .
$$
If $k=1$, the right side is $2\cdot2^{-k}|B|$. If $k\ge2$, then, since $4L\ge16$ by \eqref{eq:a},
$$
4(4L)^{1-k}=2^{-k}\cdot8\Bigl(\frac1{2L}\Bigr)^{k-1}\le2^{-k}\cdot\frac4L\le2^{-k},
$$
so the right side is at most $3\cdot2^{-k}|B|$. This is the first statement of \ref{bl:freq}. The mean of $V_{k,B}(z,h,\cdot)$ under $\P_B$ is a convex combination of the vectors $V_{k,B}(z,h,w')$ with $w'\in\supp\P_B$. Hence
$$
\bigl\|V_{k,B}(z,h,w)-\P_B\bigl(V_{k,B}(z,h,\cdot)\bigr)\bigr\|\le3\cdot2^{-k}|B| .
$$
For $\Gamma\in\mathcal K_s$ with $s<k$, the coordinate $\Gamma$ of this vector is $2^{-s-1}\bigl(S_{\Gamma,B}(z,h,w)-\P_B(S_{\Gamma,B}(z,h,\cdot))\bigr)$, which gives the second statement.
\end{stepproof}

\begin{step}[4]
Finally, we prove \ref{bl:cov}.
\end{step}

\begin{stepproof}
Let $(z',h')$ and $c$ be as in \ref{bl:cov}, and write $B=[u,v-1]$. Since $B\in\calP_k(h)$ and $B+c\in\calP_k(h')$, we have $v\in\calM_k(h)\subset\calM_{k-1}(h)$ and $v+c\in\calM_k(h')\subset\calM_{k-1}(h')$. Together with $h'_{i+c,k-1}=h_{i,k-1}$ for $i\in B$, this shows that the points of $\calM_{k-1}(h')\cap[u+c,v+c]$ are those of $\calM_{k-1}(h)\cap[u,v]$ translated by $c$. By Lemma~\ref{lem:markers}\ref{mk:children}, the elements of $\calP_{k-1}(h')$ contained in $B+c$ are the intervals $B_j+c$. The hypotheses of \ref{bl:cov} at level $k-1$ hold for each $B_j$, so $\P^{z',h'}_{B_j+c}=(\tau_{B_j,c})_*\P_{B_j}$.

Let $0\le s<k$, $\Gamma\in\mathcal K_s$, $w\in\calA_\aux^{B_j}$ and $i\in B_j^{(s)}$. Since $[i-s,i+s]\subset B_j$ and $s\le k-1$, the points $\widehat\sigma^i(z,h,w)$ and $\widehat\sigma^{i+c}(z',h',\tau_{B_j,c}w)$ have the same coordinates $z$, $(h_{\cdot,k'})_{k'\le s}$ and $\omega$ on $[-s,s]$. Hence
$$
S_{\Gamma,B_j+c}(z',h',\tau_{B_j,c}w)=S_{\Gamma,B_j}(z,h,w)\quad\text{and}\quad V_{k,B_j+c}(z',h',\tau_{B_j,c}w)=V_{k,B_j}(z,h,w).
$$
If $k=k_m$, then $s\le k-1$ for every $s<m$ by \eqref{eq:Jm}, so in the same way
$$
S_{C,B+c}(z',h')=S_{C,B}(z,h)\quad\text{for every }C\in\scrQ_s\text{ with }s<m,
$$
and by \eqref{eq:Xj}, $\sigma_\calX^{u+c}(z',h')\in\calX_k$ exactly when $\sigma_\calX^u(z,h)\in\calX_k$. If $k\notin\{k_m\ ;\ m\ge1\}$, both points belong to $\calX_k=\calX$.

Let $\Q_j$, $X_j$ and $b_j$ be the probabilities, functions and numbers of Definition~\ref{def:blocklaws} for $(z,h)$ and $B$. Since $\P^{z',h'}_{B_j+c}=(\tau_{B_j,c})_*\P_{B_j}$ and by the two displays above, those for $(z',h')$ and $B+c$ are
$$
(\tau_{B_j,c})_*\Q_j,\quad X_j\circ\tau_{B_j+c,-c}\quad\text{and}\quad b_j ,
$$
and either $\P_B=\widetilde\Q$ and $\P^{z',h'}_{B+c}=\widetilde\Q'$, or $\P_B=\overline\Q$ and $\P^{z',h'}_{B+c}=\overline\Q'$, where $\widetilde\Q'$ and $\overline\Q'$ are the probabilities chosen for $(z',h')$ and $B+c$. By Lemma~\ref{lem:choice}, $\widetilde\Q'=(\tau_{B,c})_*\widetilde\Q$ and $\overline\Q'=(\tau_{B,c})_*\overline\Q$, so $\P^{z',h'}_{B+c}=(\tau_{B,c})_*\P_B$.
\end{stepproof}
\end{proof}

The following lemma lists the properties of the family of sets $(\calX_k)_k$.

\begin{lemma}
\label{lem:tests}
The following hold.
\begin{enumerate}[label=\textup{(\roman*)},leftmargin=*]
\item\label{test:zeta} $\sum_{k\ge1}\zeta(\calX\setminus\calX_k)\le\frac18$.
\item\label{test:other} For every ergodic $\sigma_\calX$-invariant Borel probability $\eta$ on $\calX$ with $\eta\ne\zeta$, $$\eta(\calX_{k_m})\to0\text{ as }m\to\infty.$$
\end{enumerate}
\end{lemma}

\begin{proof}
We start by proving \ref{test:zeta}. If $k\notin\{k_m\ ;\ m\ge1\}$, then $\calX_k=\calX$, and $\calX_{k_m}=T_{m,k_m}$. By \eqref{eq:Jm},
$$
\sum_{k\ge1}\zeta(\calX\setminus\calX_k)=\sum_{m\ge1}\zeta\bigl(\calX\setminus T_{m,k_m}\bigr)\le\sum_{m\ge1}2^{-m-3}=\frac18 .
$$

We now prove \ref{test:other}. Since $\eta\ne\zeta$ and $\{\emptyset\}\cup\bigcup_{s\ge0}\scrQ_s$ is closed under intersections and generates the Borel $\sigma$-algebra of $\calX$, there are $s\ge0$ and $C\in\scrQ_s$ with $\eta(C)\ne\zeta(C)$. Define $\varepsilon:=|\eta(C)-\zeta(C)|>0$. By Lemma~\ref{lem:birkhoff}, for $\eta$-almost every $(z,h)$ there is $N_0$ such that
$$
\Bigl|\frac1{|I|}\sum_{i\in I}\1_C\bigl(\sigma_\calX^i(z,h)\bigr)-\eta(C)\Bigr|\le\frac\varepsilon8\quad\text{for every interval }I\ni0\text{ with }|I|\ge N_0 .
$$
Fix such $(z,h)$, and let $m>s$ be so large that $1/m<\varepsilon/2$ and $(4L)^{k_m}\ge\max\{N_0,16s/\varepsilon\}$. Let $I:=J_{k_m}(h)$, so that $|I|>(4L)^{k_m}$. By \eqref{eq:interior},
$$
\Bigl|\frac{S_{C,I}(z,h)}{|I|}-\eta(C)\Bigr|\le\frac\varepsilon8+\frac{2s}{|I|}\le\frac\varepsilon4,\quad\text{hence}\quad\bigl|S_{C,I}(z,h)-|I|\,\zeta(C)\bigr|\ge\frac34\varepsilon|I|>\frac{|I|}m,
$$
so $(z,h)\notin\mathcal T_m(I)$ by \eqref{eq:test}, that is, $(z,h)\notin T_{m,k_m}=\calX_{k_m}$. Hence $\1_{\calX_{k_m}}\to0$ $\eta$-almost everywhere as $m\to\infty$, and $\eta(\calX_{k_m})\to0$ by dominated convergence.
\end{proof}

\begin{proof}[Proof of Proposition~\ref{prop:blocklaws}]
Property \ref{P:unif} is the case $k=0$ of Definition~\ref{def:blocklaws}. Properties \ref{P:freq}, \ref{P:marg} and \ref{P:entropy} are Lemma~\ref{lem:blocklawsint}\ref{bl:freq}, \ref{bl:marg} and \ref{bl:entropy}, and \ref{P:cov} is Lemma~\ref{lem:blocklawsint}\ref{bl:cov}, since $\tau_{B,c}(w)=(w_{i-c})_{i\in B+c}$. The first statement of \ref{P:typ} is \eqref{eq:Xj} and the remark which follows it. The second statement is Lemma~\ref{lem:tests}\ref{test:zeta}, and the third follows from Lemma~\ref{lem:tests}\ref{test:other}, since $k_m\to\infty$.
\end{proof}

\section{Exponential mixing of order \texorpdfstring{$R$}{R}}
\label{sec:mixing}

In this section we prove that $(\sigma,\mu_0)$ is exponentially mixing of order $R$ for every H\"older class. The mixing follows from exact independence of the coordinates on $R$ distant intervals. Recall that $\mu_0$ is the pushforward under the map $\widehat\pi(z,h,\omega)=(z,\omega)$ of the measure $$\widehat\mu_0(dz,dh,d\omega)=\mu_{z,h}(d\omega)\,\zeta(dz,dh)$$ (see \eqref{eq:mu0}), where $\zeta=\nu\times\P_\calH$ and $\mu_{z,h}$ are the sample measures of Lemma~\ref{lem:sample}.

We use the notation of \ref{not:intervals} for intervals. For an interval $I=[a,b]$ and $c\ge0$,
$$
I^{+c}=\bigl\{i\in\Z\ ;\ |i-i'|\le c\text{ for some }i'\in I\bigr\}=[a-c,b+c]
$$
is the interval $I$ enlarged by $c$ on each side, and $\dist(I,J)=\min\{|i-j|\ ;\ i\in I,\ j\in J\}$ for intervals $I,J$.

The next proposition is the key to the mixing of $(\sigma,\mu_0)$. It says that, under $\mu_0$, blocks of coordinates which are far apart compared with their lengths are exactly independent.
\begin{proposition}
\label{prop:independence}
Define $K:=(K_R+1)(24L+1)+4DL$, where $K_R$ is the constant of Theorem~\ref{thm:BCC}\ref{BCC:ii} and $D$ is the constant of Proposition~\ref{prop:markers}. Let $\ell\ge1$ and let $I_1,\dots,I_R$ be intervals, ordered from left to right, of length at most $\ell$, with $\dist(I_s,I_{s+1})>K\ell$ for $1\le s<R$. Then under $\mu_0$ the blocks $x[I_1],\dots,x[I_R]$ are jointly independent.
\end{proposition}

\begin{proof}
Let $k\ge0$ be the least integer with $(4L)^k\ge\ell$, so that $(4L)^k\le4L\ell$, and let $r:=3(4L)^k$. Then $|I_s^{+r}|\le\ell+6(4L)^k\le(24L+1)\ell$ and, for $1\le s<R$,
\begin{equation}
\label{eq:Wdist}
\dist(I_s^{+r},I_{s+1}^{+r})\ge\dist(I_s,I_{s+1})-6(4L)^k>(K-24L)\ell\ge K_R(24L+1)\ell+4DL\ell .
\end{equation}
In particular the sets $I_s^{+r}$ are pairwise disjoint. We divide the proof into three steps.

\begin{step}[1]
We show that, for words $w_s\in\calA_\aux^{I_s}$ and $(z,h)\in\calX$,
$$
\mu_{z,h}\bigl(\omega[I_s]=w_s\text{ for }1\le s\le R\bigr)=\prod_{s=1}^R\mu_{z,h}\bigl(\omega[I_s]=w_s\bigr),
$$
and that $$\mu_{z,h}(\omega[I_s]=w_s)=\mu_{z',h'}(\omega[I_s]=w_s)$$ whenever $z[I_s^{+r}]=z'[I_s^{+r}]$ and $h_{i,k'}=h'_{i,k'}$ for all $i\in I_s^{+r}$ and $0\le k'\le k$.
\end{step}

\begin{stepproof}
Let $\calF_s$ be the set of elements of $\calP_k(h)$ which intersect $I_s$.

To apply Lemma~\ref{lem:sample}, we first show that the sets $\calF_s$ are pairwise disjoint, that their union has at most $2R$ elements, and that every element of $\calF_s$ is contained in $I_s^{+r}$. If $k\ge1$, every element of $\calP_k(h)$ has length greater than $(4L)^k\ge\ell\ge|I_s|$ by Lemma~\ref{lem:markers}\ref{mk:length}, so $\#\calF_s\le2$. If $k=0$, then $\ell=1$ and $\#\calF_s=1$. Every element of $\calF_s$ intersects $I_s$ and has length at most $3(4L)^k=r$, so it is contained in $I_s^{+r}$. Since the sets $I_s^{+r}$ are pairwise disjoint, so are $\calF_1,\dots,\calF_R$, and $\calF:=\calF_1\cup\dots\cup\calF_R$ has at most $2R$ elements.

We now prove the product formula. The event $\{\omega[I_s]=w_s\}$ depends only on the restrictions $\omega[B]$, $B\in\calF_s$. By Lemma~\ref{lem:sample} applied to $\calF$, the restrictions $\omega[B]$, $B\in\calF$, are independent under $\mu_{z,h}$ with laws $\P^{z,h}_B$. Hence the events $\{\omega[I_s]=w_s\}$, $1\le s\le R$, are independent under $\mu_{z,h}$, which is the product formula, and
\begin{equation}
\label{eq:Isprod}
\mu_{z,h}\bigl(\omega[I_s]=w_s\bigr)=\Bigl(\bigotimes_{B\in\calF_s}\P^{z,h}_B\Bigr)\bigl(\omega[I_s]=w_s\bigr).
\end{equation}

Finally we prove the last statement. Let $(z',h')\in\calX$ with $z'[I_s^{+r}]=z[I_s^{+r}]$ and $h'_{i,k'}=h_{i,k'}$ for all $i\in I_s^{+r}$ and $0\le k'\le k$. From \eqref{eq:Isprod}, it suffices to show that the elements of $\calP_k(h')$ which intersect $I_s$ are the elements of $\calF_s$, and that $$\P^{z',h'}_B=\P^{z,h}_B\ \text{for }B\in\calF_s.$$ The elements of $\calF_s$ are the intervals $[u,v-1]$ with $u<v$ consecutive points of $\calM_k(h)$ and $[u,v-1]\cap I_s\ne\emptyset$. For such $u,v$ we have $v-u\le r$ by Lemma~\ref{lem:markers}\ref{mk:length}, so $u,v\in I_s^{+r}$. Since $\calM_k(h')\cap I_s^{+r}=\calM_k(h)\cap I_s^{+r}$, the same holds with $h'$ in place of $h$, and the two families of intervals coincide. For $B\in\calF_s$ we have $B\subset I_s^{+r}$, so $z'[B]=z[B]$ and $h'_{i,k'}=h_{i,k'}$ for $i\in B$ and $0\le k'<k$, and Proposition~\ref{prop:blocklaws}\ref{P:cov} with $c=0$ gives $\P^{z',h'}_B=\P^{z,h}_B$.
\end{stepproof}

\begin{step}[2]
We show that the families $\bigl(z[I_s^{+r}],(h_{i,k'})_{i\in I_s^{+r},\,0\le k'\le k}\bigr)$, $1\le s\le R$, are jointly independent under $\zeta$.
\end{step}

\begin{stepproof}
By \eqref{eq:Wdist}, $\dist(I_s^{+r},I_{s+1}^{+r})>K_R(24L+1)\ell\ge K_R\max_{s'}|I_{s'}^{+r}|$, so Theorem~\ref{thm:BCC}\ref{BCC:ii} shows that $z[I_1^{+r}],\dots,z[I_R^{+r}]$ are jointly independent under $\nu$. By \eqref{eq:Wdist}, $\dist(I_s^{+r},I_{s'}^{+r})>4DL\ell\ge D(4L)^k$ for $s\ne s'$, so Lemma~\ref{lem:markers}\ref{mk:dep} and \ref{not:dep} show that the restrictions of $\bigl((h_{i,0},\dots,h_{i,k})\bigr)_{i\in\Z}$ to $I_1^{+r},\dots,I_R^{+r}$ are jointly independent under $\P_\calH$. Since $\zeta=\nu\times\P_\calH$, the $R$ families are jointly independent.
\end{stepproof}

\begin{step}[3]
We prove the proposition.
\end{step}

\begin{stepproof}
Let $u_s\in\{-1,1\}^{I_s}$ and $w_s\in\calA_\aux^{I_s}$ for $1\le s\le R$. By \eqref{eq:mu0} and Step~1,
$$
\mu_0\bigl(z[I_s]=u_s,\ \omega[I_s]=w_s\text{ for }1\le s\le R\bigr)=\int_\calX\prod_{s=1}^R\1_{\{z[I_s]=u_s\}}\,\mu_{z,h}\bigl(\omega[I_s]=w_s\bigr)\,d\zeta(z,h).
$$
Since $I_s\subset I_s^{+r}$, the $s$-th factor depends only on the $s$-th family of Step~2. By Step~2 the integral is the product of the integrals of the factors, and, by \eqref{eq:mu0} again, the $s$-th integral is $\mu_0(z[I_s]=u_s,\ \omega[I_s]=w_s)$.
\end{stepproof}
\end{proof}

We now show that $(\sigma,\mu_0)$ is exponentially mixing of order $R$ for every H\"older class.
\begin{theorem}
\label{thm:mixing}
Let $K$ be as in Proposition~\ref{prop:independence} and $0<\theta<1$, and define $c_\theta:=\frac{|\log\theta|}{4(K+1)}$ and $C_\theta:=2R\theta^{-1}$. Then
$$
\Bigl|\int\prod_{s=1}^RF_s\circ\sigma^{t_s}\,d\mu_0-\prod_{s=1}^R\int F_s\,d\mu_0\Bigr|\le C_\theta e^{-c_\theta\Delta(t_1,\dots,t_R)}\prod_{s=1}^R\norm{F_s}_\theta
$$
for all $F_1,\dots,F_R\in\calF_\theta(\calA^\Z)$ and all integers $0\le t_1<\dots<t_R$. In particular $(\sigma,\mu_0)$ is exponentially mixing of order $R$ for every H\"older class (see Section~\ref{sec:mixing-def}).
\end{theorem}

\begin{proof}
Write $\Delta:=\Delta(t_1,\dots,t_R)$. If $\Delta<4(K+1)$, the left side is at most $$2\prod_s\norm{F_s}_\infty\le2\theta^{-1}e^{-c_\theta\Delta}\prod_s\norm{F_s}_\theta,$$ since $e^{c_\theta\Delta}<e^{|\log\theta|}=\theta^{-1}$. Assume $\Delta\ge4(K+1)$, and define $\ell:=\lfloor\Delta/(4(K+1))\rfloor\ge1$ and $J_s:=[t_s-\ell,t_s+\ell]$. The intervals $J_s$ have length $2\ell+1$, and since $$(K+1)(2\ell+1)\le\frac\Delta2+K+1\le\frac34\Delta<\Delta,$$
we obtain
$$
\dist(J_s,J_{s+1})\ge\Delta-2\ell>(K+1)(2\ell+1)-2\ell>K(2\ell+1).
$$
Fix $a_*\in\calA$. For $x\in\calA^\Z$ let $x^{(\ell)}$ agree with $x$ on $[-\ell,\ell]$ and equal $a_*$ elsewhere, and define $F_{s,\ell}(x):=F_s(x^{(\ell)})$. By \eqref{eq:cylinderapprox}, $\norm{F_{s,\ell}}_\infty\le\norm{F_s}_\infty$ and $\norm{F_s-F_{s,\ell}}_\infty\le\Lip_\theta(F_s)\,\theta^{\ell+1}$. The function $F_{s,\ell}\circ\sigma^{t_s}$ depends only on $x[J_s]$, so Proposition~\ref{prop:independence}, applied with $2\ell+1$ in place of $\ell$, and the invariance of $\mu_0$ give
$$
\int\prod_{s=1}^RF_{s,\ell}\circ\sigma^{t_s}\,d\mu_0=\prod_{s=1}^R\int F_{s,\ell}\,d\mu_0 .
$$
For any $a_1,\dots,a_R,b_1,\dots,b_R$, the telescoping identity is
$$
\prod_{s=1}^Ra_s-\prod_{s=1}^Rb_s=\sum_{s=1}^R\Bigl(\prod_{i<s}a_i\Bigr)(a_s-b_s)\Bigl(\prod_{i>s}b_i\Bigr).
$$
Applying it pointwise with $a_s=F_s\circ\sigma^{t_s}$ and $b_s=F_{s,\ell}\circ\sigma^{t_s}$, and then with $a_s=\mu_0(F_s)$ and $b_s=\mu_0(F_{s,\ell})$, we obtain
$$
\Bigl|\int\prod_{s=1}^RF_s\circ\sigma^{t_s}\,d\mu_0-\prod_{s=1}^R\int F_s\,d\mu_0\Bigr|\le2R\,\theta^{\ell+1}\prod_{s=1}^R\norm{F_s}_\theta\le2R\,e^{-c_\theta\Delta}\prod_{s=1}^R\norm{F_s}_\theta ,
$$
since $\ell+1>\Delta/(4(K+1))$.
\end{proof}

\section{Entropy gap}
\label{sec:entropy}

In this section we prove that $\mu_0$ has entropy larger by at least $\frac q2\log2$ than every other ergodic measure of $(Y,\sigma)$, and we deduce the symbolic theorem. The lower bound for $h_{\mu_0}(\sigma)$ uses the first case of Proposition~\ref{prop:blocklaws}\ref{P:entropy}, and the upper bound for the other ergodic measures uses its second case. Both use Proposition~\ref{prop:blocklaws}\ref{P:typ}.

For $k\ge1$, $(z,h)\in\calX$ and $B\in\calP_k(h)$ with children $B_1,\dots,B_n$ (see Lemma~\ref{lem:markers}\ref{mk:children}), we define
$$
\delta^{z,h}_B:=\sum_{j=1}^nH_{\P^{z,h}_B}\bigl(\scrP_{B_j}\bigr)-H_{\P^{z,h}_B}\bigl(\scrP_B\bigr),
$$
where, as in Proposition~\ref{prop:blocklaws}\ref{P:entropy}, $\scrP_F$ is the partition of $\calA_\aux^B$ into the cylinders over $F\subset B$. Since $\scrP_B=\scrP_{B_1}\vee\dots\vee\scrP_{B_n}$, we have $\delta^{z,h}_B\ge0$ by \eqref{eq:shannon}. By Proposition~\ref{prop:blocklaws}\ref{P:marg} with $T=\{j\}$, we have $\P^{z,h}_B\bigl(\{w\ ;\ w[B_j]=u\}\bigr)=\P^{z,h}_{B_j}(\{u\})$ for every $u\in\calA_\aux^{B_j}$, so
$$
H_{\P^{z,h}_B}\bigl(\scrP_{B_j}\bigr)=H_{\P^{z,h}_{B_j}}\bigl(\scrP_{B_j}\bigr),
$$
where on the right $\scrP_{B_j}$ is the partition of $\calA_\aux^{B_j}$ into points. For $i\in\Z$ the probability $\P^{z,h}_{\{i\}}$ is uniform on $\calA_\aux^{\{i\}}$ by Proposition~\ref{prop:blocklaws}\ref{P:unif}, so $H_{\P^{z,h}_{\{i\}}}(\scrP_{\{i\}})=q\log2$. Induction on $k$ gives, for $B\in\calP_k(h)$,
\begin{equation}
\label{eq:telescope}
H_{\P^{z,h}_B}\bigl(\scrP_B\bigr)=q\log2\,|B|-\sum_{k'=1}^{k}\ \sum_{\substack{B'\in\calP_{k'}(h)\\ B'\subset B}}\delta^{z,h}_{B'} .
\end{equation}

The next proposition shows that $\mu_0$ has large entropy.  Its auxiliary coordinates contribute at least $7\log2/8$ to its entropy, up to a small error, out of the maximal possible contribution $q\log2$.

\begin{proposition} We have that
\label{prop:lower}
$$
h_{\mu_0}(\sigma)\ge\frac78\,q\log2-\frac{\log\frac43}{4L-1}.
$$
\end{proposition}

\begin{proof}
We use the relative entropy $h_{\widehat\rho}(\widehat\sigma\mid\calX)$ of Section~\ref{sec:entropy-def}, with $X'=\calX$ and $X=\calX\times\calA_\aux^\Z$, and we write $\calB_\calX$ for the Borel $\sigma$-algebra of $\calX$. For $k\ge0$ the formula $\widehat\mu^{(k)}(dz,dh,d\omega):=\mu^{(k)}_{z,h}(d\omega)\,\zeta(dz,dh)$, where $\zeta=\nu\times\P_\calH$ is the probability on $\calX$ of Section~\ref{sec:base}, defines a probability, by Lemma~\ref{lem:sample}\ref{sm:cont}, as $\widehat\mu_0$ in Definition~\ref{def:Y}, and by Lemma~\ref{lem:sample}\ref{sm:cov} and the invariance of $\zeta$ it is $\widehat\sigma$-invariant, as in the proof of Lemma~\ref{lem:Yhat}\ref{Y:mu}. We divide the proof into three steps.

\begin{step}[1]
We show that, for $k\ge1$,
$$
h_{\widehat\mu^{(k)}}(\widehat\sigma\mid\calX)\ge\frac78\,q\log2-\frac{\log\frac43}{4L-1}.
$$
\end{step}

\begin{stepproof}
For $1\le k'\le k$ and $(z,h)\in\calX$ define $\phi_{k'}(z,h):=\delta^{z,h}_{J}/|J|$, where $J:=J_{k'}(h)$ is the element of $\calP_{k'}(h)$ which contains $0$ (see Section~\ref{sec:blocklaws}). By Proposition~\ref{prop:blocklaws}\ref{P:cov} with $c=-i$, applied to the block and to its children, $$\frac{\delta^{z,h}_{B'}}{|B'|}=\phi_{k'}(\sigma_\calX^i(z,h))$$ for every $B'\in\calP_{k'}(h)$ and $i\in B'$, and by Proposition~\ref{prop:blocklaws}\ref{P:cov} with $c=0$, $\phi_{k'}$ is locally constant.

Let $n\ge1$ and $(z,h)\in\calX$, and let $\calF$ be the set of elements of $\calP_k(h)$ contained in $[0,n-1]$. Since $\bigcup\calF\subset[0,n-1]$, the partition $\scrP^n$ of $\calA_\aux^\Z$ refines $\bigvee_{B\in\calF}\scrP_B$, and under $\mu^{(k)}_{z,h}$ the partitions $\scrP_B$, $B\in\calF$, are independent and $\omega[B]$ has law $\P^{z,h}_B$. Hence, by \eqref{eq:shannon},
$$
H_{\mu^{(k)}_{z,h}}(\scrP^n)\ge H_{\mu^{(k)}_{z,h}}\Bigl(\bigvee_{B\in\calF}\scrP_B\Bigr)=\sum_{B\in\calF}H_{\P^{z,h}_B}\bigl(\scrP_B\bigr).
$$
As in Step~2 of the proof of Theorem~\ref{thm:faithful}, at most $6(4L)^k$ points of $[0,n-1]$ lie outside $\bigcup\calF$. By \eqref{eq:telescope}, since the blocks $B'\subset B\in\calF$ are disjoint for each $k'$ and contained in $[0,n-1]$, and since $\phi_{k'}\ge0$,
$$
\sum_{B\in\calF}H_{\P^{z,h}_B}\bigl(\scrP_B\bigr)\ge(n-6(4L)^k)\,q\log2-\sum_{k'=1}^{k}\ \sum_{i=0}^{n-1}\phi_{k'}\bigl(\sigma_\calX^i(z,h)\bigr).
$$
Integrating with respect to $\zeta$, using \eqref{eq:relsample} and the invariance of $\zeta$, and dividing by $n$,
$$
\frac1nH_{\widehat\mu^{(k)}}\bigl(\scrP^n_\aux\mid\pi_\calX^{-1}\calB_\calX\bigr)\ge\Bigl(1-\frac{6(4L)^k}n\Bigr)q\log2-\sum_{k'=1}^{k}\zeta(\phi_{k'}).
$$
Letting $n\to\infty$ gives $h_{\widehat\mu^{(k)}}(\widehat\sigma\mid\calX)\ge q\log2-\sum_{k'=1}^{k}\zeta(\phi_{k'})$.

We now bound $\zeta(\phi_{k'})$. Let $J:=J_{k'}(h)$ with children $J_1,\dots,J_n$. If $(z,h)\in\calX_{k'}$, then $\sigma_\calX^{\min J}(z,h)\in\calX_{k'}$ by Proposition~\ref{prop:blocklaws}\ref{P:typ}, so $\delta^{z,h}_J\le\log\frac43$ by Proposition~\ref{prop:blocklaws}\ref{P:entropy}, and $|J|>(4L)^{k'}$, so $\phi_{k'}(z,h)\le\log\frac43\,(4L)^{-k'}$. In general, $$\delta^{z,h}_J\le\sum_jH_{\P^{z,h}_J}(\scrP_{J_j})\le\sum_j|J_j|\,q\log2=|J|\,q\log2$$ by \eqref{eq:shannon}, so $\phi_{k'}\le q\log2$. Hence, by Proposition~\ref{prop:blocklaws}\ref{P:typ},
$$
\sum_{k'=1}^{k}\zeta(\phi_{k'})\le\log\frac43\sum_{k'\ge1}(4L)^{-k'}+q\log2\sum_{k'\ge1}\zeta(\calX\setminus\calX_{k'})\le\frac{\log\frac43}{4L-1}+\frac18\,q\log2 ,
$$
which gives the claim.
\end{stepproof}

\begin{step}[2]
We show that
$$
h_{\widehat\mu_0}(\widehat\sigma\mid\calX)\ge\frac78\,q\log2-\frac{\log\frac43}{4L-1}.
$$
\end{step}

\begin{stepproof}
Let $n\ge1$ and choose $k\ge1$ with $(4L)^k\ge2n+1$. Since $[0,n-1]\subset[-n,n]$, Lemma~\ref{lem:sample}\ref{sm:stab} gives, for every $(z,h)\in\calX$ and every $w\in\calA_\aux^{[0,n-1]}$,
$$
\mu_{z,h}\bigl(\{\omega\ ;\ \omega[0,n-1]=w\}\bigr)=\mu^{(k)}_{z,h}\bigl(\{\omega\ ;\ \omega[0,n-1]=w\}\bigr).
$$
Hence $H_{\mu_{z,h}}(\scrP^n)=H_{\mu^{(k)}_{z,h}}(\scrP^n)$. By \eqref{eq:relsample}, \eqref{eq:relentropy} and Step~1,
$$
\frac1nH_{\widehat\mu_0}\bigl(\scrP^n_\aux\mid\pi_\calX^{-1}\calB_\calX\bigr)=\frac1nH_{\widehat\mu^{(k)}}\bigl(\scrP^n_\aux\mid\pi_\calX^{-1}\calB_\calX\bigr)\ge h_{\widehat\mu^{(k)}}(\widehat\sigma\mid\calX)\ge\frac78\,q\log2-\frac{\log\frac43}{4L-1}.
$$
Since $n$ is arbitrary, \eqref{eq:relentropy} gives the claim.
\end{stepproof}

\begin{step}[3]
Finally we prove the proposition.
\end{step}

\begin{stepproof}
Let $n\ge1$. The partition $\scrP^n$ of $Y$ into the cylinders over $[0,n-1]$ refines the partition of $Y$ according to $\omega[0,n-1]$, whose preimage under $\widehat\pi$ is $\scrP^n_\aux$. Since $\mu_0=\widehat\pi_*\widehat\mu_0$, \eqref{eq:shannon} gives
$$
H_{\mu_0}(\scrP^n)\ge H_{\widehat\mu_0}(\scrP^n_\aux)\ge H_{\widehat\mu_0}\bigl(\scrP^n_\aux\mid\pi_\calX^{-1}\calB_\calX\bigr).
$$
Dividing by $n$, letting $n\to\infty$ and using Step~2,
$$
h_{\mu_0}(\sigma)\ge h_{\widehat\mu_0}(\widehat\sigma\mid\calX)\ge\frac78\,q\log2-\frac{\log\frac43}{4L-1}.
$$
\end{stepproof}
\end{proof}

We now show that every ergodic $\sigma$-invariant probability $\rho\ne\mu_0$ on $Y$ has entropy at most $(2+\frac q8)\log2$.

\begin{proposition}
\label{prop:upper}
Let $\rho$ be an ergodic $\sigma$-invariant Borel probability on $Y$ with $\rho\ne\mu_0$. Then
$$
h_\rho(\sigma)\le2\log2+\frac{9R}{2L}\,q\log2\le\Bigl(2+\frac q8\Bigr)\log2 .
$$
\end{proposition}

\begin{proof}
We divide the proof into four steps.

\begin{step}[1]
We show that there is an ergodic $\widehat\sigma$-invariant Borel probability $\widehat\rho$ on $\widehat Y$ with $\widehat\pi_*\widehat\rho=\rho$.
\end{step}

\begin{stepproof}
The map $\widehat\pi:\widehat Y\to Y$ is continuous and onto, and $\widehat\pi\circ\widehat\sigma=\sigma\circ\widehat\pi$. Hence there is a Borel probability $\lambda_0$ on $\widehat Y$ with $\widehat\pi_*\lambda_0=\rho$, obtained by extending the functional $f\circ\widehat\pi\mapsto\rho(f)$, $f\in\mathcal C^0(Y)$, with the Hahn--Banach theorem and applying the Riesz representation theorem \cite[Theorem 2.14]{Rudin1987}.

By the Krylov--Bogolyubov theorem \cite[Theorem 6.9]{Walters1982}, every weak$^*$ limit point $\lambda$ of the probabilities $\frac1N\sum_{i=0}^{N-1}\widehat\sigma^i_*\lambda_0$ is $\widehat\sigma$-invariant. Since $\widehat\pi_*\widehat\sigma^i_*\lambda_0=\sigma^i_*\rho=\rho$ for every $i$, we have $\widehat\pi_*\lambda=\rho$.

Let $\lambda=\int\lambda_y\,d\lambda(y)$ be the ergodic decomposition of $\lambda$ (see \cite[Chapter 8]{Walters1982}). Then
$$
\rho=\widehat\pi_*\lambda=\int\widehat\pi_*\lambda_y\,d\lambda(y),
$$
and every $(\sigma,\widehat\pi_*\lambda_y)$ is ergodic, as a factor of the ergodic system $(\widehat\sigma,\lambda_y)$. Since $(\sigma,\rho)$ is ergodic, the uniqueness of the ergodic decomposition gives $\widehat\pi_*\lambda_y=\rho$ for $\lambda$-almost every $y$, and we let $\widehat\rho:=\lambda_y$ for such a $y$.
\end{stepproof}

\begin{step}[2]
Let $\eta:=(\pi_\calX)_*\widehat\rho$. We show that $(\sigma_\calX,\eta)$ is ergodic and $\eta\ne\zeta$, where $\zeta=\nu\times\P_\calH$ (see Section~\ref{sec:base}).
\end{step}

\begin{stepproof}
The system $(\sigma_\calX,\eta)$ is a factor of the ergodic system $(\widehat\sigma,\widehat\rho)$ through $\pi_\calX$, so it is ergodic. If $\eta=\zeta$, then Corollary~\ref{cor:faithful} gives $\widehat\rho=\widehat\mu_0$, so $\rho=\widehat\pi_*\widehat\mu_0=\mu_0$, a contradiction.
\end{stepproof}

\begin{step}[3]
We show that $h_{\widehat\rho}(\widehat\sigma\mid\calX)\le\frac{9R}{2L}\,q\log2$.
\end{step}

\begin{stepproof}
Let $(\widehat\rho_{z,h})_{(z,h)\in\calX}$ be the sample measures of $\widehat\rho$, that is,
$$
\widehat\rho(dz,dh,d\omega)=\widehat\rho_{z,h}(d\omega)\,\eta(dz,dh).
$$
Since $\widehat\rho(\widehat Y)=1$, we have $\widehat\rho_{z,h}(\widehat Y_{z,h})=1$ for $\eta$-almost every $(z,h)$. Fix $k\ge1$, $n\ge1$ and such $(z,h)$.

We count the cylinders of $\scrP^n$ which intersect $\widehat Y_{z,h}$. Let $\calF:=\{B\in\calP_k(h)\ ;\ B\subset[0,n-1]\}$, and split $[0,n-1]$ into the two sets
$$
U:=\bigcup_{B\in\calF}B\quad\text{and}\quad V:=[0,n-1]\setminus U .
$$
By Step~2 of the proof of Theorem~\ref{thm:faithful}, $\#V\le6(4L)^k$. If $\omega\in\widehat Y_{z,h}$, then $\omega[B]\in\supp\P^{z,h}_B$ for every $B\in\calF$, and $\omega[V]\in\calA_\aux^V$. Hence the set $\calC_{z,h}$ of $C\in\scrP^n$ with $C\cap\widehat Y_{z,h}\ne\emptyset$ satisfies
$$
\#\calC_{z,h}\le\prod_{B\in\calF}\#\supp\P^{z,h}_B\cdot\#\calA_\aux^V\le\prod_{B\in\calF}\#\supp\P^{z,h}_B\cdot2^{6q(4L)^k} .
$$

We now bound $\#\supp\P^{z,h}_B$. Define $\chi_k(z,h):=q$ if $(z,h)\in\calX_k$, and $\chi_k(z,h):=\frac{9Rq}{2L}+3(4L)^{1-k}$ otherwise. Let $B\in\calF$. By Proposition~\ref{prop:blocklaws}\ref{P:typ}, either $\sigma_\calX^i(z,h)\in\calX_k$ for every $i\in B$, or $\sigma_\calX^i(z,h)\notin\calX_k$ for every $i\in B$. In the first case $\#\supp\P^{z,h}_B\le\#\calA_\aux^B=2^{q|B|}$. In the second case, by Proposition~\ref{prop:blocklaws}\ref{P:entropy} and $|B|>(4L)^k$,
$$
\log_2\#\supp\P^{z,h}_B\le\frac{9Rq}{2L}|B|+12L\le\Bigl(\frac{9Rq}{2L}+3(4L)^{1-k}\Bigr)|B| .
$$
In both cases
$$
\log_2\#\supp\P^{z,h}_B\le\sum_{i\in B}\chi_k\bigl(\sigma_\calX^i(z,h)\bigr).
$$
Since $\chi_k\ge0$ and the elements of $\calF$ are disjoint subsets of $[0,n-1]$, Jensen's inequality for the concave function $\log$ (as in Section~\ref{sec:entropy-def}) and $\widehat\rho_{z,h}(\widehat Y_{z,h})=1$ give
$$
H_{\widehat\rho_{z,h}}(\scrP^n)\le\log\#\calC_{z,h}\le\Bigl(6q(4L)^k+\sum_{i=0}^{n-1}\chi_k\bigl(\sigma_\calX^i(z,h)\bigr)\Bigr)\log2 .
$$
Integrating with respect to $\eta$, using \eqref{eq:relsample} and the invariance of $\eta$, dividing by $n$ and letting $n\to\infty$,
$$
h_{\widehat\rho}(\widehat\sigma\mid\calX)\le\eta(\chi_k)\log2\le\Bigl(q\,\eta(\calX_k)+\frac{9Rq}{2L}+3(4L)^{1-k}\Bigr)\log2 .
$$
By Step~2 and Proposition~\ref{prop:blocklaws}\ref{P:typ}, $\liminf_{k\to\infty}\eta(\calX_k)=0$, so letting $k\to\infty$ along a sequence with $\eta(\calX_k)\to0$ gives the claim.
\end{stepproof}

\begin{step}[4]
Finally, we prove the proposition.
\end{step}

\begin{stepproof}
Since $\rho$ is a factor of $\widehat\rho$, $h_\rho(\sigma)\le h_{\widehat\rho}(\widehat\sigma)$ (see \cite[Chapter 4]{Walters1982}). By the Abramov--Rokhlin formula \eqref{eq:abramov} for $\pi_\calX$, Lemma~\ref{lem:baseentropy} and Step~3,
$$
h_\rho(\sigma)\le h_{\widehat\rho}(\widehat\sigma)=h_\eta(\sigma_\calX)+h_{\widehat\rho}(\widehat\sigma\mid\calX)\le2\log2+\frac{9R}{2L}\,q\log2 ,
$$
and $\frac{9R}{2L}\le\frac18$ by \eqref{eq:a}.
\end{stepproof}
\end{proof}

We now combine Propositions~\ref{prop:lower} and~\ref{prop:upper} with the results of Sections~\ref{sec:extension} and~\ref{sec:mixing}.

The result is the symbolic version of Theorem~\ref{thm:A}.
\begin{theorem}
\label{thm:Bformal}
Let $q\ge12$ and $\calA:=\{-1,1\}\times\{-1,1\}^q$. Let $Y=\widehat\pi(\widehat Y)\subset\calA^\Z$ be the set of Definition~\ref{def:Y}, let $\mu_0=\widehat\pi_*\widehat\mu_0$ be the probability of Definition~\ref{def:mu0}, and let $\nu$ be the probability on $\Sig$ of Theorem~\ref{thm:BCC}. For $x=(z,\omega)\in\calA^\Z$ write $\pi_1(x):=z$ and $Z_i(\pi_1x)=z_i$. Then the following hold.
\begin{enumerate}[label=\textup{(\roman*)},leftmargin=*]
\item\label{B:law} $Y$ is a subshift of $\calA^\Z$ with $\pi_1(Y)=\Sig$, $\mu_0$ is an ergodic $\sigma$-invariant probability on $\calA^\Z$ with $\mu_0(Y)=1$, and $(\pi_1)_*\mu_0=\nu$. Hence, under $\mu_0$, the process $(Z_i\circ\pi_1)_{i\in\Z}$ has law $\nu$ and satisfies the conclusions \ref{BCC:i}, \ref{BCC:iv} and \ref{BCC:v} of Theorem~\ref{thm:BCC}. In particular, the laws of $\frac1{\sqrt n}\sum_{i=0}^{n-1}Z_i\circ\pi_1$ under $\mu_0$ have no Gaussian weak subsequential limit, and $(\calA^\Z,\sigma,\mu_0)$ is not exponentially mixing of order $2R+2$ on $\calF_\theta(\calA^\Z)$ for any $0<\theta<1$.
\item\label{B:gap} $h_{\mu_0}(\sigma)\ge\frac q2\log2$, and every ergodic $\sigma$-invariant probability $\rho$ on $Y$ with $\rho\ne\mu_0$ satisfies
$$
h_\rho(\sigma)\le h_{\mu_0}(\sigma)-\frac q2\log2 .
$$
\item\label{B:mixing} For every $0<\theta<1$, with the constants $C_\theta$ and $c_\theta$ of Theorem~\ref{thm:mixing},
$$
\Bigl|\mu_0\Bigl(\prod_{s=1}^RF_s\circ\sigma^{t_s}\Bigr)-\prod_{s=1}^R\mu_0(F_s)\Bigr|\le C_\theta e^{-c_\theta\Delta(t_1,\dots,t_R)}\prod_{s=1}^R\norm{F_s}_\theta
$$
for all $F_1,\dots,F_R\in\calF_\theta(\calA^\Z)$ and all integers $0\le t_1<\dots<t_R$.
\item\label{B:unique} Let $\psi:\calA^\Z\to\R$ be continuous with $\sup_Y\psi-\inf_Y\psi<\frac q2\log2$. Then $\mu_0$ is the unique $\sigma$-invariant probability $\rho$ on $Y$ which maximises $h_\rho(\sigma)+\rho(\psi)$, that is, the unique equilibrium state of $\psi|_Y$. In particular, for $\psi\equiv0$, $\mu_0$ is the unique measure of maximal entropy of $Y$.
\end{enumerate}
\end{theorem}

\begin{proof}
We start by proving \ref{B:law}. The first sentence is Lemma~\ref{lem:Yhat} and Corollary~\ref{cor:faithful}. Since $(\pi_1)_*\mu_0=\nu$, the process $(Z_i\circ\pi_1)_{i\in\Z}$ under $\mu_0$ has law $\nu$, and \ref{BCC:i}, \ref{BCC:iv} and \ref{BCC:v} of Theorem~\ref{thm:BCC} hold for it. For $F\in\calF_\theta(\Sig)$ we have $F\circ\pi_1\in\calF_\theta(\calA^\Z)$ with $\norm{F\circ\pi_1}_\theta\le\norm F_\theta$, since $d_\theta(\pi_1x,\pi_1y)\le d_\theta(x,y)$. Hence exponential mixing of order $2R+2$ of $(\sigma,\mu_0)$ on $\calF_\theta(\calA^\Z)$ would give exponential mixing of order $2R+2$ of $(\sigma,\nu)$ on $\calF_\theta(\Sig)$, contradicting Theorem~\ref{thm:BCC}\ref{BCC:v}.

We now prove \ref{B:gap}. By \eqref{eq:a}, $\frac{\log\frac43}{4L-1}\le\log2$. By Proposition~\ref{prop:lower} and $q\ge12$,
$$
h_{\mu_0}(\sigma)\ge\Bigl(\frac78q-1\Bigr)\log2\ge\frac q2\log2 .
$$
Let $\rho\ne\mu_0$ be an ergodic $\sigma$-invariant probability on $Y$. By Proposition~\ref{prop:upper},
$$
h_{\mu_0}(\sigma)-h_\rho(\sigma)\ge\Bigl(\frac78q-1-2-\frac q8\Bigr)\log2=\Bigl(\frac34q-3\Bigr)\log2\ge\frac q2\log2 ,
$$
since $\frac q4\ge3$.

Item \ref{B:mixing} is Theorem~\ref{thm:mixing}.

Finally we prove \ref{B:unique}. By Section~\ref{sec:entropy-def}, equilibrium states of $\psi|_Y$ exist, and almost every ergodic component of an equilibrium state is an equilibrium state. Let $\rho\ne\mu_0$ be an ergodic $\sigma$-invariant probability on $Y$. Since $\sup_Y\psi-\frac q2\log2<\inf_Y\psi$, \ref{B:gap} gives
$$
h_\rho(\sigma)+\rho(\psi)\le h_{\mu_0}(\sigma)-\frac q2\log2+\sup_Y\psi<h_{\mu_0}(\sigma)+\inf_Y\psi\le h_{\mu_0}(\sigma)+\mu_0(\psi),
$$
so $\rho$ is not an equilibrium state. Hence almost every ergodic component of an equilibrium state is $\mu_0$, so every equilibrium state equals $\mu_0$. For $\psi\equiv0$ this says that $\mu_0$ is the unique measure of maximal entropy of $Y$.
\end{proof}

\section{Smooth realisation}
\label{sec:realisation}

In this section we prove Theorem~\ref{thm:A}. Let $M$ be a compact smooth manifold with $\dim M\ge3$ and $\varphi:M\to\R$ continuous. We choose $q\ge12$ with $\frac q2\log2>\sup_M\varphi-\inf_M\varphi$, and we realise the subshift $Y$ of Theorem~\ref{thm:Bformal} as an invariant set of a diffeomorphism of $M$ in two steps. First, an affine horseshoe of the disc conjugates $\calA^\Z$ to an invariant set of a smooth diffeomorphism of the disc (Lemma~\ref{lem:horseshoe}). Second, the thickening lemma and the neutral isolation construction of Lin and Tian (Lemma~\ref{lem:thick} and Proposition~\ref{prop:LT}) embed this diffeomorphism into a diffeomorphism $T$ of $M$ whose other invariant probabilities are supported on fixed points. The entropy gap of Theorem~\ref{thm:Bformal}\ref{B:gap} then makes the image of $\mu_0$ the unique equilibrium state of $(T,\varphi)$, and the remaining properties of Theorem~\ref{thm:A} are transferred from Theorem~\ref{thm:Bformal}.

In this section we make use of the following definition.
\begin{definition}
\label{def:isotopy}
Let $N$ be a smooth manifold, possibly with boundary, and let $U\subset N$. A \emph{smooth isotopy} of $N$ relative to $U$ is a $\mathcal C^\infty$ map
$$
I:(y,s)\in N\times[0,1]\mapsto I_s(y)\in N
$$
such that every $I_s$ is a diffeomorphism of $N$ and $I_s(y)=I_0(y)$ for all $y\in U$ and $s\in[0,1]$. Two diffeomorphisms $F,G$ of $N$ are \emph{smoothly isotopic relative to $U$}, and we say that $G$ is \emph{joined to $F$ by a smooth isotopy relative to $U$}, if there is such an $I$ with $I_0=F$ and $I_1=G$. When $U=\emptyset$ we simply say that $F$ and $G$ are \emph{isotopic}.
\end{definition}

\subsection{A horseshoe}
We code $\calA^\Z$ by an affine horseshoe of the disc. Then H\"older functions on the manifold pull back to H\"older functions on $\calA^\Z$.
 
\begin{lemma}
\label{lem:horseshoe}
Let $\mathbb D$ be the closed unit disc, let $N\ge2$ and $\lambda>2N$. There exist a $\mathcal C^\infty$ diffeomorphism $B:\mathbb D\to\mathbb D$, a compact $B$-invariant set $\Lambda_N\subset\operatorname{int}\mathbb D$ and a homeomorphism $\pi:\{1,\dots,N\}^\Z\to\Lambda_N$ with the following properties.
\begin{enumerate}[label=\textup{(\roman*)},leftmargin=*]
\item\label{H:iso} $B$ is equal to the identity on a neighbourhood of $\partial\mathbb D$, and $B$ is joined to the identity by a smooth isotopy relative to a neighbourhood of $\partial\mathbb D$.
\item\label{H:conj} $B\circ\pi=\pi\circ\sigma$.
\item\label{H:est} $|\pi(x)-\pi(y)|\le2\lambda^{-m}$ for every $m\ge0$ and all $x,y\in\{1,\dots,N\}^\Z$ with $x[-m,m]=y[-m,m]$.
\end{enumerate}
\end{lemma}
 
\begin{proof}
This is Smale's horseshoe \cite{Smale1967} with $N$ branches, in the affine model of \cite[Section 2.5c]{KatokHasselblatt1995}. We start by constructing $B$ and proving \ref{H:iso}. In a square $Q_0\subset\operatorname{int}\mathbb D$ of side one, choose $N$ disjoint horizontal strips $H_1,\dots,H_N$ of height $\lambda^{-1}$ and $N$ disjoint vertical strips $V_1,\dots,V_N$ of width $\lambda^{-1}$, which is possible because $N\lambda^{-1}<\frac12$. Let $B$ map each $H_a$ affinely onto $V_a$, contracting the horizontal direction by $\lambda^{-1}$ and expanding the vertical direction by $\lambda$, with derivative $\operatorname{diag}(\lambda^{-1},\lambda)$ or its product with the rotation by $\pi$, so that $B$ preserves orientation. Such a map extends to a $\mathcal C^\infty$ diffeomorphism of $\mathbb D$ equal to the identity near $\partial\mathbb D$ (see \cite[Appendix A.1]{Buzzi2014} and \cite[Lemma 2.2]{LinTian2026} for the disc-supported realisation). Every orientation-preserving diffeomorphism of $\mathbb D$ equal to the identity near $\partial\mathbb D$ is smoothly isotopic to the identity relative to a neighbourhood of the boundary, by Smale's theorem on the contractibility of this group \cite{Smale1959}.
 
We now prove \ref{H:conj}. The set $\Lambda_N:=\bigcap_{n\in\Z}B^{-n}(H_1\cup\dots\cup H_N)$ is coded by itineraries, that is, $\pi(x)$ is the unique point with $B^n(\pi(x))\in H_{x_n}$ for all $n$, and then $B(\pi(x))$ is the unique point whose itinerary is $\sigma x$.
 
Finally we prove \ref{H:est}. If $x[-m,m]=y[-m,m]$, then $\pi(x)$ and $\pi(y)$ lie in $\bigcap_{|n|\le m}B^{-n}(H_{x_n})$, a rectangle of width at most $\lambda^{-m}$ and height at most $\lambda^{-m-1}$, which gives the estimate.
\end{proof}
 
\subsection{Neutral isolation}
The horseshoe carries many invariant measures besides those on the coded copy of $Y$. Neutral isolation removes them. It embeds the disc, thickened to a disc of dimension $\dim M-1$, in a solid torus in such a way that every other invariant measure lives on fixed points. We use the following two results of Lin and Tian, where a closed smooth $m$-disc is a smooth manifold with boundary diffeomorphic to the closed unit ball of $\R^m$.

\begin{lemma}[{Lin--Tian \cite[Lemma 3.1]{LinTian2026}}]
\label{lem:thick}
Let $m\ge1$ and $r\ge0$ be integers, let $\mathcal D$ be a closed smooth $m$-disc, and let $(I_s)_{s\in[0,1]}$ be a smooth isotopy of $\mathcal D$ relative to $\partial \mathcal D$ with $I_0=\id$ (see Definition~\ref{def:isotopy}). There are a closed smooth $(m+r)$-disc $\mathcal D^\sharp$, a smooth embedding $\kappa:\mathcal D\to \mathcal D^\sharp$ with $\kappa(\operatorname{int}\mathcal D)\subset\operatorname{int}\mathcal D^\sharp$, and a smooth isotopy $(\bar I_s)_{s\in[0,1]}$ of $\mathcal D^\sharp$ relative to $\partial \mathcal D^\sharp$ with $\bar I_0=\id$ and $\bar I_s\circ\kappa=\kappa\circ I_s$ for every $s\in[0,1]$.
\end{lemma}

\begin{proposition}[{Lin--Tian \cite[Proposition 3.3]{LinTian2026}}]
\label{prop:LT}
Let $m\ge2$, let $\mathcal D$ be a closed smooth $m$-disc, let $B:\mathcal D\to \mathcal D$ be a $\mathcal C^\infty$ diffeomorphism joined to the identity by a smooth isotopy relative to $\partial \mathcal D$, and let $Q\subset\operatorname{int}\mathcal D$ be a nonempty compact $B$-invariant set such that $B|_Q$ has no fixed point. Define $K_0:=Q\times\{\pi/2\}$. There is a $\mathcal C^\infty$ diffeomorphism $H$ of $\mathcal D\times S^1$ with the following properties.
\begin{enumerate}[label=\textup{(\roman*)},leftmargin=*]
\item\label{LT:diffeo} $H$ is equal to the identity on a neighbourhood of $\partial \mathcal D\times S^1$, and $H$ is smoothly isotopic to the identity relative to such a neighbourhood.
\item\label{LT:B} $H(y,\pi/2)=(B(y),\pi/2)$ for every $y\in Q$, and $K_0\cap\operatorname{Fix}(H)=\emptyset$.
\item\label{LT:inv} Every $H$-invariant Borel probability $m$ on $\mathcal D\times S^1$ satisfies $m\bigl(\operatorname{Fix}(H)\cup K_0\bigr)=1$.
\end{enumerate}
\end{proposition}
 
\subsection{Proof of Theorem~\ref{thm:A}}

We deduce Theorem~\ref{thm:A} from the following more precise statement.

\begin{theorem}[Theorem~\ref{thm:A}, precise form]
\label{thm:Aprecise}
Let $M$ be a compact smooth manifold of dimension at least three, let $\varphi:M\to\R$ be continuous, and let $R\ge2$ be an integer. There exist a $\mathcal C^\infty$ diffeomorphism $T:M\to M$, isotopic to the identity and equal to the identity outside a coordinate ball, and $f\in \mathcal C^\infty(M,[-1,1])$ with the following properties.
\begin{enumerate}[label=\textup{(\roman*)},leftmargin=*]
\item\label{A:unique} $(T,\varphi)$ has a unique equilibrium state $\mu$.
\item\label{A:every} The stationary process $(f\circ T^n)_{n\in\Z}$ under $\mu$ has law $\nu$. In particular $\mu(f)=0$, the process $(f\circ T^n)_{n\in\Z}$ is a Rademacher process and is $(2R+1)$-tuplewise independent under $\mu$, and the laws of $S_nf/\sqrt n$ under $\mu$ have no Gaussian weak subsequential limit.
\item\label{A:nomixing} There are $c,C>0$ and $N_0$ such that for every $N\ge N_0$ there are integers $0\le n_1<\dots<n_{2R+2}<N$ with
$$
\Delta(n_1,\dots,n_{2R+2})\ge\frac NC\quad\text{and}\quad-\mu\Bigl(\prod_{j=1}^{2R+2}f\circ T^{n_j}\Bigr)\ge cN^{-R-1}.
$$
Moreover, $(T,\mu)$ is not exponentially mixing of order $R+1$ for H\"older observables.
\item\label{A:mixing} $(T,\mu)$ is exponentially mixing of order $R$ for H\"older observables, that is, for every $0<\beta\le1$ there are $C_\beta,c_\beta>0$ with
$$
\Bigl|\mu\Bigl(\prod_{j=1}^R\phi_j\circ T^{n_j}\Bigr)-\prod_{j=1}^R\mu(\phi_j)\Bigr|\le C_\beta e^{-c_\beta\Delta(n_1,\dots,n_R)}\prod_{j=1}^R\norm{\phi_j}_{\mathcal C^\beta}
$$
for all $\phi_1,\dots,\phi_R\in \mathcal C^\beta(M)$ and all integers $0\le n_1<\dots<n_R$.
\end{enumerate}
\end{theorem}

\begin{proof}[Proof of Theorem~\ref{thm:A}]
Let $T$, $f$ and $\mu$ be given by Theorem~\ref{thm:Aprecise}. Then $\mu$ is the unique equilibrium state of $(T,\varphi)$, $\mu(f)=0$ and \ref{A:clt} hold by \ref{A:every}, \ref{A:nomix} holds by \ref{A:nomixing}, and \ref{A:mix} holds by \ref{A:mixing}.
\end{proof}

\begin{proof}[Proof of Theorem~\ref{thm:Aprecise}]
Let $M$, $\varphi$ and $R$ be as in Theorem~\ref{thm:Aprecise}. We choose an integer $q\ge12$ such that the entropy gap $\frac q2\log2$ of Theorem~\ref{thm:Bformal}\ref{B:gap} satisfies
$$
\frac q2\log2>\sup_M\varphi-\inf_M\varphi .
$$
Let $\calA:=\{-1,1\}\times\{-1,1\}^q$, and let $Y\subset\calA^\Z$ and $\mu_0$ be the subshift and the probability of Definitions~\ref{def:Y} and~\ref{def:mu0} for this $q$. We divide the proof into five steps.

\begin{step}[1]
We construct a $\mathcal C^\infty$ diffeomorphism $T:M\to M$, constants $\lambda>1$ and $C_{\rm geo}>0$, a continuous map $\iota:\calA^\Z\to M$ and a compact set $P_M\subset M$ with the following properties.
\begin{itemize}[leftmargin=*]
\item $T$ is isotopic to the identity (see Definition~\ref{def:isotopy}) and equal to the identity outside a coordinate ball.
\item $T\circ\iota=\iota\circ\sigma$ on $Y$, and $\iota|_Y$ is a homeomorphism onto its image $K:=\iota(Y)$.
\item $P_M$ consists of fixed points of $T$, $P_M\cap K=\emptyset$, and every $T$-invariant Borel probability on $M$ is supported on $K\cup P_M$.
\item For all $m\ge0$ and $x,y\in\calA^\Z$,
\begin{equation}
\label{eq:geo}
d_M\bigl(\iota(x),\iota(y)\bigr)\le C_{\rm geo}\,\lambda^{-m}\quad\text{whenever }x[-m,m]=y[-m,m].
\end{equation}
\end{itemize}
\end{step}

\begin{stepproof}
The alphabet $\calA$ has $N:=2^{1+q}$ symbols, and we identify it with $\{1,\dots,N\}$. Choose $\lambda>2N$, and let $B$, $\Lambda_N$ and $\pi:\calA^\Z\to\Lambda_N$ be given by Lemma~\ref{lem:horseshoe}. Then $Q:=\pi(Y)$ is a compact $B$-invariant set, and by Lemma~\ref{lem:horseshoe}\ref{H:conj}, $B|_Q$ is conjugate by $\pi$ to $\sigma|_Y$, which has no fixed point by Lemma~\ref{lem:nofixed}.

We now thicken $\mathbb D$. By Lemma~\ref{lem:horseshoe}\ref{H:iso}, there is a smooth isotopy $(I_s)_{s\in[0,1]}$ of $\mathbb D$ relative to $\partial\mathbb D$ with $I_0=\id$ and $I_1=B$. Lemma~\ref{lem:thick} with $m=2$ and $r=\dim M-3$ gives a closed smooth $(\dim M-1)$-disc $\mathcal D^\sharp$, an embedding $\kappa:\mathbb D\to\mathcal D^\sharp$ and an isotopy $(\bar I_s)_{s\in[0,1]}$ of $\mathcal D^\sharp$ relative to $\partial\mathcal D^\sharp$ with $\bar I_s\circ\kappa=\kappa\circ I_s$. Let $\bar B:=\bar I_1$ and $\bar Q:=\kappa(Q)$. Then $\bar B\circ\kappa=\kappa\circ B$, so $\bar Q\subset\operatorname{int}\mathcal D^\sharp$ is a compact $\bar B$-invariant set, and $\bar B|_{\bar Q}$ is conjugate by $\kappa$ to $B|_Q$, so it has no fixed point. Let $H$ and $K_0=\bar Q\times\{\pi/2\}$ be given by Proposition~\ref{prop:LT} with $m=\dim M-1$, applied to $\bar B$ and $\bar Q$.

Fix a coordinate ball $U\subset M$ and a smooth embedding $\Psi:\mathcal D^\sharp\times S^1\to U$, which exists since the product of the closed unit ball of $\R^{\dim M-1}$ with $S^1$ embeds in $\R^{\dim M}$. Let $V:=\Psi(\mathcal D^\sharp\times S^1)$, and define
$$
T:=\Psi\circ H\circ\Psi^{-1}\ \text{on }V\quad\text{and}\quad T:=\id\ \text{on }M\setminus V,
$$
and
$$
\iota:x\in\calA^\Z\mapsto\Psi\bigl(\kappa(\pi(x)),\pi/2\bigr)\in M .
$$

Since $H$ is equal to the identity near $\partial\mathcal D^\sharp\times S^1$ and is smoothly isotopic to the identity relative to such a neighbourhood by Proposition~\ref{prop:LT}\ref{LT:diffeo}, $T$ is a $\mathcal C^\infty$ diffeomorphism, isotopic to the identity, and equal to the identity outside the coordinate ball $U$. By Proposition~\ref{prop:LT}\ref{LT:B}, $\bar B\circ\kappa=\kappa\circ B$ and Lemma~\ref{lem:horseshoe}\ref{H:conj}, $T\circ\iota=\iota\circ\sigma$ on $Y$, and $\iota|_Y$ is a continuous injection of a compact space, hence a homeomorphism onto $K=\Psi(K_0)$. The set $P_M:=\Psi(\operatorname{Fix}(H))\cup(M\setminus V)$ consists of fixed points of $T$, is compact, and is disjoint from $K$ by Proposition~\ref{prop:LT}\ref{LT:B}, and every $T$-invariant probability is supported on $K\cup P_M$ by Proposition~\ref{prop:LT}\ref{LT:inv}. Finally, let $C_{\rm geo}:=2\Lip$, where $\Lip<\infty$ is the Lipschitz constant of the smooth map $(y,\theta)\mapsto\Psi(\kappa(y),\theta)$ from $\mathbb D\times S^1\subset\R^2\times\R^2$, with the Euclidean metric, to $(M,d_M)$. Then \eqref{eq:geo} follows from Lemma~\ref{lem:horseshoe}\ref{H:est}.
\end{stepproof}

\begin{step}[2]
We show that $\mu_0$ is the unique equilibrium state of $\psi|_Y$, where $\psi:=\varphi\circ\iota:\calA^\Z\to\R$.
\end{step}

\begin{stepproof}
The function $\psi$ is continuous, and by the choice of $q$,
$$
\sup_{\calA^\Z}\psi-\inf_{\calA^\Z}\psi\le\sup_M\varphi-\inf_M\varphi<\frac q2\log2 ,
$$
so the claim is Theorem~\ref{thm:Bformal}\ref{B:unique}.
\end{stepproof}

\begin{step}[3]
We prove \ref{A:unique}.
\end{step}

\begin{stepproof}
Let $m$ be a $T$-invariant probability on $M$. By Step~1, $m$ is supported on $K\cup P_M$, so $m=t\,\iota_*\rho+(1-t)m'$ with $t\in[0,1]$, $\rho$ shift-invariant on $Y$ and $m'$ supported on $P_M$. Entropy is affine, $h_{\iota_*\rho}(T)=h_\rho(\sigma)$ by the conjugacy, $h_{m'}(T)=0$, and $\int\varphi\,d\iota_*\rho=\rho(\psi)$. Hence
$$
h_m(T)+m(\varphi)=t\bigl(h_\rho(\sigma)+\rho(\psi)\bigr)+(1-t)\,m'(\varphi).
$$
By Theorem~\ref{thm:Bformal}\ref{B:gap} and the choice of $q$,
$$
P_Y(\psi)\ge h_{\mu_0}(\sigma)+\mu_0(\psi)\ge\frac q2\log2+\min_M\varphi>\max_M\varphi\ge m'(\varphi).
$$
Therefore $P_T(\varphi)=P_Y(\psi)$, and $m$ is an equilibrium state of $(T,\varphi)$ exactly when $t=1$ and $\rho$ is an equilibrium state of $\psi|_Y$, that is, by Step~2, exactly when $m=\iota_*\mu_0$. Hence $\mu:=\iota_*\mu_0$ is the unique equilibrium state of $(T,\varphi)$.
\end{stepproof}

\begin{step}[4]
We prove \ref{A:every} and \ref{A:nomixing}.
\end{step}

\begin{stepproof}
The sets $\iota(\{x\in Y\ ;\ z_0=1\})$ and $\iota(\{x\in Y\ ;\ z_0=-1\})$ are disjoint compact subsets of $M$. Choose $f\in \mathcal C^\infty(M,[-1,1])$ equal to $1$ on a neighbourhood of the first and to $-1$ on a neighbourhood of the second. Then $f\circ T^n\circ\iota(x)=z_n$ for $x\in Y$ and $n\in\Z$. By Theorem~\ref{thm:Bformal}\ref{B:law}, $(z_n)_{n\in\Z}$ has law $\nu$ under $\mu_0$, so $(f\circ T^n)_{n\in\Z}$ has law $\nu$ under $\mu$, and the properties listed in \ref{A:every} and \ref{A:nomixing} are those of Theorem~\ref{thm:BCC}\ref{BCC:i}, \ref{BCC:iv} and \ref{BCC:v}.

Assume that $(T,\mu)$ is exponentially mixing of order $R+1$ for H\"older observables. Fix $0<\beta\le1$, and let $C_{R+1,\beta},c_{R+1,\beta}>0$ be the constants of \eqref{eq:allorders} with $s=R+1$. The function $\phi:=f\cdot(f\circ T)$ is smooth, so $\phi\in\mathcal C^\beta(M)$. Since $(f\circ T^n)_{n\in\Z}$ has law $\nu$ under $\mu$,
$$
\mu(\phi)=\E_\nu[Z_0Z_1]=0\quad\text{and}\quad\mu\Bigl(\prod_{\ell=1}^{R+1}\phi\circ T^{t_\ell}\Bigr)=\E_\nu\Bigl[\prod_{\ell=1}^{R+1}Z_{t_\ell}Z_{t_\ell+1}\Bigr]
$$
for all integers $0\le t_1<\dots<t_{R+1}$. Applying the estimate to $\phi_1=\dots=\phi_{R+1}=\phi$ at the times of Lemma~\ref{lem:nomixR1} gives
$$
144^{-R-1}(8R+8)^{-2(R+1)k}\le C_{R+1,\beta}\,e^{-2c_{R+1,\beta}(8R+8)^{k-1}}\norm\phi_{\mathcal C^\beta}^{R+1}
$$
for every $k\ge2$, which is impossible for large $k$.
\end{stepproof}

\begin{step}[5]
Finally we prove \ref{A:mixing}.
\end{step}

\begin{stepproof}
For $\phi\in \mathcal C^\beta(M)$, \eqref{eq:geo} gives
$$
|\phi\circ\iota(x)-\phi\circ\iota(y)|\le\norm\phi_{\mathcal C^\beta}C_{\rm geo}^\beta\lambda^{-\beta m}\quad\text{when }x[-m,m]=y[-m,m],
$$
so $\phi\circ\iota\in\calF_{\theta}(\calA^\Z)$ with $\theta:=\lambda^{-\beta}$ and $\norm{\phi\circ\iota}_\theta\le(1+C_{\rm geo}^\beta\theta^{-1})\norm\phi_{\mathcal C^\beta}$. Since $\mu(\prod_j\phi_j\circ T^{n_j})=\mu_0(\prod_j\phi_j\circ\iota\circ\sigma^{n_j})$, Theorem~\ref{thm:Bformal}\ref{B:mixing} gives \ref{A:mixing}.
\end{stepproof}
\end{proof}

\begin{proof}[Proof of Corollary~\ref{cor:mme}]
Apply Theorem~\ref{thm:A} to $\varphi\equiv0$, whose equilibrium states are the measures of maximal entropy.
\end{proof}

\section{Lean formalisation}
\label{sec:lean}
Theorems~\ref{thm:A}, \ref{thm:Bformal} and~\ref{thm:Aprecise} and Corollary~\ref{cor:mme} have been formalised in Lean 4.35.0 using Mathlib. The source code is available at \textbf{https://github.com/Manzatto/equilibrium-states-finite-order-mixing-clt-counterexample}. The results are collected in the file \texttt{MainTheorems.lean} of the folder \texttt{FiniteOrderEq/Unique}, in the namespace \texttt{FiniteOrderEq.Unique}. The main declarations are the following.
\begin{itemize}
\item Theorem~\ref{thm:A}: \texttt{TheoremA}.
\item Corollary~\ref{cor:mme}: \texttt{Corollary1\_1}.
\item Theorem~\ref{thm:Bformal}: \texttt{Theorem7\_3}.
\item Theorem~\ref{thm:Aprecise}: \texttt{Theorem8\_5}.
\item Lemma~\ref{lem:holderC1} \texttt{Lemma2\_1}.
\item Lemma~\ref{lem:nomixR1}: \texttt{Lemma3\_5}.
\end{itemize}
Every result cited in this article and used in the proofs above is also formalised. This includes the bounded differences inequality of McDiarmid \cite{McDiarmid1989} used in Section~\ref{sec:products}, the horseshoe of Smale \cite{Smale1967}, and Lemma~3.1 and Proposition~3.3 of Lin and Tian \cite{LinTian2026} in every dimension. The process of Theorem~\ref{thm:BCC} and the marker sets of Proposition~\ref{prop:markers} are taken from the formalisation of \cite{BassolsCastro2026}, and the Kolmogorov--Sinai entropy, the pressure and the variational principle from the formalisation of \cite{BassolsCastro2025}. The correspondence between the article and the formal statements is documented in \texttt{docs/StatementAudit.md}.

The development contains no \texttt{sorry} placeholders or additional axioms. An audit of the project declarations verifies that their axiom dependencies are confined to \texttt{propext}, \texttt{Classical.choice} and \texttt{Quot.sound}. The compiler and dependency versions are fixed by \texttt{lean-toolchain} and \texttt{lake-manifest.json}. Compilation and axiom dependencies can be checked from the repository root by running
\begin{verbatim}
lake exe cache get
lake build
lake env lean Audit.lean
\end{verbatim}

\vspace{1em}
\paragraph{\textbf{Statement on AI use.}} OpenAI's ChatGPT 6.0 (Astra) was used as a research assistant during the development of this work.
AI assistance was also used to develop and revise the Lean formalisation of the main results. The resulting proofs were checked by Lean's kernel, with the axiom dependencies described above.
The author takes full responsibility for the correctness of all mathematical statements and proofs, the accuracy of the references, and the final content of the manuscript.

\section*{Acknowledgements}
  MMC is thanks Bernat Bassols Cornudella for the useful discussions. MMC is supported by the S\~ao Paulo Research Foundation (FAPESP, grant no. 2025/26997-9).
 
\bibliographystyle{abbrvnat}
\bibliography{refs}

@article{AbramovRokhlin1962,
  author  = {Abramov, Leonid M. and Rokhlin, Vladimir A.},
  title   = {The entropy of a skew product of measure-preserving transformations},
  journal = {Vestnik Leningradskogo Universiteta},
  volume  = {17},
  number  = {7},
  pages   = {5--13},
  year    = {1962},
  note    = {In Russian; English translation in American Mathematical Society Translations, Series 2, 48 (1966), 255--265}
}

@misc{BassolsCastro2026,
  author        = {Bassols Cornudella, Bernat and Castro, Matheus M.},
  title         = {Finite-order exponential mixing does not imply the central limit theorem},
  year          = {2026},
  eprint        = {2610.17162},
  archiveprefix = {arXiv},
  url           = {https://doi.org/10.48550/arXiv.2610.17162},
  note          = {arXiv preprint}
}

@article{BenOvadia2018,
  author  = {Ben Ovadia, Snir},
  title   = {Symbolic dynamics for non-uniformly hyperbolic diffeomorphisms of compact smooth manifolds},
  journal = {Journal of Modern Dynamics},
  volume  = {13},
  pages   = {43--113},
  year    = {2018},
  doi     = {10.3934/jmd.2018013}
}

@article{BjorklundGorodnik2020,
  author  = {Bj{\"o}rklund, Michael and Gorodnik, Alexander},
  title   = {Central limit theorems for group actions which are exponentially mixing of all orders},
  journal = {Journal d'Analyse Math\'ematique},
  volume  = {141},
  number  = {2},
  pages   = {457--482},
  year    = {2020},
  doi     = {10.1007/s11854-020-0106-7}
}

@article{Bowen1974,
  author  = {Bowen, Rufus},
  title   = {Some systems with unique equilibrium states},
  journal = {Mathematical Systems Theory},
  volume  = {8},
  number  = {3},
  pages   = {193--202},
  year    = {1974},
  doi     = {10.1007/BF01762666}
}

@book{Bowen1975,
  author    = {Bowen, Rufus},
  title     = {Equilibrium States and the Ergodic Theory of {A}nosov Diffeomorphisms},
  series    = {Lecture Notes in Mathematics},
  volume    = {470},
  publisher = {Springer},
  address   = {Berlin},
  year      = {1975},
  doi       = {10.1007/BFb0081279}
}

@article{Bradley1989,
  author  = {Bradley, Richard C.},
  title   = {A stationary, pairwise independent, absolutely regular sequence for which the central limit theorem fails},
  journal = {Probability Theory and Related Fields},
  volume  = {81},
  number  = {1},
  pages   = {1--10},
  year    = {1989},
  doi     = {10.1007/BF00343735}
}

@article{Bradley2007,
  author  = {Bradley, Richard C.},
  title   = {On a stationary, triple-wise independent, absolutely regular counterexample to the central limit theorem},
  journal = {Rocky Mountain Journal of Mathematics},
  volume  = {37},
  number  = {1},
  pages   = {25--44},
  year    = {2007},
  doi     = {10.1216/rmjm/1181069318}
}

@article{Bradley2010,
  author  = {Bradley, Richard C.},
  title   = {A strictly stationary, ``causal,'' 5-tuplewise independent counterexample to the central limit theorem},
  journal = {ALEA. Latin American Journal of Probability and Mathematical Statistics},
  volume  = {7},
  pages   = {377--450},
  year    = {2010}
}

@article{BradleyPruss2009,
  author  = {Bradley, Richard C. and Pruss, Alexander R.},
  title   = {A strictly stationary, $N$-tuplewise independent counterexample to the central limit theorem},
  journal = {Stochastic Processes and their Applications},
  volume  = {119},
  number  = {10},
  pages   = {3300--3318},
  year    = {2009},
  doi     = {10.1016/j.spa.2009.05.009}
}

@article{Buzzi2014,
  author  = {Buzzi, J{\'e}r{\^o}me},
  title   = {$C^r$ surface diffeomorphisms with no maximal entropy measure},
  journal = {Ergodic Theory and Dynamical Systems},
  volume  = {34},
  number  = {6},
  pages   = {1770--1793},
  year    = {2014},
  doi     = {10.1017/etds.2013.25}
}

@article{BuzziCrovisierSarig2022,
  author  = {Buzzi, J{\'e}r{\^o}me and Crovisier, Sylvain and Sarig, Omri},
  title   = {Measures of maximal entropy for surface diffeomorphisms},
  journal = {Annals of Mathematics},
  volume  = {195},
  number  = {2},
  pages   = {421--508},
  year    = {2022},
  doi     = {10.4007/annals.2022.195.2.2}
}

@misc{BuzziCrovisierSarig2025,
  author        = {Buzzi, J{\'e}r{\^o}me and Crovisier, Sylvain and Sarig, Omri},
  title         = {Strong positive recurrence and exponential mixing for diffeomorphisms},
  year          = {2025},
  eprint        = {2501.07455},
  archiveprefix = {arXiv},
  url           = {https://doi.org/10.48550/arXiv.2501.07455},
  note          = {arXiv preprint}
}

@article{ClimenhagaFisherThompson2018,
  author  = {Climenhaga, Vaughn and Fisher, Todd and Thompson, Daniel J.},
  title   = {Unique equilibrium states for {B}onatti--{V}iana diffeomorphisms},
  journal = {Nonlinearity},
  volume  = {31},
  number  = {6},
  pages   = {2532--2570},
  year    = {2018},
  doi     = {10.1088/1361-6544/aab1cd}
}

@article{ClimenhagaThompson2016,
  author  = {Climenhaga, Vaughn and Thompson, Daniel J.},
  title   = {Unique equilibrium states for flows and homeomorphisms with non-uniform structure},
  journal = {Advances in Mathematics},
  volume  = {303},
  pages   = {745--799},
  year    = {2016},
  doi     = {10.1016/j.aim.2016.07.029}
}

@article{DenkerPhilipp1984,
  author  = {Denker, Manfred and Philipp, Walter},
  title   = {Approximation by {B}rownian motion for {G}ibbs measures and flows under a function},
  journal = {Ergodic Theory and Dynamical Systems},
  volume  = {4},
  number  = {4},
  pages   = {541--552},
  year    = {1984},
  doi     = {10.1017/S0143385700002637}
}

@article{Dolgopyat2004,
  author  = {Dolgopyat, Dmitry},
  title   = {Limit theorems for partially hyperbolic systems},
  journal = {Transactions of the American Mathematical Society},
  volume  = {356},
  number  = {4},
  pages   = {1637--1689},
  year    = {2004},
  doi     = {10.1090/S0002-9947-03-03335-X}
}

@article{DolgopyatDongKanigowskiNandori2022,
  author  = {Dolgopyat, Dmitry and Dong, Changguang and Kanigowski, Adam and N{\'a}ndori, P{\'e}ter},
  title   = {Flexibility of statistical properties for smooth systems satisfying the central limit theorem},
  journal = {Inventiones Mathematicae},
  volume  = {230},
  pages   = {31--120},
  year    = {2022},
  doi     = {10.1007/s00222-022-01121-0}
}

@article{DolgopyatKanigowskiRodriguezHertz2024,
  author  = {Dolgopyat, Dmitry and Kanigowski, Adam and Rodriguez Hertz, Federico},
  title   = {Exponential mixing implies {B}ernoulli},
  journal = {Annals of Mathematics},
  volume  = {199},
  number  = {3},
  pages   = {1225--1292},
  year    = {2024},
  doi     = {10.4007/annals.2024.199.3.5}
}

@book{Downarowicz2011,
  author    = {Downarowicz, Tomasz},
  title     = {Entropy in Dynamical Systems},
  series    = {New Mathematical Monographs},
  volume    = {18},
  publisher = {Cambridge University Press},
  address   = {Cambridge},
  year      = {2011}
}

@article{Gordin1969,
  author  = {Gordin, Mikhail I.},
  title   = {The central limit theorem for stationary processes},
  journal = {Soviet Mathematics Doklady},
  volume  = {10},
  pages   = {1174--1176},
  year    = {1969},
  note    = {Russian original: Doklady Akademii Nauk SSSR 188 (1969), no. 4, 739--741}
}

@article{GorodnikSpatzier2014,
  author  = {Gorodnik, Alexander and Spatzier, Ralf},
  title   = {Exponential mixing of nilmanifold automorphisms},
  journal = {Journal d'Analyse Math\'ematique},
  volume  = {123},
  pages   = {355--396},
  year    = {2014},
  doi     = {10.1007/s11854-014-0024-7}
}

@article{GorodnikSpatzier2015,
  author  = {Gorodnik, Alexander and Spatzier, Ralf},
  title   = {Mixing properties of commuting nilmanifold automorphisms},
  journal = {Acta Mathematica},
  volume  = {215},
  number  = {1},
  pages   = {127--159},
  year    = {2015},
  doi     = {10.1007/s11511-015-0130-0}
}

@incollection{Gouezel2015,
  author    = {Gou{\"e}zel, S{\'e}bastien},
  title     = {Limit theorems in dynamical systems using the spectral method},
  booktitle = {Hyperbolic Dynamics, Fluctuations and Large Deviations},
  series    = {Proceedings of Symposia in Pure Mathematics},
  volume    = {89},
  pages     = {161--193},
  publisher = {American Mathematical Society},
  address   = {Providence, RI},
  year      = {2015},
  doi       = {10.1090/pspum/089/01487}
}

@article{GuivarchHardy1988,
  author  = {Guivarc'h, Yves and Hardy, Jean},
  title   = {Th\'eor\`emes limites pour une classe de cha\^{\i}nes de {M}arkov et applications aux diff\'eomorphismes d'{A}nosov},
  journal = {Annales de l'Institut Henri Poincar\'e. Probabilit\'es et Statistiques},
  volume  = {24},
  number  = {1},
  pages   = {73--98},
  year    = {1988}
}

@article{Ibragimov1962,
  author  = {Ibragimov, Ildar A.},
  title   = {Some limit theorems for stationary processes},
  journal = {Theory of Probability and Its Applications},
  volume  = {7},
  number  = {4},
  pages   = {349--382},
  year    = {1962},
  doi     = {10.1137/1107036}
}

@article{Janson1988,
  author  = {Janson, Svante},
  title   = {Some pairwise independent sequences for which the central limit theorem fails},
  journal = {Stochastics},
  volume  = {23},
  number  = {4},
  pages   = {439--448},
  year    = {1988},
  doi     = {10.1080/17442508808833503}
}

@book{KatokHasselblatt1995,
  author    = {Katok, Anatole and Hasselblatt, Boris},
  title     = {Introduction to the Modern Theory of Dynamical Systems},
  series    = {Encyclopedia of Mathematics and its Applications},
  volume    = {54},
  publisher = {Cambridge University Press},
  address   = {Cambridge},
  year      = {1995}
}

@article{LeBorgne1999,
  author  = {Le Borgne, St\'ephane},
  title   = {Limit theorems for non-hyperbolic automorphisms of the torus},
  journal = {Israel Journal of Mathematics},
  volume  = {109},
  pages   = {61--73},
  year    = {1999},
  doi     = {10.1007/BF02775026}
}

@misc{LinTian2026,
  author        = {Lin, Wanshan and Tian, Xueting},
  title         = {Realizing prescribed entropy functions by smooth diffeomorphisms of closed manifolds of dimension at least three},
  year          = {2026},
  eprint        = {2608.15825},
  archiveprefix = {arXiv},
  url           = {https://doi.org/10.48550/arXiv.2608.15825},
  note          = {arXiv preprint}
}

@incollection{Liverani1996,
  author    = {Liverani, Carlangelo},
  title     = {Central limit theorem for deterministic systems},
  booktitle = {International Conference on Dynamical Systems (Montevideo, 1995)},
  series    = {Pitman Research Notes in Mathematics Series},
  volume    = {362},
  pages     = {56--75},
  publisher = {Longman},
  address   = {Harlow},
  year      = {1996}
}

@misc{LuoYang2026,
  author        = {Luo, Chiyi and Yang, Dawei},
  title         = {Strong positive recurrence for potential and exponential mixing of equilibrium states of surface diffeomorphisms},
  year          = {2026},
  eprint        = {2602.04275},
  archiveprefix = {arXiv},
  url           = {https://doi.org/10.48550/arXiv.2602.04275},
  note          = {arXiv preprint}
}

@article{Misiurewicz1973,
  author  = {Misiurewicz, Micha{\l}},
  title   = {Diffeomorphism without any measure with maximal entropy},
  journal = {Bulletin de l'Acad\'emie Polonaise des Sciences. S\'erie des Sciences Math\'ematiques, Astronomiques et Physiques},
  volume  = {21},
  pages   = {903--910},
  year    = {1973}
}

@article{Newhouse1989,
  author  = {Newhouse, Sheldon E.},
  title   = {Continuity properties of entropy},
  journal = {Annals of Mathematics},
  volume  = {129},
  number  = {2},
  pages   = {215--235},
  year    = {1989},
  doi     = {10.2307/1971492}
}

@article{Parry1964,
  author  = {Parry, William},
  title   = {Intrinsic {M}arkov chains},
  journal = {Transactions of the American Mathematical Society},
  volume  = {112},
  pages   = {55--66},
  year    = {1964},
  doi     = {10.1090/S0002-9947-1964-0161372-1}
}

@book{ParryPollicott1990,
  author    = {Parry, William and Pollicott, Mark},
  title     = {Zeta Functions and the Periodic Orbit Structure of Hyperbolic Dynamics},
  series    = {Ast\'erisque},
  volume    = {187--188},
  publisher = {Soci\'et\'e Math\'ematique de France},
  address   = {Paris},
  year      = {1990}
}

@article{Pruss1998,
  author  = {Pruss, Alexander R.},
  title   = {A bounded $N$-tuplewise independent and identically distributed counterexample to the {CLT}},
  journal = {Probability Theory and Related Fields},
  volume  = {111},
  pages   = {323--332},
  year    = {1998},
  doi     = {10.1007/s004400050170}
}

@article{Ratner1973,
  author  = {Ratner, Marina},
  title   = {The central limit theorem for geodesic flows on $n$-dimensional manifolds of negative curvature},
  journal = {Israel Journal of Mathematics},
  volume  = {16},
  pages   = {181--197},
  year    = {1973},
  doi     = {10.1007/BF02757869}
}

@article{Rosenblatt1956,
  author  = {Rosenblatt, Murray},
  title   = {A central limit theorem and a strong mixing condition},
  journal = {Proceedings of the National Academy of Sciences of the United States of America},
  volume  = {42},
  number  = {1},
  pages   = {43--47},
  year    = {1956},
  doi     = {10.1073/pnas.42.1.43}
}

@book{Ruelle1978,
  author    = {Ruelle, David},
  title     = {Thermodynamic Formalism},
  series    = {Encyclopedia of Mathematics and its Applications},
  volume    = {5},
  publisher = {Addison-Wesley},
  address   = {Reading, MA},
  year      = {1978}
}

@article{Sarig1999,
  author  = {Sarig, Omri M.},
  title   = {Thermodynamic formalism for countable {M}arkov shifts},
  journal = {Ergodic Theory and Dynamical Systems},
  volume  = {19},
  number  = {6},
  pages   = {1565--1593},
  year    = {1999},
  doi     = {10.1017/S0143385799146820}
}

@article{Sarig2013,
  author  = {Sarig, Omri M.},
  title   = {Symbolic dynamics for surface diffeomorphisms with positive entropy},
  journal = {Journal of the American Mathematical Society},
  volume  = {26},
  number  = {2},
  pages   = {341--426},
  year    = {2013},
  doi     = {10.1090/S0894-0347-2012-00758-9}
}

@article{Sinai1960,
  author  = {Sinai, Yakov G.},
  title   = {The central limit theorem for geodesic flows on manifolds of constant negative curvature},
  journal = {Doklady Akademii Nauk SSSR},
  volume  = {133},
  number  = {6},
  pages   = {1303--1306},
  year    = {1960},
  note    = {In Russian; English translation in Soviet Mathematics Doklady 1 (1960), 983--987}
}

@article{Sinai1972,
  author  = {Sinai, Yakov G.},
  title   = {Gibbs measures in ergodic theory},
  journal = {Russian Mathematical Surveys},
  volume  = {27},
  number  = {4},
  pages   = {21--69},
  year    = {1972},
  doi     = {10.1070/RM1972v027n04ABEH001383}
}

@article{Smale1959,
  author  = {Smale, Stephen},
  title   = {Diffeomorphisms of the 2-sphere},
  journal = {Proceedings of the American Mathematical Society},
  volume  = {10},
  number  = {4},
  pages   = {621--626},
  year    = {1959},
  doi     = {10.2307/2033664}
}

@article{Smale1967,
  author  = {Smale, Stephen},
  title   = {Differentiable dynamical systems},
  journal = {Bulletin of the American Mathematical Society},
  volume  = {73},
  number  = {6},
  pages   = {747--817},
  year    = {1967},
  doi     = {10.1090/S0002-9904-1967-11798-1}
}

@book{Walters1982,
  author    = {Walters, Peter},
  title     = {An Introduction to Ergodic Theory},
  series    = {Graduate Texts in Mathematics},
  volume    = {79},
  publisher = {Springer-Verlag},
  address   = {New York},
  year      = {1982}
}

@phdthesis{Weakley2013,
  author = {Weakley, Le Mai},
  title  = {Some strictly stationary, $N$-tuplewise independent counterexamples to the central limit theorem},
  school = {Indiana University},
  year   = {2013},
  note   = {ProQuest LLC, Ann Arbor, MI}
}

@article{Yomdin1987,
  author  = {Yomdin, Yosef},
  title   = {Volume growth and entropy},
  journal = {Israel Journal of Mathematics},
  volume  = {57},
  number  = {3},
  pages   = {285--300},
  year    = {1987},
  doi     = {10.1007/BF02766215}
}

@article{Young1998,
  author  = {Young, Lai-Sang},
  title   = {Statistical properties of dynamical systems with some hyperbolicity},
  journal = {Annals of Mathematics},
  volume  = {147},
  number  = {3},
  pages   = {585--650},
  year    = {1998},
  doi     = {10.2307/120960}
}

@article{DooleyZhang2015,
  author  = {Dooley, Anthony H. and Zhang, Guohua},
  title   = {Local entropy theory of a random dynamical system},
  journal = {Memoirs of the American Mathematical Society},
  volume  = {233},
  number  = {1099},
  year    = {2015}
}

@book{Kifer1986,
  author    = {Kifer, Yuri},
  title     = {Ergodic Theory of Random Transformations},
  series    = {Progress in Probability and Statistics},
  volume    = {10},
  publisher = {Birkh{\"a}user},
  address   = {Boston},
  year      = {1986},
  doi       = {10.1007/978-1-4684-9175-3}
}

@book{Rudin1987,
  author    = {Rudin, Walter},
  title     = {Real and Complex Analysis},
  edition   = {Third},
  publisher = {McGraw-Hill},
  address   = {New York},
  year      = {1987}
}

@article{DownarowiczHuczek2012,
  author  = {Downarowicz, Tomasz and Huczek, Dawid},
  title   = {Faithful zero-dimensional principal extensions},
  journal = {Studia Mathematica},
  volume  = {212},
  number  = {1},
  pages   = {1--19},
  year    = {2012},
  doi     = {10.4064/sm212-1-1}
}

@article{EfronStein1981,
  author  = {Efron, Bradley and Stein, Charles},
  title   = {The jackknife estimate of variance},
  journal = {The Annals of Statistics},
  volume  = {9},
  number  = {3},
  pages   = {586--596},
  year    = {1981},
  doi     = {10.1214/aos/1176345462}
}

@article{Hoeffding1948,
  author  = {Hoeffding, Wassily},
  title   = {A class of statistics with asymptotically normal distribution},
  journal = {Annals of Mathematical Statistics},
  volume  = {19},
  pages   = {293--325},
  year    = {1948},
  doi     = {10.1214/aoms/1177730196}
}

@book{ODonnell2014,
  author    = {O'Donnell, Ryan},
  title     = {Analysis of Boolean Functions},
  publisher = {Cambridge University Press},
  address   = {Cambridge},
  year      = {2014},
  doi       = {10.1017/CBO9781139814782}
}

@book{Rockafellar1970,
  author    = {Rockafellar, R. Tyrrell},
  title     = {Convex Analysis},
  series    = {Princeton Mathematical Series},
  volume    = {28},
  publisher = {Princeton University Press},
  address   = {Princeton},
  year      = {1970}
}

@article{Kolmogorov1958,
  author  = {Kolmogorov, Andrei N.},
  title   = {A new metric invariant of transient dynamical systems and automorphisms in {L}ebesgue spaces},
  journal = {Doklady Akademii Nauk SSSR},
  volume  = {119},
  pages   = {861--864},
  year    = {1958},
  note    = {In Russian}
}

@article{Sinai1959,
  author  = {Sinai, Yakov G.},
  title   = {On the notion of entropy of a dynamical system},
  journal = {Doklady Akademii Nauk SSSR},
  volume  = {124},
  pages   = {768--771},
  year    = {1959},
  note    = {In Russian}
}

@incollection{McDiarmid1989,
  author    = {McDiarmid, Colin},
  title     = {On the method of bounded differences},
  booktitle = {Surveys in Combinatorics, 1989},
  series    = {London Mathematical Society Lecture Note Series},
  volume    = {141},
  pages     = {148--188},
  publisher = {Cambridge University Press},
  address   = {Cambridge},
  year      = {1989}
}

@book{BoucheronLugosiMassart2013,
  author    = {Boucheron, St\'ephane and Lugosi, G\'abor and Massart, Pascal},
  title     = {Concentration Inequalities: A Nonasymptotic Theory of Independence},
  publisher = {Oxford University Press},
  address   = {Oxford},
  year      = {2013}
}

@misc{BassolsCastro2025,
  author        = {Bassols Cornudella, Bernat and Castro, Matheus M.},
  title         = {Conditioned stochastic stability of equilibrium states on uniformly hyperbolic sets},
  year          = {2025},
  eprint        = {2506.15503},
  archiveprefix = {arXiv},
  url           = {https://doi.org/10.48550/arXiv.2506.15503},
  note          = {arXiv preprint}
}

@article{PesinSentiZhang2016,
  author  = {Pesin, Yakov and Senti, Samuel and Zhang, Ke},
  title   = {Thermodynamics of towers of hyperbolic type},
  journal = {Trans. Amer. Math. Soc.},
  volume  = {368},
  number  = {12},
  pages   = {8519--8552},
  year    = {2016},
  doi     = {10.1090/tran/6599}
}

@article{ClimenhagaPesinZelerowicz2019,
  author  = {Climenhaga, Vaughn and Pesin, Yakov and Zelerowicz, Agnieszka},
  title   = {Equilibrium states in dynamical systems via geometric measure theory},
  journal = {Bull. Amer. Math. Soc. (N.S.)},
  volume  = {56},
  number  = {4},
  pages   = {569--610},
  year    = {2019},
  doi     = {10.1090/bull/1659}
}

@article{ClimenhagaPesinZelerowicz2020,
  author  = {Climenhaga, Vaughn and Pesin, Yakov and Zelerowicz, Agnieszka},
  title   = {Equilibrium measures for some partially hyperbolic systems},
  journal = {J. Mod. Dyn.},
  volume  = {16},
  pages   = {155--205},
  year    = {2020},
  doi     = {10.3934/jmd.2020006}
}
 
\end{document}